\documentclass[a4paper,12pt]{smfartsM}
\usepackage[paperheight=5in,paperwidth=7in]{geometry} 
\title
[Finite combinatorics implicit in    the basic definitions of topology]
{Finite combinatorics implicit in \\ the basic definitions of topology}
\thanks{Misha Gavrilovich. University of Haifa, Israel. {\tt{mi\!\!\!ishap\!\!\!p@sd\!\!\!df.org}}. 
\\ The author was supported by ISF grant 290/19.}
\date{2024}
\usepackage{cfr-lm}
\usepackage{MnSymbol}
\usepackage{ae}

\DeclareMathAlphabet{\mathbb}{U}{msb}{m}{n}

\DeclareMathSymbol{0}{\mathalpha}{letters}{`0}
\DeclareMathSymbol{1}{\mathalpha}{letters}{`1}
\DeclareMathSymbol{2}{\mathalpha}{letters}{`2}
\DeclareMathSymbol{3}{\mathalpha}{letters}{`3}
\DeclareMathSymbol{4}{\mathalpha}{letters}{`4}
\DeclareMathSymbol{5}{\mathalpha}{letters}{`5}
\DeclareMathSymbol{6}{\mathalpha}{letters}{`6}
\DeclareMathSymbol{7}{\mathalpha}{letters}{`7}
\DeclareMathSymbol{8}{\mathalpha}{letters}{`8}
\DeclareMathSymbol{9}{\mathalpha}{letters}{`9}

\usepackage{smfthm}
\usepackage{vmargin}
\usepackage{enumitem}
\usepackage{caption}
\usepackage{pdfpages}
\usepackage[matrix, arrow,all,cmtip,color]{xy}
\usepackage{url}
\usepackage{epigraph}

\def\includegraphics[#1]#2{    \pdfximage width  \linewidth {#2} \pdfrefximage\pdflastximage}
\newbox\TestBox
\def\Remove #1 {\setbox\TestBox=\hbox{#1}%
        \leavevmode\rlap{\vrule height 2.5pt depth-1.75pt width\wd\TestBox}%
	        \box\TestBox\ }

\def\includegraphicss[#1]#2{ 
	\pdfximage width  #1\linewidth {#2} 
	\pdfrefximage\pdflastximage 
}

\def\Imm{\operatorname{Im}}

\makeatletter
\newcommand{\colim@}[2]{%
	  \vtop{\m@th\ialign{##\cr
	      \hfil$#1\operator@font colim$\hfil\cr
	          \noalign{\nointerlineskip\kern1.5\ex@}#2\cr
		      \noalign{\nointerlineskip\kern-\ex@}\cr}}%
		      }
		      \newcommand{\colim}{%
			        \mathop{\mathpalette\colim@{\rightarrowfill@\scriptscriptstyle}}\nmlimits@
				}
	
\newcommand{\mylim@}[2]{%
	  \vtop{\m@th\ialign{##\cr
	      \hfil$#1\operator@font lim$\hfil\cr
	          \noalign{\nointerlineskip\kern1.5\ex@}#2\cr
		      \noalign{\nointerlineskip\kern-\ex@}\cr}}%
		      }
\newcommand{\mylim}{%
			        \mathop{\mathpalette\mylim@{\leftarrowfill@\scriptscriptstyle}}\nmlimits@
				}
\makeatother 	

\def\inv{^{-1}}

\def\bblue{\em }
\def\rred{\em }

\def\Topp{\operatorname{Top}}

\def\bofff#1{_\bullet(#1^\leqslant)}\def\boffffbr#1{_\bullet((#1)^\leqslant)}\def\boffff#1{_\bullet(#1^\leqslant)}\def\bofffOne{\boffff{1}}\def\bofffTwo{\boffff{2}}\def\offfOne{_\bullet(1^\leqslant)}
\def\nleqslant{n^\leqslant}\def\Oleqslant{1^\leqslant}\def\Xton{|X|^{n}}\def\XXton{X^{n}}

   \def\rrt#1#2#3#4#5#6{\xymatrix{ {#1} \ar[r]|{} \ar@{->}[d]|{#2} & {#4} \ar[d]|{#5} \\ {#3}  \ar[r] \ar@{-->}[ur]^{}& {#6} }}

\newcommand{\bi}{\begin{itemize}}
\newcommand{\ei}{\end{itemize}}

\def\bee{\begin{enumerate}[label={(\arabic*)},ref={(\arabic*)}]} 
\newcommand{\eee}{\end{enumerate}}

\def\bii{\begin{itemize}[label={(\arabic*)},ref={(\arabic*)}]} 
\newcommand{\eii}{\end{itemize}}

\def\biii{\begin{itemize}}
\newcommand{\eiii}{\end{itemize}}

\newcommand{\bd}{\begin{itemize}\item}
\newcommand{\ed}{\end{itemize}}

\def\bbbR{\Bbb R}
\def\barR{\overline{\Bbb R}}
\def\barRR{[-1,1]} 
\def\ptt{\{\bullet\}}

\def\lra{\longrightarrow}
\def\rtt{\,\rightthreetimes\,}
\def\xra{\xrightarrow}

\def\diag{\mathrm{diag}}

\def\sing{\operatorname{sing}}

\def\MtooL{\left(\MM\atop{\downarrow\atop\vvLambda}\right)}

\def\MtooLl{\left(\MM\atop{\downarrow\atop\vvLambda}\right)^{l}}

\def\MtooLlr{\left(\MM\atop{\downarrow\atop\vvLambda}\right)^{lr}}
\def\MtooLVtoolr{\left({\MM\atop{\downarrow\atop\displaystyle\vvLambda}}\,\, {\vvLambda\atop{\downarrow\atop\bullet}}\right)^{lr}}

\def\MtooLVtoolrScriptsize{\left({\MM\atop{\downarrow\atop\vvLambda}}\,\, {\vvLambda\atop{\downarrow\atop\bullet}}\right)^{lr}}

\usepackage{hyperref}
\usepackage{endnotes}

\def\rrt#1#2#3#4#5#6#7{\xymatrix{ {#1} \ar[r]^{} \ar@{->}[d]_{#2} & {#4} \ar[d]^{#5} \\ {#3}  \ar[r] \ar@{-->}[ur]^{#7}& {#6} }}

\def\rtt{\rightthreetimes}

\def\lra{\longrightarrow}
\def\lr{{\rtt lr}}
\def\lrl{{\rtt l}}

\def\rlr{{\rtt r}}

\def\MtoL{\vvM \to \vvLambda}
\def\MtoLLtoO{\vvM \to \vvLambda,\,\vvLambda\to\ptt}

\def\sSets{\operatorname{sSets}}
\def\Fdiag{X_\bullet^{\mathfrak F\text{-\rm{diag}}}}

\def\Udiag{X^{\mathfrak U\textrm{-diag}}_\bullet} 
\def\UUU{\mathfrak U}
\def\FFF{\mathfrak F}

\def\sk{\operatorname{sk_0}}
\def\cosk{\operatorname{cosk_0}}
\def\FFF{\mathfrak F}

\def\id{\operatorname{id}}

\def\ra{\searrow}

\def\llrra{\leftrightarrow}

\def\xra{\xrightarrow}

\def\lt{<}
\def\gt{>}

\def\Dop{\Delta^{\operatorname{op}}}
\def\Ords{\Dop}
\def\Sets{\text{Sets}}
\def\Hom{\text{Hom}}

\def\PPhi{\ensuremath{\Filt}}

\newtheorem{theorem}{Theorem}
\newtheorem{conjecture}{Conjecture}
\newtheorem{question}{Question}
\newtheorem{exercise}{Exercise}
\newtheorem{definition}{Definition}
\newtheorem{remark}{RemarkT}
\newtheorem{remarkt}{Remark} 

\usepackage[ethiop,english]{babel}
\usepackage{ethiop}

 \newcommand{\ethi}{\selectlanguage{ethiop}}
\def\Filt{{\ethi\ethmath{wA}}}
\def\sFilt{s{{\ethi\ethmath{wA}}}}

\def\lt{<}
\def\gt{>}

\def\Sets{\text{Sets}}
\def\Hom{\text{Hom}}
\def\dist{\operatorname{dist}}

\def\sPhi{\ensuremath{\sFilt}}

\def\PPhi{\ensuremath{\Filt}}
\def\skiip#1{}

\def\smallblackbox{{\Huge\centerdot}}
\def\bblacksquareR{\resizebox{3.39pt}{!}{$\blacksquare$}}
\newcommand{\bblacksquare}{\mathpalette\raisebblacksquare\relax}
\newcommand{\raisebblacksquare}[2]{\raisebox{2\depth}{$#1\bblacksquareR$}}

\def\mathcalC{\left({\emptyset\atop{\displaystyle\Downarrow\atop \displaystyle \bullet }}\right)}
\def\mathcalG{\left({\bullet\phantom{,\bullet}}\atop{{\!\!\!\Downarrow}\atop{\bullet,\bullet}}\right)}
 \def\closedsubspace{{{\raisebox{2pt}{$\star_x\leftrightarrow\star_y\leftrightarrow\star_z$} \searrow \raisebox{-2pt}{
$\smallblackbox_c$} }\atop{\displaystyle\Downarrow}}
 \atop{\displaystyle  \star_{x=y}\leftrightarrow \star_{z=c} }}

\def\mathcalDense{\protect{\hphantom{\bullet_c\searrow}{\!\smallblackbox \atop  \displaystyle\Downarrow }\atop { 
\raisebox{3pt}{$\bullet$}\searrow\raisebox{-3pt}{$\smallblackbox$}  }}} 
\def\mathcalDensebrackets{\left(\mathcalDense\right)}

\def\nondenseimagebrackets{\mathcalDensebrackets}
\def\cl{\operatorname{cl}}

\def\pr{\operatorname{pr}}

\def\bqqq{\begin{quote}}
\def\eqqq{\end{quote}}

\newcommand{\vvLambda}{{\mathpalette\makevvLambda\relax}}
\newcommand{\makevvLambda}[2]{%
	  \raisebox{\depth}{\scalebox{1}[-1]{$\mathsurround=0pt#1\mathsf{V}$}}%
	  }

\newcommand{\vvM}{{\mathpalette\makevvM\relax}}
\newcommand{\makevvM}[2]{%
	  \raisebox{\depth}{\scalebox{1}[-1]{$\mathsurround=0pt#1\mathsf{W}$}}%
	  }

\def\VV{\mathsf{V}}
\def\MM{{\vvM}}
\def\MMtoLL{\vvM\lra \vvLambda}

\def\Oc{\ensuremath{\left\{ 
	{\ensuremath{\bullet_o}}\rotatebox{-12}{\ensuremath{\to}}  \raisebox{-2pt}{\ensuremath{\bblacksquare_c}}\right\}
}}

\def\Outoc{\ensuremath{\left\{ 
	{\ensuremath{\bullet_u}}\rotatebox{-12}{\ensuremath{\to}}  \raisebox{-2pt}{\ensuremath{\bblacksquare_c}}\right\}
}}

\def\oc{ \ensuremath{\left\{{\ensuremath{\bullet_o}}\leftrightarrow {\ensuremath{\bblacksquare_{c}}}\right\} 
}}

\def\ouc{ \ensuremath{\left\{{\ensuremath{\bullet_u}}\leftrightarrow {\ensuremath{\bblacksquare_{c}}}\right\} 
}}

\def\OcTooc{\ensuremath{\left\{ 
	{\ensuremath{\bullet_o}}\rotatebox{-12}{\ensuremath{\to}}  \raisebox{-2pt}{\ensuremath{\bblacksquare_c}}\right\}
\lra
\left\{{\ensuremath{\bullet_o}}\leftrightarrow {\ensuremath{\bblacksquare_{c}}}\right\} 
}}

\def\OutocTooc{\ensuremath{\left\{ 
	{\ensuremath{\bullet_u}}\rotatebox{-12}{\ensuremath{\to}}  \raisebox{-2pt}{\ensuremath{\bblacksquare_c}}\right\}
\lra
\left\{{\ensuremath{\bullet_u}}\leftrightarrow {\ensuremath{\bblacksquare_{c}}}\right\} 
}}

\begin{document}
\selectlanguage{english} \catcode`\_=8\catcode`\^=7 \catcode`\_=8

\begin{abstract} 
We explain 
how to see 
finite combinatorics of preorders 
implicit in the {\em text} 
of basic topological definitions or arguments in (Bourbaki, General topology, Ch.I),  
and define a concise combinatorial notation such that complete definitions of 
connectedness, compactness, contractibility, having a generic point, subspace, closed subspace, 
fit into $2$ or $4$ bytes. This notation is homotopy theoretic in nature, 
and is based on the following observation:

A number of basic  properties of continuous maps and topological spaces 
are defined using 
a single category-theoretic operation, 
taking left or right orthogonal complement with respect to the Quillen lifting property, 
	repeatedly applied to a simple 
	example illustrating the definition or its failure.
	Moreover, for most of these definitions this example 
	can be chosen to be  a map of finite topological spaces (=preorders)  
of size at most $5$.
This includes the properties of a space being connected, compact, contractible, discrete,
having a generic point, 
and a map having dense image, being 
the inclusion of an (open or closed) subspace, or of a component into a disjoint union,
and others. 

Our reformulations illustrate the generative power of the lifting property as a means of defining
basic 
mathematical properties starting from their simplest or typical example.
The exposition is accessible to a student. 
\end{abstract}
\maketitle
\begin{flushright}\tiny
\begin{minipage}[t]{0.84\textwidth}
{\tiny 
it is particularly difficult for an
inexperienced mathematician to orient himself among the multitude of possible directions,
to distinguish what is interesting from what is not. For example, 
it is useful 
for them 
{...}
to hear certain truisms — such as: that general topology should be
viewed as \\ an indispensable and precise language and not as a science that calls for further
research,
}
\\ {\tiny ------ Alexandre Grothendieck.  \href{https://webusers.imj-prg.fr/~leila.schneps/grothendieckcircle/GrothVietnam.pdf}{Mathematical Life in the Democratic Republic of Vietnam}. $1967$.}
\end{minipage}  
\end{flushright}\vskip-15pt\phantom.

\section{Introduction}

We explain 
how to see 
finite combinatorics of preorders 
implicit in the {\em text} 
of basic topological definitions or arguments in (Bourbaki, General topology, Ch.I),  
and define a concise combinatorial notation such that complete definitions of 
connectedness, compactness, contractibility, having a generic point, subspace, closed subspace, 
fit into $2$ or $4$ bytes. This notation is based on the following observation:

\begin{itemize}\item 
A number of basic definitions in topology 
can all be formulated by iteratively using several times 
the same commutative diagram 
starting from maps of finite topological spaces
of size say $\leqslant 5$; 
this includes
the notions of a topological space being 
{\em connected, compact}, {\em contractible}, and 
having {\em a generic point}. 
\end{itemize}
This commutative diagram defines  a binary relation on morphisms,
and is used to take the left or right orthogonal complement with respect to that binary relation.
The binary relation,  
known as the {\em lifting property} or {\em weak orthogonality of morphisms},
is used in an axiomatic approach to homotopy theory by Quillen 
to define  properties of morphisms starting from an explicit list of examples
by taking once or twice the left or right orthogonal complement  
\cite[1\S1, Def.1; 2\S2($\text{sSets}$), Def.1, p.2.2; 3\S3($\Topp$), Lemma 1,2, p.3.2]{Qui}.
In this paper we study what happens if we the orthogonal complement $\geqslant 2$  times
of a finite set of maps of finite topological spaces of small size.



 More precisely, our main result may be summarised 
 as follows:
\begin{enonce*}{Theorem}
A number of basic topological properties can defined
by repeatedly taking, in the category of all topological spaces,
the left or right orthogonal complement 
with respect to the lifting property
of an explicitly listed class of maps of finite topological spaces 
of size $\leqslant 5$. 

Examples of properties include being compact, contractible, connected, totally disconnected, 
having a generic point; 
injective, surjective; quotient, having a section, dense image; 
subspace, closed or open subspace,  and others. 
\end{enonce*}
\begin{proof} See 
 Corollary~\ref{coro:proper}, Theorems~\ref{thm:Thm2} and~\ref{coro:contrCW}, Proposition~\ref{propConnDisconn},
 and also \cite[Fig.1-3]{MR}.
\end{proof} 

In other words, 
a number of basic properties of topological spaces or continuous maps 
can be specified completely by a list of maps of finite topological spaces, 
and a sequence of bits representing taking left or right orthogonal complement. 
This leads 
to a very concise and uniform notation (computer syntax)
able to define 
various basic notions in topology within 
$2$ or $4$ bytes, 
which, we hope, 
has the flavour of universality of 
imprinting and the Hawk/Goose effect 
\cite[p.48,p.67,p.90]{Ergobrain}
and 
is a manifestation of ergo-logic of topology 
\cite[Ergo-Structures/Ergo-systems Conjecture, p.75]{Ergobrain}. 

This notation has an intuition coming from logic rather than topology or homotopy theory, as follows. 
Taking the (left or right) orthogonal complement of a class 
of maps with property $P$ is perhaps the simplest way to define a class of maps 
failing property $P$ in a way useful in a diagram chasing computation. 
Thus we may think of taking the orthogonal complement as a category-theoretic substitute for negation, 
and refer to the (left or right) orthogonal complement of (the class of maps with) property $P$
as {\em (left or right) Quillen negation} of property $P$.
Often, maps of finite spaces used to define a property $P$ are 
 examples of maps illustrating the property $P$ or its failure, or appear implicitly 
in the text of a standard definition or theorem related to property $P$. 
Thus, compactness  
can be defined by taking left-then-right Quillen negation of 4 proper map 
of spaces of size $\leqslant 3$ 
implicit in a characterisation of compact Hausdorff spaces (Corollary~\ref{coro:KtoOlr} and~\ref{coro:proper}).
Contractibility can be defined by taking left-then-right Quillen negation of two trivial Serre fibrations
of spaces of size $\leqslant 3$ and $\leqslant 5$ implicit in the definition of normality (separation axiom $T_4$)
which appears in the Urysohn Theorem \cite[IX\S$4.2$, Theorem~$2$]{Bourbaki}, cf. Theorem~\ref{thm:Thm2}
on extending maps to $\bbbR$, a basic example of a contractible space (Theorem~\ref{coro:contrCW}).
In fact, both notions can be defined using any example of a sufficiently complicated map of finite topological 
classes which is proper, resp.~a trivial Serre fibration (Problem~\ref{prob_AI}).

In this terminology, our main result can be rephrased as follows. 

\begin{enonce*}{Theorem$'$}
Often, a basic topological property can be defined 
by taking several times the left or right  Quillen negation 
of an explicit list of maps of finite spaces 
illustrating the property or its failure, 
or appearing implicitly in the text of a standard definition of the property.

Often it is enough to consider only maps of spaces of size $\leqslant 5$. 
\end{enonce*}

A calculation in \cite{MR} classifies the properties obtained in this way
by repeatedly taking the left or right orthogonal complement 
starting from the simplest possible map, $\emptyset\lra \ptt$: 
there are precisely $21$ of them, and  
$10$ of  them are defined 
in a typical introductory course of topology:
surjective, injective, subspace, closed subspace, quotient, induced topology, 
connected,
discrete, having a generic point, having a section, disjoint union.

In fact, iterated Quillen negation can also be used to define  the notions 
of a finite group being nilpotent, soluble, of odd or even order, 
a module being projective or injective, 
a map being injective or surjective \cite{DMG, LP1,LP2}. 


Finding our reformulations requires very little understanding of topology,
and can be done by either of the two rather simplistic and naive methods;
\begin{itemize}
\item (text$\rightsquigarrow$arrows) 
extract maps 
implicitly mentioned in the {\em text} of a standard topological definition, 
and assemble a commutative diagram out of these maps; 
\item (arrows$\rightsquigarrow$defs) 
pick one or two maps of finite topological spaces and then
take a few times its left or right Quillen negation (i.e. 
left or right orthogonal complement  
with respect to the lifting property). 
\end{itemize}
Both methods are feasible in practice because  sizes of 
finite spaces involved in typical definitions 
are so small, and and there are so few of them. 

Our approach raises a number of questions; we state a few of them in \S\ref{sec_conjectures}.

Evidently, our reformulations  suggest how to view  some standard argu
ments in point-set topology as computations with category-theoretic
(not always) commutative diagrams of preorders. 
In this approach, the
lifting property is viewed as a rule to add a new arrow, a computational
recipe to modify commutative diagrams.
We dicuss this 
in \S\ref{sec_AI}. 
In \S\ref{sec_model_str} we ask whether there is a model structure on 
the category of topological spaces defined entirely in terms of maps of finite topological spaces,
and whether our expression for contractibility defines the class of trivial fibrations. 

\subsection{How to discover our reformulations: Example}
We now explain our methods using as an example the notions of 
connected and injective. 
First we show how to analyse the text of standard definitions of connectedness
and injectivity, and arrive at our diagram-chasing reformulations in terms of
maps of finite topological spaces. As it happens, in both cases we arrive at the 
the same map of 
finite topological spaces $\{\bullet,\bullet\}\lra\{\bullet\}$ 
gluing the  points of a discrete space with two points. 
Note that 
this map is both an example of a map failing to injective, and 
involves a space failing to be connected. 

We then 
see that we could have started from this map: 
iteratively applying left or right Quillen negation to  this map 
leads to the notions of connectivity, injectivity, quotient, and others.

\subsubsection{Transcribing the definition of connectivity.} Consider a standard definition of connectedness.
\begin{itemize}
\item[$(*)_{\text{words}}$] A topological space $X$ is said to be {\em connected} 
iff it is not the union of two disjoint non-empty open subsets.
\end{itemize}

First we analyse the text and rephrase it in terms of maps. 
First notice that 
to represent a space $X$ as ``the union of two disjoint open subsets'' 
$X=A\cup B$, $A\cap B=\emptyset$,
is the same as to give an arrow (the characteristic function) 
$\chi:X\lra\{\bullet_A, \bullet_B\}$
where $\chi\inv(\bullet_A):=A$ and $\chi\inv(\bullet_B):=B$, 
and $\{\bullet_A, \bullet_B\}$ denotes the discrete space with two points.
Either $A$ or $B$ is empty iff this arrow factors via 
the singleton $\{\bullet\}$
as $X\lra \{\bullet\}\lra\{\bullet_A,\bullet_B\}$:
the arrow $\{\bullet\}\lra \{\bullet_A,\bullet_B\}$ 
picks whether $A=X$ or $B=X$. 

Now assemble the arrows above into a commutative diagram,
and obtain the following reformulation.

\begin{itemize}
\item[$(*)_{\text{arrows}}$] A topological space $X$ is said to be {\em connected} 
iff 
\begin{equation}\label{DefInjL} \xymatrix@C=2.39in{ X \ar[d]\ar[r]|\forall & \{\bullet_A,\bullet_B\} \\
                     \{\bullet\} \ar@{-->}[ru]|\exists & } 
\end{equation}
\end{itemize}

\begin{rema}
In fact, we would like to think that this method (text$\rightsquigarrow$arrows) 
is the major contribution of this paper: you could teach this method to a student who would then be able to 
rediscover our reformulations herself.  
However, this would a pure speculation: an exposition of the kind we do in this paper
may hide many implicit choices,
and the only way to see whether this method of translation into arrows is indeed a method
is by trying to teach it to someone. 
\end{rema}

\subsubsection{Transcribing the definition of injectivity.}
Let us now analyse the definition of injectivity.

\begin{itemize}
\item[$(**)_{\text{words}}$] 
``A function $f$ from $X$ to $Y$ is injective iff 
no pair $a\neq b\in X $ of different
points is sent to the same point $y=f(a)=f(b)\in Y$ of $Y$.''
\end{itemize}
``A function $f$ from $X$ to $Y$'' is an arrow $X \lra Y$. 
``A pair of points $a,b\in X$''
is a $\{\bullet_a, \bullet_b\}$-valued point, that is an arrow
$\{\bullet_a, \bullet_b\} \lra X$
from a two element set, or rather a discrete space with two points; 
we ignore word ``different'' for now. 
``the same point ... of $Y$''
is an arrow $\{\bullet\} \lra Y$. 
Represent ``sent to'' by an arrow
$\{\bullet_a, \bullet_b\} \lra \{\bullet_{a=b}\}$. 
Note how we added a subscript $a=b$ to indicate 
where  points $\bullet_a$ and $\bullet_b$ map to.

Now see that the arrows 
fit 
into a commutative diagram of the same shape as Eq.~\eqref{DefInjL}:

\begin{itemize}
\item[$(**)_{\text{arrows}}$] 
A function $f$ from $X$ to $Y$ is injective iff 
\begin{equation}\label{DefInjR} 
\xymatrix@C=2.39in{ 
\{\bullet_a, \bullet_b \} \ar[d] \ar[r]|\forall & X \ar[d]|f   \\
\{\bullet_{a=b}\}  \ar[r]|\forall \ar@{-->}[ru]|\exists & Y
}\end{equation}
\end{itemize}
%

\subsubsection{Taking Quillen negation (orthogonal complement) twice} 
The universal property of a quotient map is obtained from 
the map $\{\bullet_a,\bullet_b\}\lra\{\bullet_{a=b}\}$
by taking right-then-left Quillen negation (orthogonal complement):
\begin{itemize}\item
A map $q:A\lra B$ is a quotient map iff for each map $f:X\lra Y$ the following implication holds:	
\begin{equation}
\begin{gathered}\label{DefQuo}
\vcenter{\hbox{
\xymatrix@C=2.39cm{ 
\{\bullet_a, \bullet_b \} \ar[d] \ar[r]|\forall & X \ar[d]|f   \\
\{\bullet_{a=b}\}  \ar[r]|\forall \ar@{-->}[ru]|\exists & Y}
	}}
implies
\vcenter{\hbox{
\xymatrix@C=2.39cm{ 
A \ar[d]|q  \ar[r]|\forall & X \ar[d]|f   \\
B  \ar[r]|\forall \ar@{-->}[ru]|\exists & Y}
	}}
\end{gathered}\end{equation}
\end{itemize}

Taking 
right  Quillen negation (orthogonal complement) twice gives a somewhat less natural property:
\begin{itemize}\item
A map $q:A\lra B$ is surjective and topology on $A$ is induced from $B$ 
 iff for each map $f:X\lra Y$ the following implication holds:	
\begin{equation}
\begin{gathered}\label{DefCoQuo}
\vcenter{\hbox{
\xymatrix@C=2.39cm{ 
\{\bullet_a, \bullet_b \} \ar[d] \ar[r]|\forall & X \ar[d]|f   \\
\{\bullet_{a=b}\}  \ar[r]|\forall \ar@{-->}[ru]|\exists & Y}
	}}
implies
\vcenter{\hbox{
\xymatrix@C=2.39cm{ 
X \ar[d]|f  \ar[r]|\forall & A \ar[d]|q   \\
Y  \ar[r]|\forall \ar@{-->}[ru]|\exists & B}
	}}
\end{gathered}\end{equation}
\end{itemize}
Both  statements above are easy to verify by hand. 
Detailed proofs are given in \cite[Theorem~\ref{clrrl} and ~\ref{clrrr}]{MR}.
Moreover, \cite[Fig.1]{MR} shows that by iteratively applying Quillen negation to injectivity (i.e.
iteratively taking left and right orthogonal complements of the class of injective maps)
one gets precisely $17$ different classes of maps. Most of these classes are quite natural, 
and can be used to define the notions of dense image, subspace, closed subspace, 
separation axioms $T_0$ and $T_1$, and having a generic point.

\subsection{Finite combinatorics implicit in basic definitions of topology}
We start by demonstrating an example of our combinatorial notation for properties (=classes) of continuous maps. 
Then we write up our expressions for compactness and contractibility, thereby demonstrating 
the finite combinatorics implicit in these notions.

\subsubsection{Our notation and iterated Quillen negation} 
The reader should be able to guess our notation if we say that 
the classes of maps defined by \eqref{DefInjR}, \eqref{DefQuo}, \eqref{DefCoQuo} are denoted by
\def\abTopt{\left(\bullet_a,\bullet_b \atop { \displaystyle {\Downarrow \atop \displaystyle \bullet_{a=b} }} \right) }
$$\abTopt^r\,\,\,\,\,\,\,\, \abTopt^{rl}\,\,\,\,\,\,\,\,\abTopt^{rr}$$
Using this notation we are able to say that \cite[Theorem~\ref{clrr}]{MR} says that 
\def\oTopt{\left( \emptyset \atop { \displaystyle {\Downarrow \atop \displaystyle \bullet }} \right)} 
$$\abTopt^r=\oTopt^{lrr},$$
and that \cite[Fig.1-3]{MR} gives a finite graph of all the classes (properties) of continuous maps 
represented by expressions of form 
$\mathcalC^w$ where $w\in\{l,r\}^{<\omega}$ is a word in letters $l,r$.

\subsubsection{Finite combinatorics implicit in the notions of compactness and contractibility}
We now give more complicated examples of combinatorial expressions in our notation
followed by a brief explanation. 
We hope the reader shall be able to guess the meaning of these expressions;
their syntax 
is introduced in detail in \S\ref{not:mintsGE}. 
These expressions appear in 
 Corollary~\ref{coro:proper}, Theorem~\ref{thm:Thm2} and Theorem~\ref{coro:contrCW}, 
 also \cite[Thm.~\ref{crllrrll} and~\ref{crll}]{MR}. 

%

connectedness is $\mathcalC^{rll}=\mathcalG^{lr}$ 

dense image is 
        $\mathcalC^{rllrrll}= \left(\closedsubspace\right)^{ll}=\nondenseimagebrackets^l$, 

compactness is 
	$ \left(
	{ {\{\filledstar_a\leftrightarrow \filledstar_b\}}\atop{\Downarrow
	\atop{ \{\filledstar_{a=b} \} } }}   \,\,\,
	{\left\{{\ensuremath{\bullet_o}}\rotatebox{-12}{\ensuremath{\to}}  \raisebox{-2pt}{\ensuremath{\bblacksquare_c}}\right\}\atop{\Downarrow\atop{\{\bullet_{o=c}\}}}}
	\,\,\,
	{{{\{\bblacksquare_c\}}\atop{\Downarrow
		\atop{\left\{{\ensuremath{\bullet_o}}\rotatebox{-12}{\ensuremath{\to}}  \raisebox{-2pt}{\ensuremath{\bblacksquare_c}}\right\}\  
	 }}}}%
	{{	\{\bblacksquare_a\raisebox{6pt}{}\rotatebox{12}{\ensuremath{\leftarrow}}\raisebox{2pt}{\ensuremath{\bullet_u}\raisebox{4pt}
				{}\rotatebox{-13}{\ensuremath{\rightarrow}}} \bblacksquare_b\} }\atop{\Downarrow
	\atop{ \{\filledstar_{a=u=b} \} } }}   \,\,\,
	\right)^{lr}
\subset
\left( \left( 
	{{{\bullet_o\hphantom{\to\bblacksquare_c}}\atop{\Downarrow
		\atop{ \displaystyle{\ensuremath{\bullet_o}}\rotatebox{-12}{\ensuremath{\to}}  \raisebox{-2pt}{\ensuremath{\bblacksquare_c}}  }}}}%
		\right)_{<5} \right)^{lr}$

contractibility is 
$ \left(
{
{\left\{\underset{\bblacksquare_a}{}{\swarrow} \overset{\raisebox{1pt}{$\bullet_u$}}{}{\searrow} \underset{\bblacksquare_b}{}
\right\}
}
\atop
{\displaystyle \Downarrow \atop
{\left\{ 
\overset{\raisebox{1pt}{$\bullet_{a=x=\!b}$}}{}
\right\}
}
}
}
\phantom{AA}
{
{\left\{\underset{\bblacksquare_a}{}{\swarrow} \overset{\raisebox{1pt}{$\bullet_u$}}{}{\searrow} \underset{\bblacksquare_x}{}{\swarrow}\overset{\raisebox{1pt}{$\bullet_v$}}{}{\searrow}\underset{\bblacksquare_b}{}\right\}
}
\atop
{\displaystyle \Downarrow \atop
{\left\{\underset{\bblacksquare_a}{}{\swarrow}\overset{\raisebox{1pt}{$\bullet_{u=x=v}$}}{}{\searrow} \underset{\bblacksquare_b}{}\right\}
}
}
}
\right)^{lr}
$

In fact, conjecturally the last two expressions denote the class of proper maps and the class of trivial Serre fibrations. 
The upper subscript $l$ or $r$ 
denotes taking left or right orthogonal complement with respect to the lifting property. 
The expression in brackets denotes a list of maps of finite preorders; 
we hope that positioning of the vertices makes the maps  visually apparent 
but we also use subscripts to indicate what maps to what. 
Here we view a finite preorder 
as a topological space such that $x \lra y$ iff $y \in \cl x$.
Equivalently, we may view a preorder as a category where each diagram commutes;
a monotone map of preorders becomes then a functor.
Hence we may say that in our notation we define topological properties
in terms of functors between categories with finitely many objects.

\subsubsection{Ergo-logic of Gromov/AI}
In \cite[\S4.3, Question 4.12]{mintsGE} 
we suggest it appears worthwhile to try to develop a formalism, 
or rather a very
short program (kilobytes) based on such a formalism 
which supports reasoning in
elementary topology in terms of the combinatorial expressions in our notation such as
the ones listed above.
In terms of \cite[II\S25,p.168]{Ergobrain}, the interesting structure 
built within itself by such an ergo-learner program from texts on topology, 
would be based on maps of finite preorders of small size, or, equivalently,
functors between categories with finitely many objects where each diagram commutes.

In fact, our reformulations 
arose in an attempt to understand ideas 
of Misha Gromov \cite{Ergobrain} 
about ergo-logic/ergo-structure/ergo-systems. Oversimplifying,
ergo-logic is a kind of reasoning which helps to understand how to generate proper
concepts, ask interesting questions, and, more generally, produce interesting rather
than useful or correct behaviour. He conjectures there is a related class of mathematical, essentially combinatorial, structures, called ergo-structures or ergo-systems,
and that this concept might eventually help to understand complex biological behaviour including learning and create mathematically interesting models of these
processes.
We hope our observations may eventually help to uncover an essentially combinatorial reasoning behind elementary topology, and thereby suggest an example of
an ergo-structure.

%
%
%
%
\subsection{Structure of the paper}
The exposition is accessible to a student, and, 
 we hope,
can be used as a source of
exercises demonstrating the expressive power of category theory 
as a concise and uniform language for mathematics,
on well-motivated examples from an introductory course of topology.
In \S\ref{sec_teaching_exercises} we explain this in detail.
In \S\ref{sec_synopsis} we give a synopsis of our reformulations of connectedness, compactness, and contractibility.
In \S\ref{sec_conjectures} we list several open questions motivated by our approach. 
In \S\ref{not:mintsGE} we define  our notation for maps of finite topological spaces.
The purpose of \S\S\ref{sec_connected}-\ref{sec_Kompakt} 
is to explain how to find our reformulations by analysing  
the text of \cite{Bourbaki}.

In \S\ref{sec_connected} we introduce category-theoretic background 
while trying to avoid technicalities associated 
with  our notation for finite topological spaces 
defined in \S\ref{not:mintsGE}.
We do so by 
by analysing excerpts from \cite[I\S11.1-5]{Bourbaki} discussing connectedness.

In \S\ref{sec_T4}  we see that Separation Axiom $T_4$(normality) 
implicitly describes a map of finite topological spaces (Eq.~\eqref{T4Odv})
which happens to be a trivial Serre fibration, and 
in fact is a finite model of the barycentric subdivision of the interval.
The same map is also implicit in 
the standard proof of the Urysohn theorem 
on  extending a real-valued map from a closed subset of a normal space  (Theorem~\ref{thm:Thm2}), 
and this leads to a diagram-chasing reformulation of the proof 
leading to a somewhat weaker statement characterising contractibility and, conjecturally,
the class of trivial fibrations (Theorem~\ref{coro:contrCW}).

In \S\ref{sec_Kompakt} we 
analyse the statement of a theorem characterising compactness in terms of extending maps
from certain dense subspaces, and 
see that that the maps of finite spaces 
implicit in its formulation 
are all closed and thereby proper (Eq.\eqref{KtoOlr}, \eqref{KtoOlrSimplified}), 
and, in turn, we use them to define 
in a diagram-chasing manner compactness for Hausdorff spaces,
and, conjecturally, the class of proper maps.
Our reformulation is very close to the definition 
of compactness and proper maps in terms of limits of ultrafilters (Corollary~\ref{coro:proper}).


The exposition
is in the form of 
an analysis of excerpts from the {\em text} of 
\cite[Ch.I]{Bourbaki} 
and aims to be self-contained and accessible
to a first year student who has taken 
an introductory course in  naive set
theory, topology, and who has heard a definition of a category. A more
sophisticated reader may find it more illuminating to recover our for
mulations herself from reading either the abstract, or the abstract and
the observation of the next subsection. 
The displayed formulae 
provide complete formulations
of our theorems to such a reader.
Our approach naturally leads to a more general observation that
in basic point-set topology, a number of arguments are computations
based on symbolic diagram chasing with finite preorders. 

Our 
analysis of \cite[Ch.I]{Bourbaki} 
is similar in spirit to the  analysis of Ibn Sina by \cite{Hodges05} 
and follows the paradigm  of J.Lukasiewicz ``Look for formal systems'' as described in \cite[p.15]{Hodges09}.

\subsection{Use in teaching}\label{sec_teaching_exercises}

We suggest the following is a series of illuminating exercises 
demonstrating the expressive power of category theory as 
a concise and uniform language of mathematics. They are accessible
to a student familiar with basic definitions of topology and
the notion of a commutative diagram in category theory. 

\subsubsection{Warm-up}
First, ask a student to reformulate in terms of maps to/from finite topological spaces 
the notions of 
{\em connectedness, injectivity, induced topology, a map having dense image}.
A student will notice  that all the reformulations 
(\eqref{DefCon}, \eqref{def:Engel:inj}-
\eqref{def:Engel:dense})
result in the same commutative diagram. Point out 
that this diagram appears in a prominent way in an axiomatic approach 
to homotopy theory by \cite{Qui}. 

To demonstrate unity of mathematics, further point out that 
the same diagram can be used to define a few notions in algebra as well: 
for a student familiar with the notions involved  
it is an illuminating exercise to reformulate the notions of
a finite group being of odd order, soluble, perfect, or nilpotent; 
a module being projective or injective, 
a homomorphism being injective, surjective, epi- or mono-, 
and a few others \cite{DMG, LP1, LP2}. 

\subsubsection*{Compactness} 
Ask a student to reformulate in terms of maps 
the definition of compactness in terms of ultrafilters.
Point out 
that \cite[I\S$10.2$,Th.1d)]{Bourbaki} (cf.~Proposition~\ref{propo_Theorem1d})
uses this reformulation to define
the class of proper maps, a property of maps generalizing compactness.

A 1960s theorem \cite[Theorem $3.2.1$]{Engelking} 
gives a necessary and sufficient condition to extend a map to a compact Hausdorff space
from a dense subset. 
Ask a student 
to reformulate this theorem in terms of maps to finite spaces. 
The student then may notice that all the maps of finite spaces
contained in the resulting diagrams 
\eqref{DefCompUltras}-\eqref{quasicompactUltrafilter}  
and in its reformulation Proposition~\ref{thm:Engel},
are closed and therefore proper. 
Relating this to the observation above leads 
to a characterisation 
\eqref{KtoOlr}-\eqref{KtoOlrSimplified} 
of compactness
for Hausdorff spaces
in terms of maps of finite spaces.
Finally, use this to show that  
Stone-Cech compactification is an instance of a universal construction
similar to something in homotopy theory, namely Axiom M2 of Quillen's model categories, cf.~Conjecture~\ref{conjM2}.

\subsubsection*{Contractibility} 
Ask a student to reformulate Separation Axiom $T_4$ (normality) 
in terms of maps to finite topological spaces. 
Then ask a student to consider Urysohn Theorem
\cite[IX\S$4.2$,Theorem $2$]{Bourbaki} saying that it is always possible
to extend a map to the real line from a closed subspace of a normal space,
and recall its proof. 
Point out that the proof constructs by induction finer and finer approximations, 
and this corresponds to dividing the interval $[0,1]$ into smaller and smaller
sub-intervals. Then point out that sub-diving the interval into two parts
is exactly the geometric meaning of the map of finite spaces 
occurring in the reformulation \eqref{T4Odv} of normality, 
cf.~\S\ref{def:MtoL:above}. 
Explain how this observation leads to a diagram-chasing
reformulation of the proof leading to a somewhat weaker statement 
Theorem~\ref{thm:Thm2}, and a characterisation of contractible spaces
in terms of maps of finite spaces.

Finally, a student who have heard of homotopy theory, might find it amusing
that the map in \eqref{T4Odv} is a trivial Serre fibration, and 
that the reformulation of Urysohn theorem \cite[IX\S$4.2$,Theorem $2$]{Bourbaki} 
leads to a conjectural characterisation
of trivial fibrations (cf.~Conjecture~\ref{conj:CW-WC}).

\subsection{A synopsis of the paper\label{sec:synopsis}\label{sec_synopsis}} We now sketch how to express the notions of a space being {connected, contractible, and compact} 
in terms of finite topological spaces and iterated Quillen negation/lifting property. 

A common pattern is as follows. Often a definition of a property of a topological space or map 
describes a finite system 
 of distinguished subsets with certain properties, 
and claims it 
can be extended or otherwise modified in some way. Such a system of subsets
is ``classified'' by a map to a finite topological space whose topology 
reflects the required properties of the subsets.
This shows  one or two  finite topological spaces, and hopefully a map, implicit in the definition, 
or perhaps in the {\em text} of the definition. 
Now pick a morphism related to the definition in this or other way, 
and take several times its orthogonal with respect to the Quillen lifting property. 
A common coincidence is that in this way you can define the class of morphisms having the property 
you started from. A couple of  examples of such  coincidences are sketched in \S\ref{sec:synopsis} below
and later described in detail. In particular, we see that the morphisms implicit in characterisations 
of compactness and contractibility are proper and trivial Serre fibrations, resp.

In this subsection we assume the reader is familiar 
with the lifting property. 
For a class $P$ of maps we denote by  $P^{l}$ and $P^r$ the class of maps
which have the left, resp. right, lifting property with respect to each map in $P$ (see Def.~\ref{def:rtt}). 
For a word $w\in \{l,r\}^{<\omega}$, we say that {\em a map $f$ or a property (=class) $P$ of maps $w$-defines property $Q$} iff $Q=P^w$, 
and that {\em a map $f$ or a property (=class) $P$ of maps $w$-defines a property $Q$ for spaces with property $R$} 
iff for each space $X$ with property $R$ it holds 
 ${ X\atop{\downarrow\atop\ptt} }\in \{f\}^{w} $ iff $X$ has property $Q$.
We also say that {\em a property $Q$ is Quillen-definable in terms of $P$} iff $Q=P^w$ for some word $w\in\{l,r\}^{<\omega}$.

\subsubsection{Connectedness} A space is called {\em connected} iff it is the union of two disjoint open subsets.
To represent a space $X$ as the union of two disjoint open subsets is the same  as to give a map $X\lra \{\bullet,\bullet\}$
to the discrete space with two points \cite[I\S$11.2$,Proposition $5$]{Bourbaki}. 
A verification shows that a space $X$ is connected iff 
the map $X\lra\ptt$ lifts with respect to the map $\{\bullet,\bullet\}\lra \ptt$. 
A metrisable space $Q$ is profinite, or equivalently 
totally disconnected compact Hausdorff,  
 iff $Q\lra \ptt$ lies in 
the left-then-right orthogonal complement of the class consisting of the map $\{\bullet,\bullet\}\lra\ptt$, 
which we denote by $\bigl\{\{\bullet,\bullet\}\lra\{\bullet\}\bigr\}^{lr}$. 

Thus, the map $\{\bullet,\bullet\}\lra \ptt$ $l$-defines connectedness, and $lr$-defines 
being profinite for metrisable spaces, and totally disconnected
for metrisable compact Hausdorff spaces. 
Being totally disconnected for all spaces is $lr$-defined by 
a class consisting of $4$ maps of spaces with $\leqslant 3$ points  implicit in 
 \cite[I\S$11.5$,Proposition $9$]{Bourbaki}, cf.~Proposition~\ref{prop9}, which roughly 
 says that each topological space has a unique totally disconnected quotient 
 such that the quotient map has connected fibres.

\subsubsection{Contractibility} The definition of a {\em normal} space (Axiom $T_4$) implicitly contains a map of finite spaces
with $5$ and $3$ points, 
which we denote by $\MtoL$ (the shapes of letters 
represent the specialisation preorder of the spaces): the space $\vvLambda$ is obtained by contracting two disjoint closed subspaces $A$, $B$, and their complement $X\setminus(A\cup B)$,
in a connected Hausdorff space $X$, and $\vvM$ is obtained by contracting each of $A$, $U\setminus A$, $V\setminus B$, and $X\setminus(U\cup V)$,
for their disjoint open neighbourhoods $U\supset A$ and $V\supset B$. 
This map is a trivial Serre fibration and so is the map $\vvLambda\lra\ptt$ contracting $\vvLambda$ to a point. 
Because the notion of a trivial Serre fibration is defined by a left lifting property, diagram chasing considerations 
imply that so is also any map in $\{\MtoLLtoO\}^{lr}$. 
We conjecture that in some sense it is the class of trivial fibrations, and are able to 
show that a finite CW complex $X$ is contractible iff $X \lra\ptt\in \{\MtoLLtoO\}^{lr}$, as follows.

An analysis of the proof of the Urysohn theorem 
\cite[IX\S$4.2$,Theorem $2$, Corollary]{Bourbaki}
characterising normality in terms of  being able to always extend a real valued map from a closed subspace, 
shows it is a diagram chasing argument giving that $[-1,1]\lra\ptt\in \{\MtoLLtoO\}^{lr}$
and $\Bbb R \lra\ptt\in \{\MtoLLtoO\}^{lr}$. The class $\{\MtoLLtoO\}^{lr}$ is closed under retracts and products,
hence $X \lra\ptt\in \{\MtoLLtoO\}^{lr}$ for any retract of a Cartesian power of $\bbbR$. 
Being such a retract is equivalent to contractibility for a finite CW complex, 
hence for a contractible finite CW complex $X$ it holds $X \lra\ptt\in \{\MtoLLtoO\}^{lr}$;
the other direction follows from the fact $X\lra\ptt$ being a trivial Serre fibration
implies that the finite CW complex $X$ is contractible.

Thus, contractibility for finite CW complexes is $lr$-defined by two maps with $\leqslant 5$ points. 

\subsubsection{Compactness} The four maps implicit in a Smirnov-Vulikh-Taimanov theorem 
\cite[Thm. $3.2.1$,p.$136$]{Engelking}
giving necessary and sufficient conditions to extend 
a map to a compact Hausdorff space from an dense subspace, are all proper. Being proper
is a right lifting property \cite[I\S$10.2$, Theorem $1d$); I\S$6.5$,Example,p.$62$]{Bourbaki}, as is the definition of compactness via ultrafilters,
and this implies that for any class $P$ of proper maps $P^{lr}$ is also a class of proper maps.
If $P$ contains the $4$ proper maps implicit in the Smirnov-Vulikh-Taimanov theorem, 
then for a Hausdorff space $K$ is compact iff $K\lra\ptt\in P^{lr}$. 
In fact, it is enough to pick a single proper map
between spaces with $3$ and $2$ points, and, moreover, almost any complicated enough 
non-surjective closed map of finite spaces would do. 
Moreover, for a maps of finite spaces being proper is equivalent to the left lifting property with respect
to the simplest non-proper map $\{\bullet_o\}\lra \Oc$ sending  a point into the open point of the Sierpinski space.
This allows us to state a conjecture defining the class of proper maps in terms of the simplest counterexample
negated three times.

Thus, compactness for Hausdorff spaces is $lr$-defined by a map of spaces with $\leqslant 3$ points,
and, in fact, by any non-surjective closed map of finite topological spaces
complicated enough in some precise sense 
(Corollaries~\ref{coro:proper},\ref{coro:KtoOlr},\ref{coro:KtoOlrSimplifiedFurther}).

\subsection*{Acknowledgments and historical remarks}

It seems embarrassing to thank anyone for ideas so trivial, and we
do that in the form of historical remarks. 
All the ideas date back to the note \cite{DMG} 
which `[was written] to demonstrate 
 that
some natural definitions are lifting properties relative to the simplest
counterexample, and to suggest a way to ``extract'' these lifting prop-
erties from the text of the usual definitions and proofs'. 
Arguably, all examples in this paper would have been found thorough 
a diligent attempt by an expert to extract lifting properties or finite topological spaces  
from the text of  \cite[I,IX]{Bourbaki}.

Reformulations \cite[(*)${}_{\rtt}$ and (**)${}_{\rtt}$]{DMG} 
of surjectivity and injectivity, as well as connectedness
and (not quite) compactness
in terms of the lifting property, first
appeared in early drafts of a paper \cite{GH}  with Assaf Hasson
 as
trivial and somewhat curious examples needed to illustrate the concept of a lifting property
but were removed during preparation for publication. After 
the reformulation \cite[(**)${}_{\rtt}$]{DMG} of injectivity in terms the lifting property with respect to the simplest example of a non-injective map
came up in a conversation with Misha Gromov the author decided to try to
think seriously about such lifting properties, and in fact gave talks at logic seminars in $2012$ at Lviv
and in $2013$ at Munster and Freiburg, and $2014,2022$ at St.\ Petersburg, and $2023,2024$ 
in Jerusalem, and $2023,2024$ in Haifa.

The ideas of \cite{DMG} here have greatly influ
enced by extensive discussions with Grigori Mints, Martin Bays, and,
later, with Alexander Luzgarev and Vladimir Sosnilo. At an early stage
Xenia Kuznetsova helped to realise an earlier reformulation of com
pactness was inadequate and that labels on arrows are necessary  to
formalise topological arguments, and from then on trying to ``understand compactness''
was a main goal.
Alexandre Borovik suggested to write a note \cite{DMG} for
The De Morgan Gazette explaining the observation that `some of hu
man's ``natural proofs'' are expressions of lifting properties as applied
to ``simplest counterexample''\,'.

The reformulation of compactness  first appeared in \cite[$2.2.3-4$,pp.$14-16$]{mintsGE}.
Tietze extension lemma (also known as Urysohn theorem \cite[IX\S$4.2$,Theorem $2$]{Bourbaki}) 
also first appeared in \cite[$2.4.1$,p.$25$]{mintsGE},
but a realisation that its reformulation also defines contractibility came only 
after a  question by Sergei V. Ivanov
about it at the Alexandrov geometry seminar at St.Petersburg in $2022$. 
Tyrone Cutler clarified for us the relationship between contractibility, ANR, 
and retracts of $\Bbb R^\kappa$ for cardinal $\kappa$ \cite{Douwen-Pol}. 
Misha Levine brought a reformulation of the Brouwer theorem to our attention.
We thank Anna Erschler for comments on a preliminary version of this paper, and helpful discussions.

\vspace{5pt}
\begin{flushright}\tiny
\tiny
\begin{minipage}[t]{0.73\textwidth}
\begin{minipage}[t]{0.958\textwidth}
%
%
%
%
%
%
 if a man bred 
 to the seafaring life, 
  % ...  
 and accustomed to think and talk 
 only of matters relating to navigation, 
... 
should take it into his 
 head to philosophize concerning the faculties of the mind, 
 it cannot be doubted, but he would draw his notions from 
 the fabric of his ship, and would find in the mind, sail, 
 masts, rudder, and compass. 
 
Sensible objects of one kind or other, do no less occupy 
and engross the rest of mankind, than things relating to 
navigation, the seafaring man. For a considerable part 
of life, we can think of nothing but the objects of sense 
and to attend to objects of another nature, so as to form 
clear and distinct notions of them, is no easy matter, even 
after we come to years of reflection. 
%
\\ \tiny
------ Thomas Reid. An Inquiry into the Human Mind on the Principles of Common Sense. $1764$. 
\end{minipage}  
\end{minipage}
\end{flushright}
\phantom.\vspace{-15pt}\phantom.
\section{Connectedness}\label{sec_connected} 
%
%
%
%
\subsection{Connectedness in terms of 
commutative diagrams and the simplest counterexample} 

We see how to define connectedness in terms of 
the simplest example of a non-connected space, namely the discrete space with two points,
and some diagram chasing. This example has the advantage that we need not
use our notation for maps of finite spaces. 

\subsubsection{Connected.} We quote \cite[Definition 1, \S$11.1$]{Bourbaki}:
\noindent\newline\includegraphicss[.99]{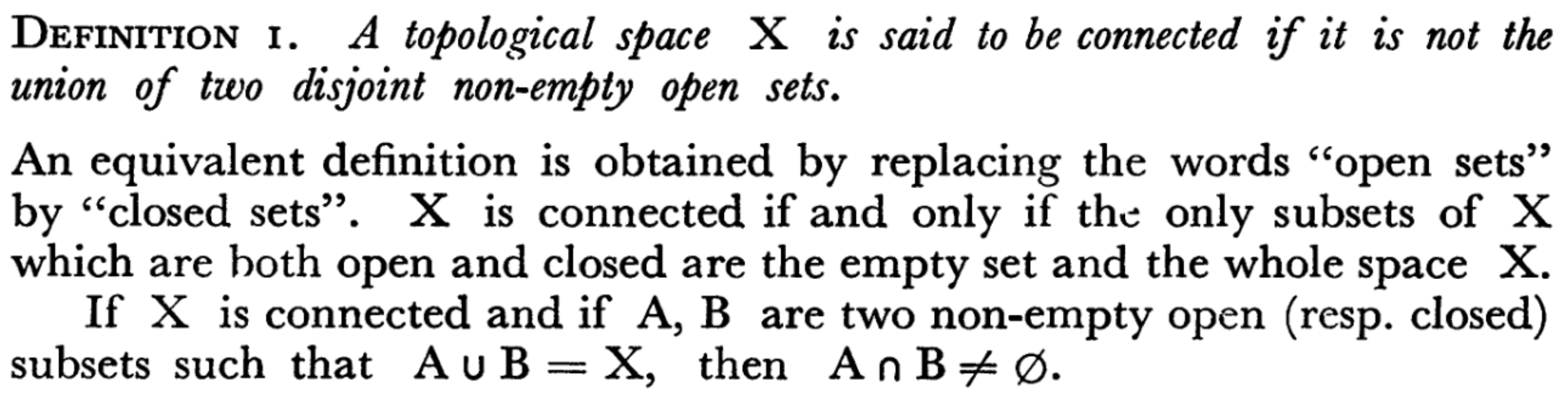} 
\newline The simplest example of a non-connected space 
is a discrete space $\{a,b\}$, with two elements, as witnessed by two disjoint open subsets: $\{a\},\{b\}\subset\{a,b\}$. 

This archetypal counterexample to connectivity is implicit in the definition above:
To give two 
disjoint open subsets $A$, $B$ of $X$ whose union is $X$ 
is the same as to give a 
continuous mapping $\chi:X\lra\{a,b\}$ of $X$ 
into a discrete space of two ele
ments $\{a, b \}$ defined by $\chi (A) = \{a\}$ and $\chi (B) = \{b\}$. 
We shall call $\chi_{A,B}:=\chi:X\lra\{a,b\}$ 
			the {\em characteristic (or indicator) function of 
			the system of two disjoint open subsets $A,B$ of $X$}. 

To say that both subsets are non-empty is to say that this map is surjective.
Bourbaki \cite[I\S$11.2$(the image of a connected set under a continuous mapping)]{Bourbaki} 
phrases this reformulation as follows:
\newline\includegraphics[test]{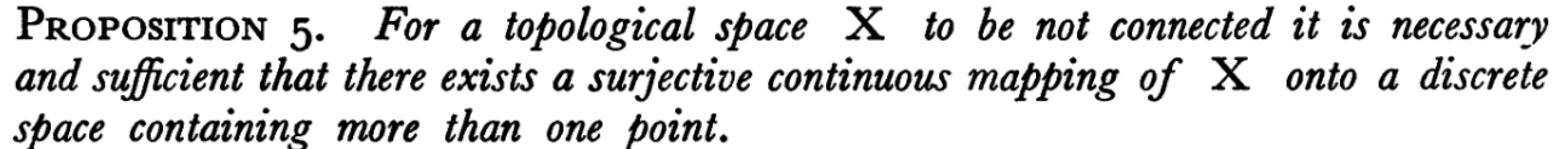}
Now we transcribe this Proposition from Bourbaki 
in terms of diagrams and maps of finite topological spaces as follows.
\begin{enonce}{Proposition}[{\cite[I\S$11.2$,Proposition $5$]{Bourbaki}}] 
For a topological space $X$ to be 
	connected it is necessary
and sufficient that 
there does not exists a surjective continuous mapping of $X$ onto a discrete
space with two points. 
In other words, $X$ is connected iff the following diagram holds:
\begin{equation}
\begin{gathered}\label{DefCon}
	\xymatrix{ {X} \ar[r]|{\forall} \ar@{->}[d]_{} & {\{a,b\}} \ar[d]^{} \\ {\{\bullet\}}  \ar[r] \ar@{-->}[ur]|{\exists}& {\{\bullet\}} } 
\end{gathered}\end{equation}
\end{enonce}	
	\begin{proof}	The top horizontal arrow 
		$\chi:X\lra \{a,b\}$ defines 
		a decomposition $X=A\cup B$ into a union of two disjoint open 
		subsets
		$A:=\chi\inv(a)$ and $B:=\chi\inv(b)$. 
		Commutativity of the upper triangle means precisely that either $A=X$ 
	(if the diagonal arrow sends $\bullet$ to $a$) 
or	$B=X$ (if it sends $\bullet$ to $b$).

		The square and the lower triangle both always commute, 
		as the singleton is the final object of $\Topp$. 
	\end{proof}

In (\ref{DefCon}) the lower triangle can be ignored, but will play a role in similar examples later.


%

	\subsubsection{Totally disconnected}
Recall that a space $Q$ is said to be {\em totally disconnected} iff
	each connected subset of $X$ consists of a single point, cf.~\cite[I\S$11.5$(Components)]{Bourbaki}. 

We can define this property by using the same commutative diagram twice:

			\begin{enonce}{Proposition}\label{propDefDisCon}
	A space $Q$ is {\em totally disconnected} iff 
for each space $X$ the following implication holds: 
		\begin{equation}
\begin{gathered}\label{DefDisCon}
\vcenter{\hbox{
	\xymatrix{ {X} \ar[r]|{\forall} \ar@{->}[d]_{} & {\{a,b\}} \ar[d]^{} \\ {\{\bullet\}}  \ar[r] \ar@{-->}[ur]|{\exists}& {\{\bullet\}} }}}
implies
\vcenter{\hbox{
	\xymatrix{ {X} \ar[r]|{\forall} \ar@{->}[d]_{} & {Q} \ar[d]^{} \\ {\{\bullet\}}  \ar[r] \ar@{-->}[ur]|{\exists}& {\{\bullet\}} }}}
\end{gathered}\end{equation}
\end{enonce}
	\begin{proof} The diagram in the conclusion says that each map $X\lra Q$ sends all of $X$ to a single point,
		namely the image of the diagonal arrow $\{\bullet\}\lra Q$. As we saw above, the diagram in the hypothesis
		says that $X$ is connected. Thus the implication says that for each map $f:X\lra Q$ from a connected space $X$ 
		its image $f(X)\subset Q$ in $Q$ is a single point. 
		Now recall that the image of a connected set is connected \cite[I\S$11.2$,Prop.4]{Bourbaki}. 
		Hence, if $Q$ is totally disconnected, the conclusion holds. 
		Conversely, if $Q$ is not totally disconnected, take $X$ to be a connected subset of $Q$ which is is not a singleton,
		then the inclusion map $X\lra Q$ witnesses the failure of the lifting property required by the conclusion. 
	\end{proof}

It is easy to see that in the same way one may 
reformulate and prove Propositions $1,3-8$ 
in \cite[I\S$11$(Connectedness)]{Bourbaki} 
by diagram chasing calculations.



\subsection{The Quillen lifting property: category theory background} 

The diagram used in both (\ref{DefCon}) and (\ref{DefDisCon}) 
has a name, and in fact 
appears in a prominent way in the theory of model categories, an axiomatic framework for homotopy theory introduced by Daniel Quillen \cite[I\S1,Def.1; 2\S2($\text{sSets}$), Def.1, p.$2.2$; 3\S3($\Topp$), Lemma 1,2, p.$3.2$]{Qui}.
Perhaps it is {\em the simplest way to define a class of morphisms without a given property in a manner useful in a diagram chasing computation},
and a remarkable number of {\em basic notions from a first year course in algebra or topology may be defined by iteratively using this diagram 
starting from an explicit list of counterexamples \cite{LP1,LP2}}.
For this reason it is convenient to talk about the left or right 
{\em  Quillen negation} of a property of morphisms.   

\subsubsection{The definition of the Quillen lifting property} 
\begin{defi}\label{def:rtt} We write $f \rtt g$ and 
	say that a morphism $f:A\lra B $ {\em has the right lifting property} or {\em weakly orthogonal} to or {\em antagonises}  a morphism $g:X\lra Y$ 
	iff
	for each $t:A\to X$  and $b:B\to Y$ such that $ g \circ t = b \circ f$ 
	there exists $d : B \lra  X $
	such that $ h \circ f = t$ 
	and $g \circ h = b $, i.e.~the following diagram holds
	\begin{equation}
\begin{gathered}\label{DefLift}
\vcenter{\hbox{
	\xymatrix@C=+2.39cm{ {A} \ar[r]|{\forall\,t} \ar@{->}[d]_{f} & {X} \ar[d]^{g} \\ {B}  \ar[r]|{\forall\,b} \ar@{-->}[ur]|{\exists\,h}& {Y} }}}
\end{gathered}\end{equation}
\end{defi}

For example, in any category $f \,\rtt\, f$ means that $f$ is an isomorphism,\footnote{Indeed, 
if one takes both horizontal arrows to be identity, then the diagram becomes
the standard definition of an isomorphism. Conversely, take $h:=t\circ f\inv$ if $f\inv$ is well-defined.}
and in $\Sets$ or $\Topp$ $\emptyset\lra\ptt \,\rtt\, g$ means that $g$ is a surjection,\footnote{Indeed,
both the square and the upper triangle always commute because $\emptyset$ is the initial object.
Hence, the diagram says that each map $g:\ptt\lra Y$ factors via $X$, i.e.~for each $y:=g(\bullet)\in Y$ 
there is $x\in X$ such that $g(x)=y$, which is precisely the standard definition of surjectivity.}
and $\{\bullet,\bullet\}\lra\ptt \,\rtt\, g$ means that $g$ is injective.\footnote{Indeed, 
to give a map $t:\{\bullet,\bullet\}\lra X $ is the same as to give two points $x,y\in X$.
The map fits into a commutative square iff $g(x)=g(y)$. There is a lifting iff $x=y$. Thus,
the lifting property says that  $g(x)=g(y)$ implies $x=y$, which is precisely the standard definition of injectivity.}
In $\Sets$ $f \,\rtt\, \{\bullet,\bullet\}\lra\ptt $ also means that $f$ is injective,\footnote{Indeed, 
$h(f(x_1))=t(x_1)$ and $h(f(x_2)=t(x_2)$ and we may pick $t:X\lra\{\bullet,\bullet\}$ such that $x_1\neq x_2$ whenever $x_1\neq x_2$.}
whereas in $\Topp$ the same lifting property defines injectivity on $\pi_0$ for ``nice'' spaces.\footnote{Indeed, 
if both from $X$ and $Y$ have finitely many connected components, we may consider the factorisation
$X\lra\pi_0(X)\lra\{\bullet\}$ and replace $f:X\lra Y$
 by $\pi_0(f):\pi_0(X)\lra \pi_0(Y)$ thereby reducing the question to $\Sets$.}

Note all these  notions are defined in terms of their simplest counterexample; in fact this is a common pattern
in basic definitions of topology and elsewhere. 
\begin{defi} The {\em left, resp.~right, Quillen negation} or {\em orthogonal} of property $P$ of morphisms is 
	$$P^l := \{ f : \forall g \in P\, f\rtt g \},\ \ 
	P^r := \{ g : \forall f \in P\, f\rtt g \}$$
\end{defi}
Let $$P^{lr}:=(P^l)^r, ...$$ 

By a diagram chasing calculation 
it is clear that 
$P\cap P^l$ and $P\cap P^r$ consist only of isomorphisms, and
$P\subset P^{lr} $ and $P\subset P^{rl}$, and $P^l=P^{lrl}$, and $P^r=P^{rlr}$.
The class $P^r$ 
is always closed under retracts\footnote{Recall that a morphism $f:A\lra B$ is a retract of a morphism $g:X\lra Y$ iff there is a commutative diagram 
\begin{equation}\begin{gathered}\label{Def:retract}
	\xymatrix{ {A} \ar@/^1.20pc/[rr]|{\operatorname{id}}  \ar[r] \ar@{->}[d]|{f} & {X} \ar[r] \ar[d]|g & A \ar[d]|{f} \\ {B} 
	\ar@/_1.20pc/[rr]|{\operatorname{id}}   \ar[r] &  {Y} \ar[r] & B  }
\end{gathered}\end{equation}
}, pullbacks,  products (whenever they exist in the category), and  composition of morphisms, and contains all isomorphisms (that is, invertible morphisms) of the underlying category. 
Meanwhile, $P^l$
is closed under retracts, pushouts,  coproducts, and transfinite composition (filtered colimits) of morphisms (whenever they exist in the category), and also contains all isomorphisms.

\subsubsection{An intuition of the Quillen lifting property as a category-theoretic negation}
Taking the Quillen negation is perhaps the simplest way to define a class of morphisms
without a given property in a manner useful in a diagram chasing computation,
and a number of basic definitions can be obtained by taking several times
the Quillen negation of their simplest counterexamples. 
In this paper, we see that {\em a number of standard definitions in topology
are in fact iterated Quillen negations of their simplest (counter)examples
which are maps of finite topological spaces.}

\subsection{Connectivity defined in terms of the simplest counterexample} 
\label{subsec:conn}
In this \S{} we take a simple example of non-connected space $\{\bullet,\bullet\}$ and a related morphism $\{\bullet,\bullet\}\lra\{\bullet\}$,
and show how taking its Quillen negations defines the notions of being {\em connected}, {\em totally disconnected} (for metrisable compact Hausdorff spaces),
and {\em a retract of $\{\bullet,\bullet\}^\kappa$ for some cardinal $\kappa$}. 
Later in Proposition~\ref{prop9} we shall uncover  more complicated morphisms (with $\leqslant 3$ points) 
implicit in \cite[I\S$11.5$,Proposition $9$]{Bourbaki} and that this Proposition hides a weak factorisation system which, in particular,
leads to a definition of totally disconnected for arbitrary spaces. 

\subsubsection{Connected and totally disconnected as Quillen negations}
Using the notation of Quillen negation, Eq.~(\ref{DefCon}) is expressed concisely as
\bd
A space $X$ is connected iff $X\lra\{\bullet\} \,\in\, \{ \{\bullet,\bullet\}\lra \{\bullet\}\}^{l}$. 
\ed

Eq.~(\ref{DefDisCon}) implies that 
\bd
A space $Q$ is 
totally disconnected if $Q\lra\{\bullet\} \,\in\, \{ \{\bullet,\bullet\}\lra \{\bullet\}\}^{lr}$. 
\ed

Proposition~\ref{propConnDisconn} below describes the meaning of these Quillen negations/orthogonals.
We sketch its content and proof before stating it in detail.

If a space $Q$ is ``nice'' in the sense that its connected components are open, 
then 
\bd
$\pi_0(X)\lra \pi_0(Q)$ is injective iff 
$X\lra Q \,\in\, \{ \{\bullet,\bullet\}\lra \{\bullet\}\}^{l}$. 
\ed

A diagram chasing argument\footnote{This follows from the lifting property
        $ Q\xra{\{ X\lra \{\bullet\}\}^{l} }  \prod\limits_{f:Q\lra X} X \,\rtt\, Q \xra{\{ X\lra \{\bullet\}\}^{lr} } \ptt $}
shows that $Q$ is a retract of 
a (possibly infinite) Cartesian power of $X$ 
iff $Q\lra\{\bullet\} \,\in\, \{ X\lra \{\bullet\} \}^{lr} $. 

It is a standard fact that each closed subspace of the Cantor space $\{\bullet,\bullet\}^\omega$
is a retract of it, and that
a compact Hausdorff metrisable space is totally disconnected iff it is a retract of 
the Cantor space $\{\bullet,\bullet\}^\omega$. 
As left Quillen negations are closed under retracts and products, this means that
\bd
A metrisable 
space $Q$ is profite, equiv.
totally disconnected compat Hausdorff,
iff $$Q\lra\{\bullet\} \,\in\, \{ \{\bullet,\bullet\}\lra \{\bullet\}\}^{lr}$$ 
\ed
%

\begin{enonce}{Proposition}\label{propConnDisconn} The following holds in the category of all topological spaces.
\bee
\item  \label{propConnDisconnlXtoO}A space $X$ is connected iff $X\lra\{\bullet\} \,\in\, \{ \{\bullet,\bullet\}\lra \{\bullet\}\}^{l}$. 
\item  \label{propConnDisconnlXtoQ} 
	$X\lra Q  \,\in\, \{ \{\bullet,\bullet\}\lra \{\bullet\}\}^{l}$
	iff   
	\begin{itemize} 
			\item 
	           each clopen subset of $X$ is the preimage of a clopen subset of $Q$ 
	\end{itemize}
	In particular, this holds if either  \bi\item $X$ is connected 
	\item $X\lra Q$ is a quotient map with connected fibres. \ei
		If $\pi_0(Q)$ is discrete, 
		this is equivalent to $\pi_0(X)\lra \pi_0(Q)$ being injective.
\item \label{propConnDisconnlrQtoO}   
A space $Q$ is a retract of $\{\bullet,\bullet\}^\kappa$ for some ordinal $\kappa$
iff $$Q\lra\{\bullet\}\,\in\, \{\{\bullet,\bullet\}\lra \{\bullet\}\}^{lr}$$
In particular, such a space $Q$ is necessarily totally disconnected compact Hausdorff.
\item  \label{propConnDisconnlrOtoOmetr}
A metrisable space $Q$ is profinite, or equivalently
totally disconnected compact Hausdorff, 
	iff $$Q\lra\{\bullet\} \,\in\, \{ \{\bullet,\bullet\}\lra \{\bullet\}\}^{lr}$$

\item  \label{propConnDisconnM2llr} Each morphism $f:X\lra Y$ decomposes as $f=f_l\circ f_{lr}$ 
	where $f_l\in \{\{\bullet,\bullet\}\lra\{\bullet\}\}^l$ and $f_{lr} \in \{\{\bullet,\bullet\}\lra\{\bullet\}\}^{lr}$.
In notation, there exists a decomposition 
$$
		X \xra  [f_{l}]{ \{\{\bullet,\bullet\}\lra\{\bullet\}\}^{l} } Q \xra  [f_{lr}]{ \{\{\bullet,\bullet\}\lra\{\bullet\}\}^{lr} }  Y $$ 

\eee\end{enonce}
		\begin{proof} (1,2). The first two items are  merely a reformulation in the new notation
	of the lifting property $X\lra\{\bullet\}\,\rtt\,\{\bullet,\bullet\}\lra\{\bullet\}$ and $X\lra Q\,\rtt\,\{\bullet,\bullet\}\lra\{\bullet\}$.
        If $X$ is connected, this lifting property holds trivially as there are no non-trivial maps $X\lra \{\bullet,\bullet\}$. 
	If the fibres of $\tau:X\lra Q$ are connected, then necessarily for connected subset $X_a$ of $X$ 
	it holds $X_a=\tau\inv\tau(X_a)$. If, moreover, $\tau:X\to Q$ is a quotient map, then $X_a=\tau\inv\tau(X)$ being clopen 
	implies that $\tau(X_a)$ is also clopen.\footnote{The same argument appears in the proof of  
	\cite[I\S$11$(Connectedness), Proposition 7]{Bourbaki}.}  

	(3). One direction follows from the fact that a left Quillen negation is closed under retracts and products.
	
	The other direction follows from a diagram chasing calculation, as follows. 
         Take $I:=\prod\limits_{f:Q\lra\{\bullet,\bullet\}} \{\bullet,\bullet\}$, and $Q\lra I$ to be the obvious map 
	 which on the $f$-th coordinate is given by $f$. 
	The lifting property $Q\xra\iota I \,\rtt\, \{\bullet,\bullet\}\lra\{\bullet\}$ holds by construction of $I$ 
	(the lifting diagonal is given by the projection corresponding to the horizontal map $I\lra \{\bullet,\bullet\}$).
	The retraction $I\lra Q$ is constructed using the lifting property $I\lra Q \,\rtt\, Q \lra \{\bullet\}$.
	The following diagrams represent this argument: 
	$$\begin{gathered}	\xymatrix{ Q \ar[d]^{\iota}_{ \{\{\bullet,\bullet\}\lra \{\bullet\}\}^{l}} \ar[r]|{f} &  \{\bullet,\bullet\}  \ar[d]
	\\  
	\prod\limits_{f:Q\lra\{\bullet,\bullet\}} \{\bullet,\bullet\}
		\ar[r] \ar@{-->}[ru]|{\operatorname{pr}_{\!f}}  & { \{\bullet\}} } %
		\xymatrix{ \\ & \text{ means } &   Q\lra I \,\in\, \{ \{\bullet,\bullet\}\lra\{\bullet\} \}^l & \\ }
	\end{gathered}$$ $$ 
	\begin{gathered}
	\\ \xymatrix{ Q \ar[d]^{\iota}_{ \{\{\bullet,\bullet\}\lra \{\bullet\}\}^{l}} \ar[r]|{\operatorname{id}} & Q \ar[d]^{ \{\{\bullet,\bullet\}\lra \{\bullet\}\}^{lr} }   \\  
	\prod\limits_{f:Q\lra\{\bullet,\bullet\}} \{\bullet,\bullet\}
	\ar[r] \ar@{-->}[ru]  & { \{\bullet\}} } 
		\xymatrix{ \\ & \text{ means } &   Q \text{ is a retract of }I & \\ }
	\end{gathered}$$%

	(4) We need to show that a metrisable space is totally disconnected compact Hausdorff iff it is a retract of some power of $\{\bullet,\bullet\}^\kappa$.
	
Such a retract is the image of a map of a compact Hausdorff space into itself, hence it is a closed subset of $\{\bullet,\bullet\}^\kappa$,
	hence also compact Hausdorff. It is totally disconnected because each subset of a totally disconnected space is totally disconnected, 
as	a non-trivial connected subset of a subset would also be a non-trivial connected subset of the space itself.

        Now assume that  $Q$ is metrisable and totally disconnected compact Hausdorff, and hence
	zero-dimensional.  This means the clopen subsets form a basic of topology, 
	i.e. 
	each point has a clopen neighbourhood of diameter $<\varepsilon$ for each $\varepsilon>0$. 
	For each $n>0$ use compactness to find a finite covering $Q$ by clopen neighbourhoods of diameter $<\varepsilon=2^{-n}$. 
	Such a clopen neighbourhood corresponds to a continuous function $Q \lra \{\bullet,\bullet\}$, and thus such a finite covering gives us a map 
	$\chi_n:Q\lra \{\bullet,\bullet\}^{N_n}$ for some finite $N_n$ such that $\chi_n(x)\neq \chi_n(y)$ whenever $\dist(x,y)>2^{-n}$. 
	Hence, the map $Q\lra \prod_n \{\bullet,\bullet\}^{N_n}=\{\bullet,\bullet\}^\omega$ is injective. As $Q$ is compact Hausdorff, its image is closed,
	and the topology on $Q$ is the subspace topology. 

	In conclusion, $Q$ is a closed subspace of the Cantor space $\{\bullet,\bullet\}^\omega$, and it is 
	a standard fact that  a closed subspace of the Cantor space 
	is necessarily its retract. 
	For completeness we sketch a proof. Define a metric on  $\{\bullet,\bullet\}^\omega$ such that $\dist(x,y)=\dist(x',y)$ implies $x=x'$, e.g.~as follows:
	$\dist( (x_n)_n, (y_n)_n) := \sum\limits_{n>0} \frac {|x_n-y_n|}{10^n}$.
	This property of the metric implies  means for each $x$ there is a unique nearest point of $Q$ in this metric,
	and hence the contraction $r:\{\bullet,\bullet\}^\omega\lra Q$ taking each point $x$ into the nearest point of $Q$ is well-defined. 
	To see that $r$ is continuous, it is enough to show that for an arbitrary sequence  $(c_n)_n\lra c$ converging to a point $c\in Q$
	we also have that $r(c_n)\lra r(c)$. Pick a convergent subsequence $(r(c_{n_i}))_i \lra c'$ by compactness. 
	Then $\dist(c',c)= \lim_i \dist(r(c_{n_i}),c_{n_i}))=\lim_i \dist(A,c_{n_i})= \dist(A,c)=\dist(r(c),c)$,
	and hence $c'=r(c)$ by the property of the metric. 

(5). 
Let us first show this for $Y=\{\bullet\}$. 
         Take $Q:=\prod\limits_{f:X\lra\{\bullet,\bullet\}} \{\bullet,\bullet\}$, and $X\lra Q$ to be the obvious map 
	 which on the $f$-th coordinate is given by $f$: by definition, 
	 this map glues together  two points $x,y\in X$ 
	 iff $y$ lies in the {\em quasi-component of $x$}, i.e.~the intersection of all the closed open subsets containing $x$.  
	The lifting property $X\xra\iota I \,\rtt\, \{\bullet,\bullet\}\lra\{\bullet\}$ holds by construction of $I$: a horizontal map $f:X\lra\{\bullet,\bullet\}$
	lifts to the diagonal map $I\lra \{\bullet,\bullet\}$ given by the projection on the $f$-th coordinate.
On the other hand, $Q\lra \{\bullet\}\,\in\,  
	\{\{\bullet,\bullet\}\lra\{\bullet\}\}^{lr}$ as Quillen negations are closed under products.

Now let $Y$ be arbitrary. We claim that 
$$
		X \xra  [f_{l}]{ \{\{\bullet,\bullet\}\lra\{\bullet\}\}^{l} } Q\times Y \xra  [f_{lr}]{ \{\{\bullet,\bullet\}\lra\{\bullet\}\}^{lr} }  Y $$ 
As $(\cdot)^{r}$-classes are closed under pullbacks, $Q \times Y \lra Y \,\in\, \{\{\bullet,\bullet\}\lra\{\bullet\}\}^{lr}$.
		To see that $X \lra Q \times Y \,\in\, \{\{\bullet,\bullet\}\lra\{\bullet\}\}^{l}$, note that 
		to find a decomposition for $\tau:X \lra \{\bullet,\bullet\}$ 
		though $X \lra Q \times Y$ as required by the lifting property $X\lra Q\times Y \,\rtt\, \{\bullet,\bullet\}\lra\{\bullet\}$,
it is enough to find such a decomposition through $X \lra Q \times Y \lra Q$, which we know is in $\{\{\bullet,\bullet\}\lra\{\bullet\}\}^{l}$ by construction. 
          
		The following diagrams represent the arguments above:
\begin{equation}
	\xymatrix@C=2cm{ \cdot \ar[d]\ar[r]|{f\times g}  & X_{lr}\times Y \ar[d] \\ \cdot\ar[r]|h \ar@{-->}[ru]|{\tilde f \times h}   & Y }
	\xymatrix@C=2cm{ X \ar[d]_{\therefore\{\{\bullet,\bullet\}\lra\{\bullet\}\}^{l} } 
	\ar[rr] \ar@/_2pc/@{-->}[rr] \ar[rd]|{\{\{\bullet,\bullet\}\lra\{\bullet\}\}^{l} } & & \{\bullet,\bullet\} \ar[d] \\ 
	X_{lr}\times Y  \ar@{->}[r]  & X_{lr} \ar@{-->}[ru] \ar[r] & 
	\cdot } \label{M2llrConnectedProof}
\end{equation}
\end{proof}

\begin{rema} In item \ref{propConnDisconnlXtoO} we may take  any non-injective map $f$ with discrete fibres instead of $\{\bullet,\bullet\}\lra\ptt$. 
	Items \ref{propConnDisconnlrQtoO} and \ref{propConnDisconnM2llr} of Proposition~\ref{propConnDisconn}
	apply to any orthogonal of form $(F\lra \ptt)^{lr}$; we never use in the proof that $F=\{\bullet,\bullet\}$. 
\end{rema}

\subsection{Connected components  
as an instance of a 
weak factorisation system}\label{subsec:prop9}
We now show that \cite[I\S$11.5$,Proposition~9]{Bourbaki}, quoted below, describes an instance of a weak factorisation system
generated by 4 maps of finite topological spaces. The proposition says that the equivalence relation ``to lie 
in the same connected component'' is well-defined, has closed equivalence classes, and the quotient is a totally disconnected space. 
The 4 maps are simple (counter)examples failing
the properties of being connected, surjective, being a quotient, and having connected fibres.
%
\newline\includegraphics[test]{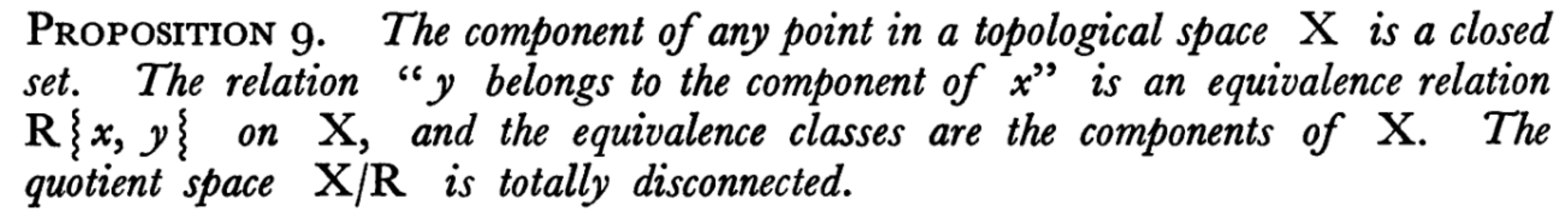}
We find it remarkable that \cite{Bourbaki} chooses  to state explicitly an instance of a weak factorisation system
related to a notion they discuss. As everywhere in the paper, here our exposition focuses on how to ``read off'' our category-theoretic reformulation
from the text of Bourbaki, and more mathematically inclined readers may prefer to skip directly to Proposition~\ref{prop9}, and guess rather than read 
our notation for finite spaces in Appendix~\S\ref{not:mintsGE} as necessary, as we try to give an explicit definition
of each finite space we use. 


\subsubsection{Surjective and quotient maps}

Let $R\subset X\times X$ be an equivalence relation on a topological space $X$.


``Let $\phi$ be the  canonical mapping $\phi:X\lra X/R$. 
By definition (\cite[I\S2, no. 4, Proposition $6$ and its corollary]{Bourbaki}) 
 the {\em open} (resp.~{\em closed}) subsets  in $X/R$
are the subsets $U$ such that $\phi\inv(U)$ is {\em open} (resp.~{\em closed}) in $X$.'' \cite[I\S$3.4$, just before Proposition $6$]{Bourbaki}:
These words describe the following lifting property:
	$$X \xra \phi  X/R\,\rtt\, \,  \OutocTooc \, $$
\def\rrt#1#2#3#4#5#6{\xymatrix@C=2.39cm{ {#1} \ar[r]|{\forall \tilde\xi} \ar@{->}[d]|{#2} & {#4} \ar[d]|{#5} \\ {#3}  \ar[r]|{\forall\xi} \ar@{-->}[ur]^{\exists}& {#6} }}
$$ \rrt  {X}  {\phi} {X/R}  {\Outoc} {} {\ouc } $$
Indeed, the lower arrow $\xi:X/R\lra \ouc$ varies though 
precisely the arbitrary subsets of $X/R$ (take $U:=\xi\inv(\bullet_u)$;
the topology on $\ouc$ is indiscrete and thus there are no continuity requirements on $\xi:X/R\lra\ouc$); 
the upper arrow $\tilde\xi:X\lra \Outoc$ describes precisely the open subsets of $X/R$ (take $\tilde U:=\tilde\xi\inv(\bullet_u)$;
now the point $\bullet_u$ is open in $\Outoc$ meaning that its preimage is open is the only continuity restriction on $\tilde\xi:X\lra \Outoc$).
The commutativity of the square means precisely that $\tilde U=\phi\inv(U)$. 
In short, the commutative square describes those subsets $U$ of $X/R$ such that $\phi\inv(U)$ is open. 
Finally, note that the diagonal $X\lra \Outoc$ describes open (equiv., closed) subsets of $X/R$, hence 
the diagram 
says precisely the phrase of Bourbaki 
``the {\em open} (resp.~{\em closed}) subsets  in $X/R$
are the subsets $U$ such that $\phi\inv(U)$ is {\em open} (resp.~{\em closed}) in $X$.'' \cite[I\S$3.4$, just before Proposition $6$]{Bourbaki}.

The canonical mapping $\phi:X\lra X/R$ is surjective. 
Surjectivity  can also be defined by a left lifting property: a map $X\lra Y$ is surjective iff\footnote{
Indeed, the square commutes iff the preimage of $\star_a$ in $Y$ 
contains the image of $X$.
The lower triangle commutes iff the preimage of $\star_a$ in $Y$ is the whole of $Y$. 
Thus, the diagram says that for each map $X/R\lra \{\star_a \leftrightarrow \star_b\}$, the preimage of $\star_a$ is the whole of $Y$ 
whenever it contains the image of $X$, i.e.~that the image of $X$ is the whole of $Y$.}
$$X \lra Y \,\rtt\, \{\star_a\} \lra \{\star_a \leftrightarrow \star_b\}$$
   \def\rrt#1#2#3#4#5#6{\xymatrix{ {#1} \ar[r]|{} \ar@{->}[d]|{#2} & {#4} \ar[d]|{#5} \\ {#3}  \ar[r] \ar@{-->}[ur]^{}& {#6} }}
$$ \rrt  {X}  {} {Y}  {\{\star_a\}} {} {\{\star_a \leftrightarrow  \star_b\} } $$

We summarise the remarks above in the following proposition. 

\begin{prop}\label{PropQuotientl}
\bee\item 
A map $A \xra f B$ is surjective iff 
$X \xra f B \,\rtt\, \{a\} \lra \{a\llrra b \}$ 
\item The following two properties of a map $A \xra f B $ are equivalent: 
\bi\item a subset $Z$ of $B$ is open iff its preimage $f\inv(Z)\subset A$ is open; in particular $\Imm f$ is both open and closed. 
\item $X \xra f B \,\rtt\, \OcTooc $ 
	\ei
\eee

In particular, the class
	$$ \left\{\, \{a\} \lra \{a\llrra b \} \, ,\,  \OcTooc \,\right\}^l $$
is precisely the class of canonical mappings $X\lra X/R$ where $R$ is an equivalence relation on a topological space $X$.
\end{prop}
\begin{proof} (1) View  the image $f(A)\subset B$ as its characteristic  map to $\chi:B\lra \{a \llrra b\}$, $\chi(f(A)):=\{a\}$, $ \chi(B\setminus f(A)):=\{b\}$.
	This maps factors though $a$ iff $B\setminus f(A)=\emptyset$. 
	(2) A subset $Z$ of $B$ gives rise to its characteristic map $\chi:B\lra\{\bullet_u\llrra \bblacksquare_c\}$ where $\chi(Z):=\{\bullet_u\}$ 
	and $\chi(B\setminus Z):=\{\bblacksquare_c\}$. 
	This map factors via $\Oc$ iff the subset is open. The map $A \lra \oc $ represents the open subset of $A$ which is the preimage of $Z$.
        (3) In view of the items above, this says 
	a map is of this form iff it is a surjective map such that a subset is open iff its preimage is open.
\end{proof}

\subsubsection{Separation axiom $T_1$.} Saying that ``the component of any point in a topological space $X$ is closed''
is equivalent to saying that each point of $X/R$ is closed, i.e.~it satisfies separation Axiom $T_1$. 
Axiom $T_1$ can be defined using perhaps the simplest example of a non-$T_1$ space which is the Alexandroff-Sierpinski space 
$\left\{{\ensuremath{\bullet}}\rotatebox{-12}{\ensuremath{\to}}  \raisebox{-2pt}{\ensuremath{\bblacksquare}}\right\}$
where the point $\bullet$ is not closed:

\begin{prop}\label{propTone} A space $X$ satisfies the separation axiom $T_1$, i.e.~each point of $X$ is closed, iff
	$$
	\left\{{\ensuremath{\bullet}}\rotatebox{-12}{\ensuremath{\to}}  \raisebox{-2pt}{\ensuremath{\bblacksquare}}\right\}
	\lra \{\bullet\} 
	\,\rtt\,
	X\lra \{\bullet\}$$
\end{prop}\begin{proof}A verification shows that 
to give a map $\left\{{\ensuremath{\bullet}}\rotatebox{-12}{\ensuremath{\to}}  \raisebox{-2pt}{\ensuremath{\bblacksquare}}\right\} \xra f X$ 
	is the same as to give a pair of points $x\in \cl(y)$ in $X$  where $x:=f(\bblacksquare)$ and $y:=f(\bullet)$. 
	This map fits into a lifting square iff $x=y$.
\end{proof}
\begin{coro} An equivalence relation $R$ on a topological space $X$ has closed equivalence classes iff
	each point of $X/R$ is closed, i.e.
	$$
	\left\{{\ensuremath{\bullet}}\rotatebox{-12}{\ensuremath{\to}}  \raisebox{-2pt}{\ensuremath{\bblacksquare}}\right\}
	\lra \{\bullet\} 
	\,\rtt\,
	X/R\lra \{\bullet\}$$
\end{coro}
\begin{proof} Immediate by the definition of the quotient topology on $X/R$.\end{proof}

\subsubsection{Connected fibres} We want to define the class of maps with connected fibres.
By definition, a map $f:X\lra Y$ does not have connected fibres iff there is  $y\in Y$ such that the fibre $f\inv(y)=A\cup B$ 
is a union of two disjoint open (resp.~closed) subsets of the fibre $f\inv(y)\subset X$ with induced topology. 
Assume that $y$ is a closed point and that $X$ and $Y$ are connected, 
as would be true in many natural examples; in this case the fibre $f\inv(y)$ and its components $A$ and $B$
are closed subsets of $X$. Contracting in $X$ the subsets 
$X\setminus f\inv(y)$, $A$, and $B$, and contracting in $Y$ the subset $Y\setminus\{y\}$, 
leads to a simple example of a map of finite topological spaces which does not have connected fibres\footnote{If point $y$ is neither closed nor open, doing so
might 
lead to a map 
of finite indiscrete spaces with connected fibres.}
$$\left\{ 
 \raisebox{-2pt}{\ensuremath{\bblacksquare_a}\!\! \rotatebox{12}{\ensuremath{\leftarrow}}}\, 
	{\ensuremath{\bullet_u}}\rotatebox{-12}{\ensuremath{\to}}  \raisebox{-2pt}{\ensuremath{\bblacksquare_b}}\right\}
\lra
\left\{{\ensuremath{\bullet_u}}\rotatebox{-12}{\ensuremath{\to}}  \raisebox{-2pt}{\ensuremath{\bblacksquare_{a=b}}}\right\} 
$$
and a commutative square with does not admit a lifting:
$$\xymatrix@C=2.39cm{ X \ar[r] \ar[d]|f  & \left\{
	 \raisebox{-2pt}{\ensuremath{\bblacksquare_a}\!\! \rotatebox{12}{\ensuremath{\leftarrow}}}\,
	         {\ensuremath{\bullet_u}}\rotatebox{-12}{\ensuremath{\to}}  \raisebox{-2pt}{\ensuremath{\bblacksquare_b}}\right\}
		 \ar[d] 
		 \\
Y \ar[r] &  
\left\{{\ensuremath{\bullet_u}}\rotatebox{-12}{\ensuremath{\to}}  \raisebox{-2pt}{\ensuremath{\bblacksquare_{a=b}}}\right\} 
}
\xymatrix@C=0.64cm{& }
\xymatrix@C=2.39cm{ X \ar[r] \ar[d]|f  & \left\{
	 \raisebox{-2pt}{\ensuremath{\bblacksquare_A}\!\! \rotatebox{12}{\ensuremath{\leftarrow}}}\,
		 {\ensuremath{\bullet_{X\setminus f\inv(y)}}}\rotatebox{-12}{\ensuremath{\to}}  \raisebox{-2pt}{\ensuremath{\bblacksquare_B}}\right\}
		 \ar[d] 
		 \\
Y \ar[r] &  
\left\{{\ensuremath{\bullet_{Y\setminus\{y\}}}}\rotatebox{-12}{\ensuremath{\to}}  \raisebox{-2pt}{\ensuremath{\bblacksquare_{y}}}\right\} 
\\& f\inv(y)=A\sqcup B } $$
Here 
$\left\{{\ensuremath{\bullet_u}}\rotatebox{-12}{\ensuremath{\to}}  \raisebox{-2pt}{\ensuremath{\bblacksquare_{a=b}}}\right\} 
$ denotes the space with one open point $\bullet_u$ and one closed point $ \bblacksquare_{a=b}$,
and 
$\left\{
	 \raisebox{-2pt}{\ensuremath{\bblacksquare_a}\!\! \rotatebox{12}{\ensuremath{\leftarrow}}}\,
	         {\ensuremath{\bullet_u}}\rotatebox{-12}{\ensuremath{\to}}  \raisebox{-2pt}{\ensuremath{\bblacksquare_b}}\right\}$
 denotes the space with one point $\bullet_u$ and 
two closed points $ \bblacksquare_{a}$ and  $ \bblacksquare_{a}$. The map glues the two closed points together, as subscripts may suggest,
and thus the fibre above the closed point $\bblacksquare_{a=b}$ is not connected.
On the right, we use subscripts to indicate intended preimages. 
\begin{prop}\label{ConnTOfibresl}
Let $B$ be a $T_1$-space, i.e.~each point of $B$ is closed. 
Then a quotient map $A \xra f B$ has connected fibres iff
$$A  \xra f B \,\rtt\,
	\left\{
		 \raisebox{-2pt}{\ensuremath{\bblacksquare_a}\!\! \rotatebox{12}{\ensuremath{\leftarrow}}}\,
		         {\ensuremath{\bullet_u}}\rotatebox{-12}{\ensuremath{\to}}  \raisebox{-2pt}{\ensuremath{\bblacksquare_b}}\right\}
		 \lra
		 \left\{{\ensuremath{\bullet_u}}\rotatebox{-12}{\ensuremath{\to}}  \raisebox{-2pt}{\ensuremath{\bblacksquare_{a=b}}}\right\}
$$\end{prop}
\begin{proof}
	To give a commutative square is the same as to give a closed subset $Z\subset B$ and a decomposition of its preimage $f\inv(Z)\subset A$
	as a union of two closed subsets $f\inv(Z)=\tilde Z_1\cup \tilde Z_2$. 
	To give a lifting map is the same as to give a decomposition of $Z$ into two closed subsets $Z=Z_1\cup Z_2$
	such that $\tilde Z_1= f\inv(Z_1)$ and  $\tilde Z_2= f\inv(Z_2)$. 
	
	In particular, if the lifting property holds, then the preimage of a closed connected subset is necessarily connected.
	Hence, the fibres of $f$ above closed points are connected.

	Now assume the fibres are connected. Then each fibre over a point $z\in Z$ is contained in either $\tilde Z_1$ or $\tilde Z_2$. Hence,
	$\tilde Z_1= f\inv( Z_1)$ and $\tilde Z_2=f\inv(Z_2)$ for $Z_1:=f(\tilde Z_1)$ and $Z_2:=f(\tilde Z_2)$. 
        As $f:A\lra B$ is a quotient map, the images $Z_1$ and $Z_2$ are closed, and we get a required decomposition of $Z$.
\end{proof}
Note that the assumption that $f:A\lra B$ is a quotient map was necessary: an injective map from a discrete set of size at least two to an indiscreet set 
does not have this lifting property.

\begin{prop}\label{prop9} {\cite[I\S$11.5$,Proposition $9$]{Bourbaki}}  
can equivalently be formulated as follows: each morphism $X\lra\{\bullet\}$ admits a decomposition
	$$X \xra {(P_1)^l} X_{lr} 
	 \xra {(P_1)^{lr}} \{\bullet\} $$
	for the class $P_1$ consisting of the following 4 morphisms of finite topological spaces
	$$ 
	\{a\} \lra \{a\llrra b \} \, ,\,  \OcTooc \, 
, 	\left\{
		 \raisebox{-2pt}{\ensuremath{\bblacksquare_a}\!\! \rotatebox{12}{\ensuremath{\leftarrow}}}\,
		         {\ensuremath{\bullet_u}}\rotatebox{-12}{\ensuremath{\to}}  \raisebox{-2pt}{\ensuremath{\bblacksquare_b}}\right\}
		 \lra
		 \left\{{\ensuremath{\bullet_u}}\rotatebox{-12}{\ensuremath{\to}}  \raisebox{-2pt}{\ensuremath{\bblacksquare_{a=b}}}\right\}
	,$$ 
or, equivalently, for the class $P_2$ consisting of a single morphisms 
$$
 	\left\{
		 \raisebox{-2pt}{\ensuremath{\bblacksquare_a}\!\! \rotatebox{12}{\ensuremath{\leftarrow}}}\,
		         {\ensuremath{\bullet_u}}\rotatebox{-12}{\ensuremath{\to}}  \raisebox{-2pt}{\ensuremath{\bblacksquare_b}}
			 \right\}
		 \lra
		 \left\{\bullet_u\leftrightarrow\bblacksquare_{a=b}\leftrightarrow\filledstar_x\right\}
	$$ 
\end{prop}
\begin{proof} 
	Indeed, by Proposition~\ref{PropQuotientl} each map in $(P_1)^l$ is of form $X\lra X/R$ for some equivalence relation $R$. 
 
	The first claim in Proposition~$9$ is that each fibre of $X\lra X/R$ is closed, and by the definition of quotient topology 
	this is equivalent to saying that that each point of $X/R$ is closed, i.e.~by Proposition~\ref{propTone} 
	$$
	\left\{{\ensuremath{\bullet}}\rotatebox{-12}{\ensuremath{\to}}  \raisebox{-2pt}{\ensuremath{\bblacksquare}}\right\}
	\lra \{\bullet\} 
	\,\rtt\,
	X/R\lra \{\bullet\}$$
	A verification shows that 
	$
	\left\{{\ensuremath{\bullet}}\rotatebox{-12}{\ensuremath{\to}}  \raisebox{-2pt}{\ensuremath{\bblacksquare}}\right\}
	\xra{(P_1)^l} \{\bullet\} $, hence 
	$X_{lr} \xra {(P_1)^{lr}} \{\bullet\}$ implies that each point of $X_{lr}:=X/R$ is closed. 
	This gives us the first claim of the Proposition~9 that the equivalence classes of $R$
	are closed.
	Furthermore, it also gives us that by Proposition~\ref{ConnTOfibresl}  $X\lra X/R$ lies in $(P_1)^l$ iff 
	the map $X\lra X/R$ has connected fibres. 
        
To sum up, we see that 
	the map $X\lra X_{lr}$ is necessarily of form of a quotient map $X\lra X/R$ where $R$ is an equivalence relation on $X$ 
	with connected closed fibres.
	If $X/R$ is totally disconnected, then each equivalence class of $R$ maps to a single point, as the image of a connected set is necessarily connected.
	Hence, 
	the relation $R$ is the equivalence relation ``$y$ belongs to the connected component of $x$'' iff the space $X_{lr}=X/R$ is totally disconnected.
Thus it is only left to show that a space $Q$ is totally disconnected iff $Q\xra{(P_1)^{lr}}  \{\bullet\}$. 
We do so in Lemma~\ref{lem:totDisc} below.

A combinatorial argument shows that $P_1^l=P_2^l$, as follows. 
	The map in $P_2$ is the composition 
	of two maps from $P_1$ and a pullback of a map in $P_1$ (see the diagram below on the left), namely  
$$
 	\left\{
		 \raisebox{-2pt}{\ensuremath{\bblacksquare_a}\!\! \rotatebox{12}{\ensuremath{\leftarrow}}}\,
		         {\ensuremath{\bullet_u}}\rotatebox{-12}{\ensuremath{\to}}  \raisebox{-2pt}{\ensuremath{\bblacksquare_b}}\right\}
		 \lra
		 \left\{{\ensuremath{\bullet_u}}\rotatebox{-12}{\ensuremath{\to}}  \raisebox{-2pt}{\ensuremath{\bblacksquare_{a=b}}}\right\}
		 \lra
		 \left\{{\ensuremath{\bullet_u}}\leftrightarrow   {\ensuremath{\bblacksquare_{a=b}}}\right\}
		 \lra 
\left\{\bullet_u\leftrightarrow\bblacksquare_{a=b}\leftrightarrow\filledstar_x\right\} $$
As right Quillen negations are closed under composition and taking pullback, it implies $P_2\subset P_1^{lr}$ and therefore $P_2^l \supset P_1^{lrl}=P_1^l$. 
To see the converse, first notice that the first two maps in $P_1$ are retracts of the map in $P_2$. 
It is left to show that 
$f \,\rtt\, 
 	\left\{
		 \raisebox{-2pt}{\ensuremath{\bblacksquare_a}\!\! \rotatebox{12}{\ensuremath{\leftarrow}}}\,
		         {\ensuremath{\bullet_u}}\rotatebox{-12}{\ensuremath{\to}}  \raisebox{-2pt}{\ensuremath{\bblacksquare_b}}\right\}
		 \lra 
\left\{\bullet_u\leftrightarrow\bblacksquare_{a=b}\leftrightarrow\filledstar_x\right\} $
implies
$f \,\rtt\, 
		 \left\{{\ensuremath{\bullet_u}}\rotatebox{-12}{\ensuremath{\to}}  \raisebox{-2pt}{\ensuremath{\bblacksquare_{a=b}}}\right\}
		 \lra
		 \left\{{\ensuremath{\bullet_u}}\leftrightarrow   {\ensuremath{\bblacksquare_{a=b}}}\right\}
$.
To see this,
consider the diagram below on the right and notice that
a lifting makes the square commute iff it makes the outer diagram of solid arrows commute, 
as the vertical arrow 
$
		 \left\{{\ensuremath{\bullet_u}}\rotatebox{-12}{\ensuremath{\to}}  \raisebox{-2pt}{\ensuremath{\bblacksquare_{a=b}}}\right\}
		 \lra
		 \left\{{\ensuremath{\bullet_u}}\leftrightarrow   {\ensuremath{\bblacksquare_{a=b}}}\right\}
$ is surjective. 
%
%
%
%
$$
\xymatrix@C=1.09cm{ 
	\left\{\bullet_u\leftrightarrow\bblacksquare_{a=b}\phantom{\leftrightarrow\filledstar_x}\right\} \ar[r] \ar[d] &\left\{ \star_{u=a=b}\phantom{\leftrightarrow \filledstar_x}\right\} \ar[d] 
	\\ \left\{\bullet_u\leftrightarrow\bblacksquare_{a=b}\leftrightarrow\filledstar_x\right\} \ar[r] & \left\{ \star_{u=a=b} \leftrightarrow \star_x \right\}}
\phantom{AAAAAAA}\xymatrix{ A \ar[d]|f \ar[r] &
\left\{{\ensuremath{\bullet_u}}\rotatebox{-12}{\ensuremath{\to}}  \raisebox{-2pt}{\ensuremath{\bblacksquare_{b}}}\right\} \ar[d] \ar[r]  
&\left\{
		 \raisebox{-2pt}{\ensuremath{\bblacksquare_a}\!\! \rotatebox{12}{\ensuremath{\leftarrow}}}\,
		         {\ensuremath{\bullet_u}}\rotatebox{-12}{\ensuremath{\to}}  \raisebox{-2pt}{\ensuremath{\bblacksquare_b}}
			 \right\} \ar[d] \\
		 B \ar[r]\ar@{-->}[ru]  \ar[r]\ar@{-->}[rru]         		 & 
		 \left\{{\ensuremath{\bullet_u}} \leftrightarrow   {\ensuremath{\bblacksquare_{b}}}\right\} \ar[r] & 
		 \left\{\bullet_u\leftrightarrow\bblacksquare_{a=b}\leftrightarrow\filledstar_x\right\} }$$\end{proof}

\subsubsection{A definition of totally disconnected} The following lemma defines {\em totally disconnected} 
in terms of finite spaces. Compare this reformulation to the reformulations in Proposition~\ref{propConnDisconn}(\ref{propConnDisconnlrOtoOmetr}) 
and Proposition~\ref{propDefDisCon}. 
\begin{lemm}\label{lem:totDisc} 
	A space $Q$ is totally disconnected iff $Q\xra{(P_1)^{lr}}  \{\bullet\}$. 
\end{lemm}
\begin{proof}
Let $Q\lra Q/R_Q$ be the canonical mapping from $Q$ to its totally disconnected quotient constructed by Proposition~9.
	By the argument above we have $Q\xra{(P_1)^{l}} Q/R_Q$. The lifting property 
	$Q\xra {(P_1)^{l}} Q/R_Q \,\,\rtt\,\, Q\xra{(P_1)^{lr}}  \{\bullet\}$ constructs a map $Q/R_Q\lra Q$ such that 
	the composition $Q\lra Q/R \lra Q$ is the identity $\id_Q:Q\lra Q$. This implies that the map $Q\lra Q/R_Q$ is injective,
	i.e.~that connected components of $Q$ are precisely points; by definition this means that $Q$ is totally disconnected.

	Now assume that $Q$ is totally disconnected. We want to show that $Q\xra{(P_1)^{lr}}  \{\bullet\}$, 
	i.e. that it holds 
	$A\xra {(P_1)^{l}} B \,\,\rtt\,\, Q\lra \{\bullet\}$. 
A map $A\lra Q$ sends each connected component of $A$ into a single point as $Q$ is totally disconnected,
hence it factors though the totally disconnected quotient $A/R_A$ of $A$. Hence, to show this lifting property
it is enough to show that 
	$A\xra {(P_1)^{l}} B \,\,\rtt\,\, A/R_A\lra  \{\bullet\}$. 
By Proposition~\ref{PropQuotientl}  $A\xra {(P_1)^{l}} B $ implies that $B$ is a quotient space of $A$ 
by some equivalence relation $S$, and the map is the canonical mapping $\psi:A \lra A/S$. 
By the universal property of the quotient space, to construct a continuous mapping $A/S \lra A/R_A$ as required,
it is enough to construct a possibly not continuous mapping $A/S \lra A/R_A$ as required.
That is, it is enough to show that $x,y$ lie in the same connected component whenever $\psi(x)=\psi(y)$.  
Assume $x,y$ do not lie in the same connected component. By Proposition~9 the connected components of $x$ and of $y$
are closed disjoint subsets of $A$, hence they determine a map $A\lra 
	\left\{ \raisebox{-2pt}{\ensuremath{\bblacksquare_a}\!\! \rotatebox{12}{\ensuremath{\leftarrow}}}\,
	                         {\ensuremath{\bullet_u}}\rotatebox{-12}{\ensuremath{\to}}  \raisebox{-2pt}{\ensuremath{\bblacksquare_b}}\right\}
				                  \lra
								   \left\{{\ensuremath{\bullet_u}}\rotatebox{-12}{\ensuremath{\to}}  \raisebox{-2pt}{\ensuremath{\bblacksquare_{a=b}}}\right\}$
sending $x$ to $\bblacksquare_a$ and $y$ to $\bblacksquare_b$. The lifting property 
	$$A \xra{(P_1)^l} B \,\,\rtt\,\, \left\{ \raisebox{-2pt}{\ensuremath{\bblacksquare_a}\!\! \rotatebox{12}{\ensuremath{\leftarrow}}}\,
                         {\ensuremath{\bullet_u}}\rotatebox{-12}{\ensuremath{\to}}  \raisebox{-2pt}{\ensuremath{\bblacksquare_b}}\right\}
			                  \lra
							   \left\{{\ensuremath{\bullet_u}}\rotatebox{-12}{\ensuremath{\to}}  \raisebox{-2pt}{\ensuremath{\bblacksquare_{a=b}}}\right\}$$
	implies that  
there is a map $B \lra 
	\left\{\raisebox{-2pt}{\ensuremath{\bblacksquare_a}\!\! \rotatebox{12}{\ensuremath{\leftarrow}}}\,
	                         {\ensuremath{\bullet_u}}\rotatebox{-12}{\ensuremath{\to}}  \raisebox{-2pt}{\ensuremath{\bblacksquare_b}}\right\}
				                  \lra
\left\{{\ensuremath{\bullet_u}}\rotatebox{-12}{\ensuremath{\to}}  \raisebox{-2pt}{\ensuremath{\bblacksquare_{a=b}}}\right\}$
	sending $\psi(x)$ to $\bblacksquare_a$ and $\psi(y)$ to $\bblacksquare_b$. 
Hence, 	$\psi(x)\neq \psi(y)$, as required.
\end{proof}
\begin{rema} It is outside the scope of this short note to prove that such a decomposition 
	$(P_1)^l(P_1)^{lr}$ decomposition in $\Topp$ does in fact exist. >>TODO: find a simple proof!!
\end{rema}

\section{Contractibility and trivial fibrations}\label{sec_T4} 

We show that how to define contractibility using a map implicit in the definition of a normal space.
The  proof of this is a diagram-chasing reformulation of a standard proof of the Urysohn theorem on extending real valued 
functions from a closed subset of a normal space in \cite[IX\S$4.2$]{Bourbaki}.

The map involved happens to be a trivial Serre fibration and thereby so is each map in its $lr$-orthogonal, which we conjecture 
to define the class of trivial fibrations in some model structure.

\subsection{The definition of a normal space}

\def\includegraphicss[#1]#2{    \pdfximage width  #1\linewidth {#2}  \pdfrefximage\pdflastximage}


....
Imagine we are trying to formalise 
the text of Axiom $(\mathrm O '_\mathrm{V})$ in Theorem $1$ (Urysohn) in \cite[IX\S$4.1$]{Bourbaki}, 
and, like in the $60$s,
we care about the number of bytes in formalisation. 


\noindent\includegraphics[width=\linewidth]{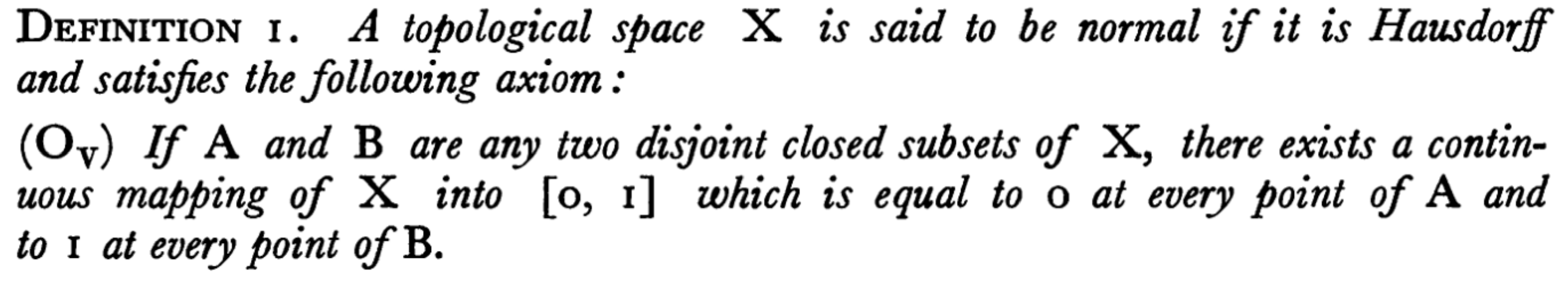}

\noindent\includegraphics[width=\linewidth]{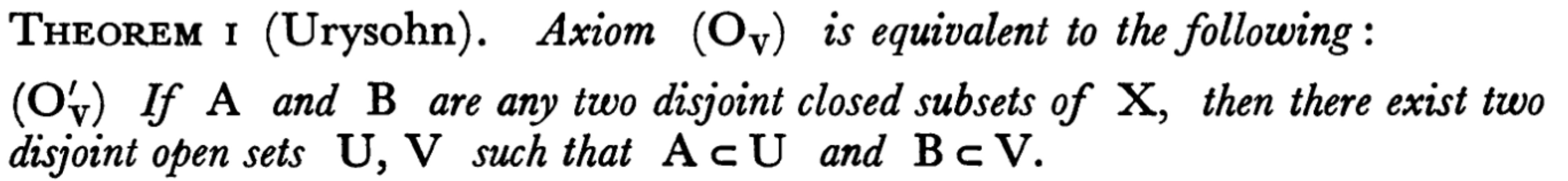}

\subsubsection{Preliminaries: the barycentric subdivision of the interval.}\label{def:MtoL:above}
We now define a map important in this section. 
Represent the closed interval as the union of two closed points (its endpoints) and an open interval: 
$[0,1]=\{0\}\cup (0,1) \cup \{0\}$. Contracting the open ``cell'' $(0,1)$ in this ``cell decomposition'' to a point  
gives a finite space with two closed points and one open point which we denote by $\vvLambda$, and by $\xi_\vvLambda:[0,1]\lra \vvLambda$ 
we denote the contracting map. Now ``subdivide'' the interval into two intervals, i.e.~take the barycentric subdivision of this cell decomposition: 
$$[0,1]=\{0\}\cup (0,\frac12)\cup \{\frac12\} \cup (\frac12,1)\cup \{1\}$$ 
Contracting the two open cells (i.e.~intervals)  $(0,\frac12)$ and $(\frac12,1)$ to points 
gives a finite space with three closed points and two open points which we denote by $\MM$, and the contracting map by $\xi_\MM:[0,1]\lra \MM$.
Let $\MtoL$ denote the corresponding map from $\MM$ to $\vvLambda$. 
The pictures below illustrate this construction. On the second picture, we indicate the open subsets by ovals. The third picture is a notation for this map. 
	\begin{minipage}[c]{0.5\linewidth} \pdfximage width  \linewidth {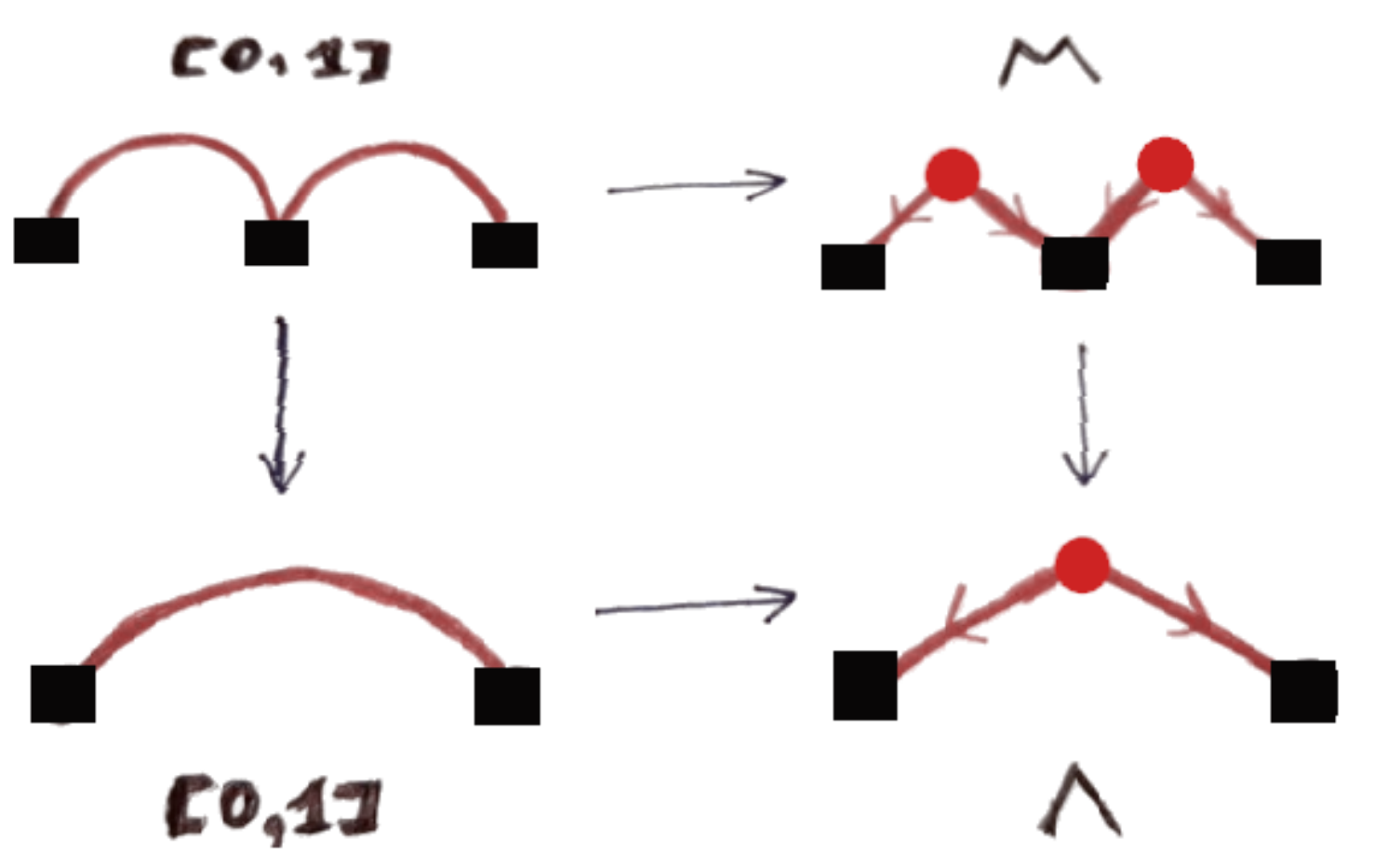} 
	\pdfrefximage\pdflastximage 
	\end{minipage}\begin{minipage}[c]{0.125\linewidth} \pdfximage width \linewidth {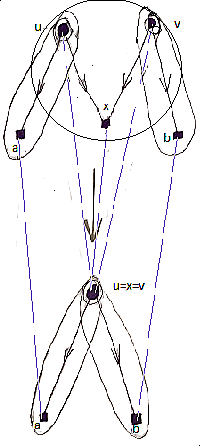} 
	\pdfrefximage\pdflastximage 
\end{minipage}\begin{minipage}[c]{0.33\linewidth}
\begin{center} $\left\{\underset{a}{}{\swarrow} \overset{\raisebox{1pt}{$u$}}{}{\searrow} \underset{x}{}{\swarrow}\overset{\raisebox{1pt}{$v$}}{}{\searrow}\underset{b}{}\right\}$\hskip2pt
\vskip5pt$\color{blue}\Downarrow$\vskip4pt
            \hskip2pt$\left\{\underset{a}{}{\swarrow}\overset{\raisebox{1pt}{$u\!=\!x\!=\!v$}}{}{\searrow} \underset{b}{}\right\}$\end{center}
\end{minipage}

The map $\MtoL$ is an acyclic Serre fibration, as we shall see later, and may be viewed as a finite model of the barycentric cell decomposition of the 1-simplex $[0,1]$. 

\subsection{A trivial Serre fibration implicit in the definition of normality (Separation Axiom $T_4$)}
We start with one way to see the map of finite spaces implicit in the definition of normality. 
The explanation following Eq.~\eqref{T4Odv} gives another way focusing on how to ``read off'' 
the map from the {\em text} of Axiom  $(\mathrm O '_\mathrm{V})$ as written in Bourbaki.

\subsubsection{A map of finite spaces implicit in the definition of normality} 
Consider a  ``natural'' 
example of the situation described in Axiom  $(\mathrm O '_\mathrm{V})$:

\noindent\includegraphics[width=\linewidth]{Bourbaki_Urysohn1.png}

As the example is natural, we assume that $X$ is Hausdorff and connected.

Let $A$ and $B$ be any two disjoint closed non-empty subsets of $X$. As we are assuming that $X$ is Hausdorff and connected, 
contracting in $X$ to a point each of subsets $A$, and $B$, and $X\setminus(A\cup B)$, 
gives a 3-point space $\vvLambda$ with one open point $x$ and two closed points $a$ and $b$
such that $a,b\in \cl x$. Maps to $\vvLambda$ ``classify'' pairs of disjoint closed subsets:
to give a map $\chi:X\lra \vvLambda$ is the same as to give two disjoint closed subsets $A:=\chi\inv(a)$
and $B:=\chi\inv(b)$ of $X$. 

Similarly, let $U$ and $V$ be two disjoint open subsets such that $A\subset U$ and $B\subset V$.
Contracting in $X$ to a point each of subsets $A$, and $B$, and $U\setminus A$, and $V\setminus B$, and $X\setminus(U\cup U)$, 
gives a 5-point space $\MM$ with two open points $u$ and $v$ and three closed points $a$, $x$, and $b$
whose topology is generated by the open subsets $\{u\}$, $\{v\}$, and $\{a,u\}$, $\{b,v\}$, and $\{a,x,v\}$, 
Similarly, maps to $\MM$ classify pairs of disjoint open subsets separating closed subsets:
to give a map $\chi:X\lra \MM$ is the same as
to give two disjoint closed subsets $A:=\chi\inv(a),B:=\chi\inv(b)$ of $X$ and 
two disjoint open subsets $U:=\chi\inv(\{a,u\}),V:=\chi\inv(\{v,b\})$ 
such that $A\subset U$ and $B\subset V$.

There is a natural map $\MtoL$ such that $X\lra \vvLambda$ is equal to the composition $X\lra \MtoL$.
Namely, it sends  $a$ to $a$, and $b$ to $b$, and the rest $u,x,v$ to $x$. 
This is the same map as defined in \S\ref{def:MtoL:above}.

Finally, we transcribe Axiom $(\mathrm O'_{\mathrm{V}})$ as: 
each map $X\lra \vvLambda$ factors as $X\lra \MM\lra \vvLambda$.

\subsubsection{Reformulating the definition of normality} 
Now we ready to state the theorem. 

\begin{enonce}{Theorem}\label{thmOne} Axiom $(\mathrm O '_\mathrm{V})$ is equivalent to the following: 
\bii\item[$(\mathrm {\tilde O}'_\mathrm{V})$] $\emptyset \lra X \,\rtt\, \MtoL $ 
\eii
\end{enonce}\begin{proof} The proof is a matter of decyphering notation.  
	Let us rewrite the lifting property explicitly as a diagram 
	denoting finite topological spaces by their specialisation preorders:\footnote{Namely, $\overset{\raisebox{1pt}{$\bullet$}}{}{\searrow}\underset{\bblacksquare}{}$ means 
	that $\bblacksquare \,\in\, \cl (\bullet)$.}
\begin{equation}\label{T4Odv}\begin{gathered}
\xymatrix
{  
 {} & 
{
\left\{\underset{\bblacksquare}{}{\swarrow} \overset{\raisebox{1pt}{$ \bullet$}}{}{\searrow} 
\underset{\bblacksquare}{}
{\swarrow}\overset{\raisebox{1pt}{$\bullet$}}{}{\searrow}\underset{\bblacksquare}{}\right\}
}
 \ar[d]^{
}
         \\
{X}\ar@{-->}[ur] |{}  \ar[r]|{} & 
{	 
\left\{\underset{\bblacksquare}{}{\swarrow} \overset{\raisebox{1pt}{$ \bullet$}}{}{\searrow} 
\underset{\bblacksquare}{}\right\}
}
}\end{gathered}\end{equation}
	The horizontal arrow $X\lra \vvLambda$ represents two disjoint closed subsets (namely, the preimages 
	of the two closed points in $\vvLambda$).

	A commutative triangle represents their disjoint neighbourhoods
	(namely, the preimages of $\underset{\bblacksquare}{}{\swarrow} \overset{\raisebox{1pt}{$ \bullet$}}{}{}$
	and $\overset{\raisebox{1pt}{$\bullet$}}{}{\searrow}\underset{\bblacksquare}{}$ under $X\lra \MM$). 
	To see this, note that the topology on $\MM$ is the coarsest topology such that
        $\underset{\bblacksquare}{}{\swarrow} \overset{\raisebox{1pt}{$ \bullet$}}{}{}$ is an open neighbourhood of the closed point $\bblacksquare$ it contains, 
	and, similarly, 
         $\overset{\raisebox{1pt}{$\bullet$}}{}{\searrow}\underset{\bblacksquare}{}$ is an open neighbourhood of the closed point $\bblacksquare$ it contains.

The diagram below is an informal representation of this argument using the subscripts to indicate the (indented) preimages
in terms of Axiom $(\mathrm O '_\mathrm{V})$:
\begin{equation}\label{T4OdvExtendedLabels}\begin{gathered}
\xymatrix
{  
 {} & 
{
\left\{\underset{\bblacksquare_A}{}{\swarrow} \overset{\raisebox{1pt}{$ \bullet_{U\setminus A}$}}{}{\searrow} 
\underset{\bblacksquare_{X\setminus(U\cup V)}}{}
{\swarrow}\overset{\raisebox{1pt}{$\bullet_{V\setminus B}$}}{}{\searrow}\underset{\bblacksquare_B}{}\right\}
}
 \ar[d]^{
	 \bblacksquare_A\longmapsto \bblacksquare_A, \bblacksquare_B\longmapsto \bblacksquare_B 
	\atop  \hphantom{..............}\bullet_{U\setminus A},
\bblacksquare_{X\setminus(U\cup V)}, 
\bullet_{X\setminus V} 
\longmapsto \bullet_{X\setminus(A\cup B)}
}
         \\
{X}\ar@{-->}[ur] |{}  \ar[r]|{} & 
{	 
\left\{\underset{\bblacksquare_A}{}{\swarrow} \overset{\raisebox{1pt}{$ \bullet_{X\setminus (A\cup B)}$}}{}{\searrow} 
\underset{\bblacksquare_B}{}\right\}
}
}\end{gathered}\end{equation}
\end{proof}

We describe how to ``read off'' the map $\MtoL$ from Axiom~$(\mathrm O'_{\mathrm V})$ in this verbose end-note.\endnote{
\subsection{Formalising ($O'_V$)---the definition of normality without the auxiliary set of real numbers.}
Let us try to formalise $(O'_V)$. 

\subsubsection{Representing two disjoint closed subsets}

A concise way to represent {\tt any two disjoint closed subsets $A$ and $B$ of $X$} is
by a {\em characteristic} map from the space $X$ to the finite topological space with two points closed and one point open. 
Denote this 3-point space by $\vvLambda$, and 
the map by $\chi_{A,B}: X\to \vvLambda$.

\subsubsection{Representing two  closed subsets and their disjoint open neighbourhoods}
Similarly, 
a concise way to represent {\sf two disjoint open sets $U,V$ such that $A\subset U$ and $B\subset V$} %
is
by a map from the space to the finite topological space whose points (rather, their preimages) 
represent the subsets involved and whose topology represents
the open/closedness requirements on the subsets.

What is that space ? Such a system of subsets determines a decomposition of $X$ as
a union of disjoint subsets $A, U\setminus A, X\setminus (U\cup V), V\setminus B, B$.
Hence, the characteristic map 
		goes to a space with 5 points we denote by 
		$$\bblacksquare_A, \bullet_{U\setminus A}, 
		\bblacksquare_{X\setminus(U\cup V)}, 
		\bullet_{V\setminus B}, \bblacksquare_B $$
The topology on this 5 point space represents the requirements on the topology of $A,B, U, V$:
namely, equip the space with the coarsest topology such that  points $\bblacksquare_A$ and $\bblacksquare_B$ 
closed, and subsets  $\{\bblacksquare_A, \bullet_{U\setminus A}\}$ 
and $\{\bblacksquare_B, \bullet_{V\setminus B}\}$ are open. 
Denote this space by $\MM$. Evidently, then, to give 
a continuous map $\chi:X\lra \MM$ 
is the same as 
to give 
two closed subsets $A$ and $B$ (namely, preimages of 
the closed points $\bblacksquare_A$ and $\bblacksquare_B$) 
with disjoint neighbourhoods $U$ and $V$ (namely, preimages of 
disjoint open neighbourhoods  $\{\bblacksquare_A, \bullet_{U\setminus A}\}$
and $\{\bblacksquare_B, \bullet_{V\setminus B}\}$ of closed points $\bblacksquare_A$ and $\bblacksquare_B$, resp.).

\subsubsection{...} To finish reformulating $(O'_V)$ we need to somehow relate 
the? maps $X\lra \vvLambda$ and $X\lra \MM$. 
The first map corresponds to the closed disjoint subsets $A,B\subset X$, 
and the second map corresponds to their disjoint open neighbourhoods $U\supset A$ and $V\subset B$. 
The following relationship between these subsets
$$X\setminus(A\cup B)=(U\setminus A)\cup (X\setminus(U\cup V)) \cup (V\setminus B)$$ 
gives an obvious ``forgetful'' map
$\MtoL$ which glues together the $\VV$ in $\MM$ to the open point of $\vvLambda$, 
and a commutative triangle
$$\xymatrix@C=+2.39cm{ & \MM \ar[d] \\
X \ar[ru]|\chi \ar[r]|{\chi_{A,B}} & \vvLambda }$$


Further, as $U$ is a neighbourhood of $A$, make $\{\bblacksquare_A, \bullet_{U\setminus A}\}$ 
open, and for the same reason make $\{\bblacksquare_B, \bullet_{V\setminus B}\}$ open. 
The topology is the coarsest topology such that the point corresponding 

namely, make points $\bblacksquare_A, \bblacksquare_{X\setminus(U\cup V)}, \bblacksquare_B$ 
to be closed, and points $\bullet_U$ and $\bullet_V$ open. 
Further, as $U$ is a neighbourhood of $A$, make $\{\bblacksquare_A, \bullet_{U\setminus A}\}$ 
open, and for the same reason make $\{\bblacksquare_B, \bullet_{V\setminus B}\}$ open. 
The topology is the coarsest topology such that the point corresponding 

The sets $A\cap V=\emptyset$, $B\cap U=\emptyset$,

In this case, it has 
two open points and three closed points:
the disjoint closed subsets, and the completion of the open neighbourhoods are sent to the closed points,
and everything in a neighbourhood but not in a closed subset, is sent to the two open points.
(Todo: maybe skip this sentence, such details are given in next paragraph. maybe explain clearer
how to represent systems of subsets of a space as maps to finite spaces. Maybe: define
{\em characteristic function $\chi_{A,B,..}$} of a {\em system} of subsets with {\em requirements}
on their topology (whether open, closed, etc) ?)

How to represent a finite topological space ? The topology of a finite topological space 
is fully determined by ts specialisation preorder $o\leqslant c$ iff $c \in \cl\, o$: 
a closed subset is equal to the union of closures of its points,
and an arbitrary (necessarily finite) union of the closures of points
in necessary closed. Hence,  a subset is closed iff it it is upward closed 
in the specialisation preorder.  

For the same reason, a map is continuous iff it preserves the specialisation preorder.

We denote the specialisation preorder by arrow $o \lra c $ iff $c \in \cl\,o$ because
we like to think of a finite topological space as a (degenerate) {\em category where all diagrams commute},
or, equivalently, there is at most one morphism from any one object to another. 
We like to think of a finite topological space as a category because we are already working with a category 
(of topological spaces), and that makes sense from programming point of view.

\subsubsection{
	Now we repeat the same with some notation.} 
A concise way to represent {\tt $A$ and $B$ ... two disjoint closed subsets of $X$} is
by {\em the characteristic} map from the space to the finite topological space with two points closed and one point open; 
see Fig.(which best to put here)?; squares denote closed points, and bullets open points.
$$X\lra 
\left\{\underset{\bblacksquare_A}{}{\swarrow} \overset{\raisebox{1pt}{$ \bullet_{X\setminus (A\cup B)}$}}{}{\searrow} 
\underset{\bblacksquare_{ B}}{}\right\}
$$
$$A\longmapsto \bblacksquare_A, B\longmapsto \bblacksquare_B, \, X\setminus(A\cup B)\longmapsto \bullet_{X\setminus(A\cup B)}
$$ 

A concise way to represent {\tt 
two disjoint open sets $U$, $V$ such that $ A \subset U$ and $B\subset V$} for  two disjoint closed subsets $A$ and $B$
is 
by an arrow (map) from $X$ to a finite topological space with 5 points, two open and three closed
$$X\lra 
\left\{\underset{\bblacksquare_A}{}{\swarrow} \overset{\raisebox{1pt}{$ \bullet_{U\setminus A}$}}{}{\searrow} 
\underset{\bblacksquare_{X\setminus(U\cup V)}}{}
{\swarrow}\overset{\raisebox{1pt}{$\bullet_{V\setminus B}$}}{}{\searrow}\underset{\bblacksquare_B}{}\right\}$$ 
 $$A\longmapsto \bblacksquare_A, B\longmapsto \bblacksquare_B,  X\setminus(U\cup V)\longmapsto \bblacksquare_{X\setminus(U\cup V)},$$
$$U\setminus A \longmapsto \bullet_{X\setminus U},
V\setminus B \longmapsto \bullet_{X\setminus V}$$

To represent that both arrows talk about the same $A$ and $B$, say that the diagram commutes:
\begin{equation}\label{T4Odv}\begin{gathered}
\xymatrix
{  
 {} & 
{
\left\{\underset{\bblacksquare_A}{}{\swarrow} \overset{\raisebox{1pt}{$ \bullet_{U\setminus A}$}}{}{\searrow} 
\underset{\bblacksquare_{X\setminus(U\cup V)}}{}
{\swarrow}\overset{\raisebox{1pt}{$\bullet_{V\setminus B}$}}{}{\searrow}\underset{\bblacksquare_B}{}\right\}
}
 \ar[d]^{
	 \bblacksquare_A\longmapsto \bblacksquare_A, \bblacksquare_B\longmapsto \bblacksquare_B 
\atop  \bullet_{U\setminus A},
\bblacksquare_{X\setminus(U\cup V)}, 
\bullet_{X\setminus V} 
\longmapsto \bullet_{X\setminus(A\cup B)}
}
         \\
{X}\ar@{-->}[ur] |{}  \ar[r]|{} & 
{	 
\left\{\underset{\bblacksquare_A}{}{\swarrow} \overset{\raisebox{1pt}{$ \bullet_{X\setminus (A\cup B)}$}}{}{\searrow} 
\underset{\bblacksquare_B}{}\right\}
}
}\end{gathered}\end{equation}

Note that $(O''_V)$ is represented by the same diagram involving the same maps, only labelled differently:
\begin{equation}\label{T4Oddv}\begin{gathered}
\xymatrix
{  
 {} & 
{
\left\{\underset{\bblacksquare_A}{}{\swarrow} \overset{\raisebox{1pt}{$ \bullet_{W\setminus A}$}}{}{\searrow} 
\underset{\bblacksquare_{\bar W\setminus W}}{}
	{\swarrow}\overset{\raisebox{1pt}{$\bullet_{V\setminus \bar W}$}}{}{\searrow}\underset{\bblacksquare_{X\setminus V}}{}\right\}
}
 \ar[d]^{
	 \bblacksquare_A\longmapsto \bblacksquare_A, \bblacksquare_{X\setminus V}\longmapsto \bblacksquare_{X\setminus V} 
\atop  \bullet_{W\setminus A},
\bblacksquare_{\bar W\setminus W}, 
\bullet_{V\setminus \bar W} 
\longmapsto \bullet_{V}
}
         \\
{X}\ar@{-->}[ur] |{}  \ar[r]|{} & 
{	 
\left\{\underset{\bblacksquare_A}{}{\swarrow} \overset{\raisebox{1pt}{$ \bullet_{V}$}}{}{\searrow} 
	\underset{\bblacksquare_{X\setminus V}}{}\right\}
}
}\end{gathered}\end{equation}

\begin{theorem}[(Urysohn)] 
	Diagrams (\ref{T4Ov}) and (\ref{T4Odv}) are equivalent.
\end{theorem}

}

\subsubsection{The map  $\protect\MtoL$  is a trivial Serre fibration}
\begin{enonce}{Lemma}\label{lem:rttMtoL}\label{lem:rttMtoLLtoO}
\label{lem:T4Odddv}  If both $A$ and $X$ have closed points, and $X$ is
completely (=hereditary)  normal \cite[IX\S4, Exercise 3]{Bourbaki} 
	(i.e. separation axioms $T_0$ and $T_6$), then 
		$$
\begin{gathered}
\vcenter{\hbox{\xymatrix@-2pc{ 
	A \xra f X\text{ is a closed subset iff } 
}} } 
\vcenter{\hbox{\xymatrix
{  
 A \ar[d]|f\ar[r] & \vvM  \ar[d]
	\\ X \ar[r] \ar@{-->}[ur] & \vvLambda} }} 
\end{gathered}$$
In particular, this equivalence holds whenever both $A$ and $X$ are metrisable. 
\end{enonce}
\begin{proof} $\implies$: 
The commutative square gives two disjoint closed subsets $B,C\subset X$, and two disjoint 
open subsets $U, V\subset A\subset X$ open in $A$ 
such that $B\subset U\cap A$ and $C\subset V \cap A$. 
The lifting arrow $X\lra \MM$ gives two disjoint open subsets $U'\supset U\cup B$ and $V'\supset V\cup C$.
$U$ and $V$ being open in $A$ means there are some open subsets $U'',V''\subset X$ 
such that $U=U''\cap A$ and $V=V''\cap B$. 
The open subset $U'' \cup (X\setminus (A\cup C))$ separates $U\cup B$ from $V\cup C$, and 
the open subset $V'' \cup (X\setminus (A\cup C))$ separates $V\cup C$ from $U\cup B$. 
A characterisation of hereditary normality says that a space is hereditary normal
iff any two (not necessarily closed) separated subsets are separated by their open disjoint neighbourhoods.\footnote{%
	It is easy to see that metric spaces have this property. Indeed, the assumption on the separated subsets means that
	each point $x$ of the separated subsets is at positive distance from the other subset; hence, define 
	the open neighbourhood of a subset by adding a ball around each $x$ of radius half the distance to the other subset.
	The triangle inequality implies these balls do not intersect.}
Take $U'$ and $V'$ to be these disjoint neighbourhoods of $U\cup B$ and $V\cup C$.

The following diagram somewhat informally represents this argument; there we use sup- and subscripts
to indicate intended preimages.  
\begin{equation}\label{T4Odv}
\begin{gathered}
\xymatrix
{   
	A \ar[d] \ar[r]  & 
{
	\left\{\underset{\bblacksquare^{B\cap A}_B}{}{\swarrow} \overset{\raisebox{1pt}{$ \bullet^{U\setminus (B\cap A)}_{U'\setminus B}$}}{}{\searrow} 
	\underset{\bblacksquare^{A\setminus(U\cup V)}_{X\setminus(U'\cup V')}}{}
{\swarrow}\overset{\raisebox{1pt}{$\bullet^{V\setminus(C\cap A)}_{V'\setminus C}$}}{}{\searrow}\underset{\bblacksquare^{C\cap A}_C}{}\right\}
} 
 \ar[d]^{
}
         \\
{X}\ar@{-->}[ur] |{}  \ar[r]|{} & 
{	 
	\left\{
	\underset{\bblacksquare_B}{}{\swarrow} 
	\overset{\raisebox{1pt}{ \ensuremath{ \bullet_{X\setminus (B\cup C)} } }}{}{\searrow} 
        \underset{\bblacksquare_C}{} 
	\right\}
	}  
}	\iff 
{ \substack {
	{\displaystyle \text{subsets }U\cup B \text{ and } V\cup C \text{ of } X } \\
\displaystyle \text{ are separated } \\ 
\displaystyle \text{by disjoint open neighbourhoods }  
      \\ \displaystyle U' \text{ and }V' }  
 }    
 \end{gathered}  \end{equation} 

\noindent$\impliedby$: 
First show that that the image of $A$ is closed by showing that for each closed point $b\in X\setminus A$ it holds that 
$\cl_B \Imm A \subset X\setminus b$. For this consider the commutative square where $A$ is sent to the closed mid-point of $\vvM$, 
the set $X\setminus b$ is sent to the open top-point of $\vvLambda$, and $b$ is sent to one of the closed points of $\vvLambda$. 
The lifting $X\lra \vvM$ has to send the closure $\cl_B \Imm A$ into the mid-point of $A$, yet $b$ is sent elsewhere. 
Hence $b\notin \cl_B \Imm A$, as required.

	Because points of $A$ are closed, it is enough to show that the topology on $A$ is induced from $X$, i.e.{} 
	each closed subset $Z$ of $A$
	is the preimage of some closed subset $Z'$ of $X$. To see this, consider a commutative square where $Z\subset A$ 
	is the preimage of the mid-point of $\vvM$ under the arrow $A\lra \vvM$. Then the preimage of the same point 
	under the lifting $X\lra \vvM$ has the required property. 

The two arguments above are represented by the following two diagrams below.
\begin{equation}\label{T4Odv}\begin{gathered}
\xymatrix
{  
	A \ar[d] \ar[r]  & 
{
	\left\{\underset{\bblacksquare^\emptyset_b}{}{\swarrow} \overset{\raisebox{1pt}{$ \bullet^{\emptyset}_{?}$}}{}{\searrow} 
	\underset{\bblacksquare^A_{Z\supset\cl_B\Imm A}}{}
{\swarrow}\overset{\raisebox{1pt}{$\bullet^\emptyset_{?}$}}{}{\searrow}\underset{\bblacksquare^\emptyset_\emptyset}{}\right\}
}
 \ar[d]^{
}
         \\
{X}\ar@{-->}[ur] |{}  \ar[r]|{} & 
{	 
	\left\{\underset{\bblacksquare_b}{}{\swarrow} \overset{\raisebox{1pt}{$ \bullet_{X\setminus b}$}}{}{\searrow} 
\underset{\bblacksquare_\emptyset}{}\right\}
}
}\implies { {\text{for each closed point }b\in X\setminus \Imm A}\atop{\cl_B\Imm A\subset X\setminus b} } \end{gathered}\end{equation}
\begin{equation}\label{T4Odv}\begin{gathered}
\xymatrix
{  
	A \ar[d] \ar[r]  {} & 
{
	\left\{\underset{\bblacksquare^\emptyset_\emptyset}{}{\swarrow} \overset{\raisebox{1pt}{$ \bullet^{A\setminus Z}_{?}$}}{}{\searrow} 
	\underset{\bblacksquare^{Z}_{Z'}}{}
{\swarrow}\overset{\raisebox{1pt}{$\bullet^\emptyset_{?}$}}{}{\searrow}\underset{\bblacksquare^\emptyset_\emptyset}{}\right\}
}
 \ar[d]^{
	 \bblacksquare_A\longmapsto \bblacksquare_A, \bblacksquare_B\longmapsto \bblacksquare_\emptyset
\atop  \bullet_{U\setminus A},
\bblacksquare_{X\setminus(U\cup V)}, 
\bullet_{X\setminus V} 
\longmapsto \bullet_{X\setminus(A\cup B)}
}
         \\
{X}\ar@{-->}[ur] |{}  \ar[r]|{} & 
{	 
\left\{\underset{\bblacksquare_\emptyset}{}{\swarrow} \overset{\raisebox{1pt}{$ \bullet_{X}$}}{}{\searrow} 
\underset{\bblacksquare_\emptyset}{}\right\}
}
}\implies {Z\text{ is the preimage of } Z' }\end{gathered}\end{equation}
\end{proof}

\subsubsection{Defining cofibrations and trivial fibrations?} 
\begin{coro} The maps $\MtoL$ and $\vvLambda\to\ptt$ are trivial Serre fibrations. In other words, 
	 $\MtooLVtoolr$ 
	is a class of trivial Serre fibrations. 
\end{coro}
\begin{proof} The lifting property $\mathbb S^n \lra \mathbb D^{n+1} \,\rtt\, \MtoL$ defining trivial Serre fibration involves only metric spaces,
	hence the previous Lemma applies. Verifying the claim for $\vvLambda \lra \ptt$ is an exercise. 
\end{proof}
The corollary leads to the question what this is the class of {\em all} trivial fibrations which we state in Conjecture~\ref{conj:CW-WC}
among a few other conjectures in \S\ref{conjectures:WC}.

\begin{coro}\label{coro:ANR-C} If both $A$ and $X$ are metrisable absolute neighbourhood retracts, then 
for a map	$f:A\lra X$ is a closed Hurewitz cofibration the following are equivalent:
\bee
\item $f:A\lra X$ is a closed inclusion
\item $f:A\lra X$ is a closed Hurewitz cofibration 
\item	
	$$ { {A} \atop {{ \,\,\,\downarrow \,f} \atop \displaystyle X } }
\,\in\,\MtooLl $$
\eee\end{coro}
\begin{proof}
By 
		        \cite[Proposition $4$.3]{NlabCF} and references therein, 
for absolute neighbourhood retracts $A$ and $X$,
 being a closed Hurewitz cofibration is equivalent  to
	being a closed subspace inclusion $A\lra X$, and by Lemma~\ref{lem:rttMtoL} 
	this is equivalent to being in $\Bigl\{\MtoL\Bigr\}^{l}$.
\end{proof}

\begin{coro} 
In the full subcategory $\Topp_{\text{\rm ANR}}$ 
	of topological spaces consisting of finite topological spaces 
	and metrisable absolute neighbourhood retracts, 
	$$\MtooL^{lr}_{\mathrm{ANR}}$$
	is a class of trivial Serre fibrations 
	containing all trivial Hurewitz fibrations of metrisable absolute neighbourhood retracts.
\end{coro}
\begin{proof} 
	Note that a map from a finite topological space to  a Hausdorff space is necessarily constant on the connected components.
	Hence, if we calculate $\Bigl\{\MtoL\Bigr\}^{lr}$ in the full subcategory of $\Topp$ consisting of finite topological spaces
	and metric absolute neighbourhood retracts, we may ignore maps to and from finite spaces in 
	$\Bigl\{\MtoL\Bigr\}^{lr}$, hence Corollary~\ref{coro:ANR-C} implies  that $\Bigl\{\MtoL\Bigr\}^{lr}$ contains all trivial Hurewitz fibrations 
	of metric absolute neighbourhood retracts. 
\end{proof}
\begin{rema}
	Michael continuous selection theory, cf.~\cite[\S$3.1$,\S$3.2$]{V}, can probably be used 
	to show that a rather large class of locally constant maps with contractible fibres lies 
	in $\Bigl\{\MtoL\Bigr\}^{lr}$ calculated in another full subcategory of $\Topp$. 
\end{rema}

\subsubsection{A definition of trivial fibrations ?}\label{conjectures:WC} 
The corollary above leads to the question whether this is the class of {\em all} trivial fibrations. 
Here we state it in its simplest form, see \S~\ref{conjectures:WC:Lw} for more. 
\begin{conj}\label{conj:CW-WC} 
	A cellular map $f$ 
	of finite CW complexes is a trivial Serre fibration iff
	$$f\,\in\,\MtooLVtoolr$$ 
%
\end{conj}

\begin{enonce}{Problem} Find 
	a model structure on the category of topological spaces where 
	$$\MtooLVtoolr\text{ is the class of trivial fibrations. }$$
\end{enonce}

\subsubsection{Replacing the interval object $\bbbR$ by $\protect\vvLambda_\omega$ ?}\label{conjectures:WC:Lw} 
Remark~\ref{rema:perfectly-normal} suggests it might be natural to consider $\vvLambda_\omega$ instead of $\bbbR$ as a path/interval object,
and modify the notion of geometric realisation in the spirit of \cite[Corollary 4.13,p.$25$]{Wofsey}.\footnote{
\cite{Wofsey} uses different conventions: for us their Hasse diagrams of preorders are upside down, e.g. there closed points are on top, and downsets (=downward closed)  there are open rather than closed.}

Somewhat modifying  the definition of $\tilde X$ mentioned after \cite[Corollary 4.13,p.$25$]{Wofsey}, define 
{\em the non-$T_0$ geometric realisation $|Y|_\omega$} 
of a finite topological space $X$ by
$$ |Y|_\omega := \mylim ( \cdots \lra \beta^{n+1} X \lra  \beta^n X \lra \cdots \lra X )$$
where $\beta^{n+1}(X):=\beta(\beta^n(X))$, and $\beta(X)$ is the set of non-empty  
chains(=linearly ordered subsets) of $X$
ordered under inclusion (i.e.~$C_1\rightarrow  C_2$ iff $C_1\supset C_2$) , 
and the map $\beta(X)\lra X$ takes a chain into its least element.\footnote{Instead of $\beta(X)$ 
one might use iteratively 
$\Hom(\left\{\overset{\raisebox{1pt}{$\bullet$}}{}{\searrow}\underset{\bblacksquare}{}\right\}, X)$ 
as the image of  $\left\{\overset{\raisebox{1pt}{$\bullet$}}{}{\searrow}\underset{\bblacksquare}{}\right\}$ 
is a non-empty chain.}
\begin{conj}
%
For each finite topological space $Y$ it holds
	$${|\Delta Y|_\omega \atop {\downarrow\atop{\displaystyle Y}} } \,\in\,\MtooLVtoolr$$ 
\end{conj}
\begin{proof}[Evidence] This is a direct analogue of  \cite[Theorem $3.5$]{Wofsey} 
where we drop the assumption that $X$ is {\em perfectly} normal. 
In terms of the lifting property \cite[Theorem $3.5$]{Wofsey}  states that for a finite space $Y$ and 
a closed subspace $A$ of a perfectly normal space $X$ 
it holds
$A\lra X \,\rtt\, |\Delta Y| \xra \pi  Y$, and that any liftings are 
``homotopic such that every stage of the homotopy is also such a lift''.
Similarly to the Urysohn theorem,  Lemma~\ref{lem:rttMtoL} should let us replace 
the condition ``$A$ is a closed subset of $X$'' for perfectly normal $X$ 
by the appropriate lifting property.  \end{proof}
\begin{conj}
	For a map $f:X\lra Y$ of finite topological spaces, its geometric realisation 
	$|f|:|\Delta X|\lra |\Delta Y|$ 
	is a trivial Serre fibration iff 
	$$ { {|X|_\omega} \atop {{ \,\,\,\,\,\,\downarrow \,|f|_\omega} \atop \displaystyle |Y|_\omega } }
\,\in\,\MtooLVtoolr $$
\end{conj}
\begin{proof} This is an analogue of Conjecture~\ref{conj:CW-WC} where we try to consider $\vvLambda_\omega$ as a path/interval object rather than $\bbbR$. 
\end{proof} 
%
%

\subsection{Extension of a continuous real-valued function}
Now we try to reformulate Theorem 2 of Urysohn \cite[IX\S$4.2$]{Bourbaki}, quoted below.
\newline\noindent\includegraphics[line]{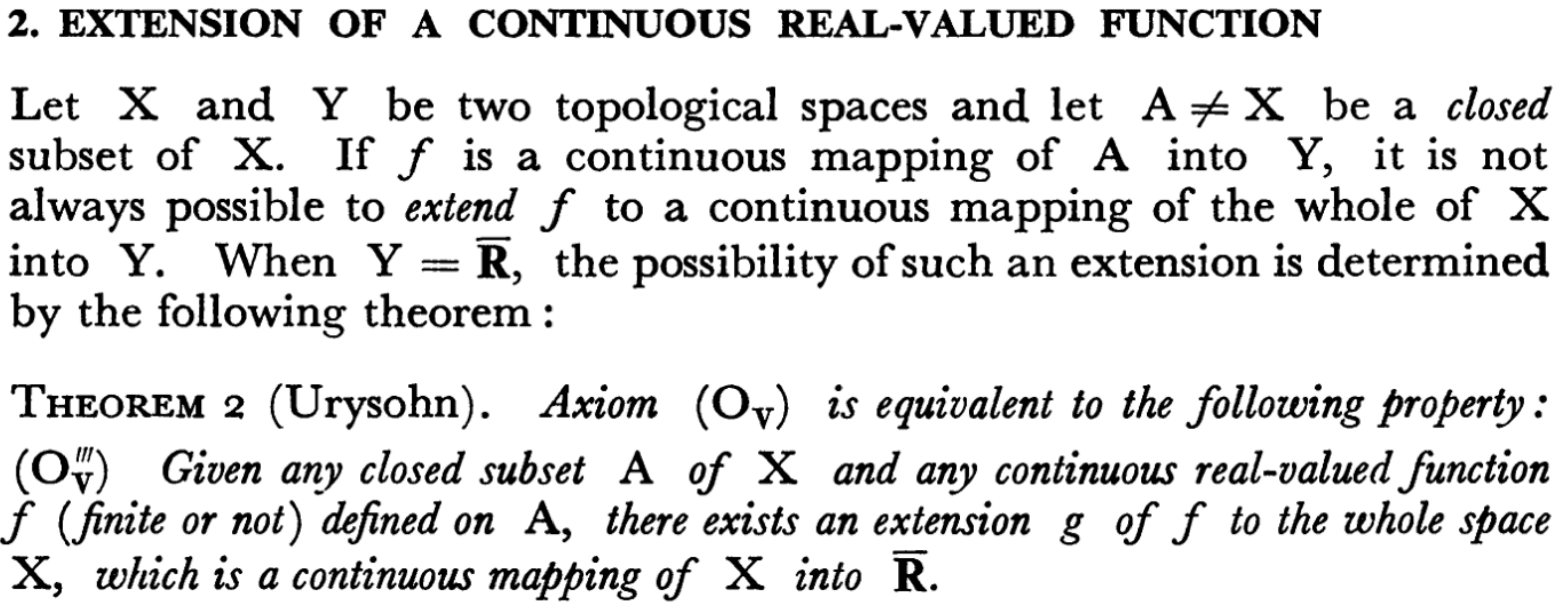}
Here in Bourbaki $\overline {\Bbb R}:=\{-\infty\}\cup \Bbb R \cup \{+\infty\}$ with the obvious uniform structure;
as a topological space, $\overline {\Bbb R}=[-1,1]$. 

First let us reformulate  $(O'''_V)$ as a lifting property:\footnote{
	A verification (\cite[\S$3.4$.4]{MR}) shows that 
being a closed subset is also a lifting property, hence we may equivalence of Axioms $( O''_V)$ and $(O'''_V)$ entirely in terms of maps of finite spaces,
as follows: \newline For each space $X$ the following are equivalent:
\bi\item[$(\bar O'_V)$ ] $\emptyset\lra X \,\rtt\, \MtoL $ 
\item[$(\bar{\bar O}'''_V)$] for each map $A\lra X$
the following implication holds:
\newline $A\lra X \,\rtt\, \left\{ \overset{\raisebox{1pt}{$ \star_u\llrra \star_v$}}{}{\searrow} \underset{\bblacksquare_a}{}\right\} \lra \left\{ \star_u \llrra \star_v=\bblacksquare_a \right\}$ 
implies $A\lra X \,\rtt\, \barR \lra \{\bullet\}$
\ei}
\begin{itemize} 
\item[$(\bar O'''_V)$] 
		Given any closed subset $A$ of $X$, it holds 
$$
\xymatrix@C=+5.39cm
{  
A  \ar[d]|{(\text{closed subset})}\ar[r] &  { \overline{ \Bbb  R}} 
\\ X  \ar@{-->}[ur]  &  }$$
\end{itemize}

\def\bbbR{\Bbb R}
\def\barR{\overline{\Bbb R}}
\def\barRR{[-1,1]} 

\subsubsection{Reformulating the proof of Theorem 2} \cite{Bourbaki} construct $g:=\lim_n g_n$ as the limit of a sequence of functions $g_n:X\lra [-1,1] $
such that for each $x\in X$ $|g_m(x)-g_n(x)|\leqslant (2/3)^{\min(m,n)+1}$ and thereby $|g(x)-g(n)|\leqslant (2/3)^n$; 
at the inductive step $g_{n+1}$ is constructed using Axiom 
$(\mathrm O'_{\mathrm V})$ (\cite[II\S$4.2$,Lemma $1$]{Bourbaki})
from $f$ and $g_{n}$.
Thus, only approximate values of $g_n$'s are important, and in terms of  finite spaces we may prefer to think of 
of the $g_n$'s as a sequence of functions $X\lra \vvLambda_n$  to finite topological spaces representing 
finer and finer decomposition's of $[-1,1]$. 
Subdividing in two an interval in a decomposition corresponds to a map of finite topological spaces 
which is a pullback of $\MtoL$ corresponding to the subdivision of the interval. 
As a right lifting property is preserved under pullbacks, we may use the lifting property $A\lra X \,\rtt\, \MtoL$ to 
construct $g_{n+1}$ from $g_n$.

Reformulating this sketch of a proof in a diagram-chasing manner leads to the following theorem.

\begin{enonce}{Theorem}[Urysohn]\label{thm:Thm2} 
	$${\barRR \atop {\downarrow \atop {\{\bullet\}}}} \,\in\, \MtooLVtoolr$$
\end{enonce}

\begin{proof} 
%
	Decompose the interval into smaller intervals as 
	$$[-1,1]=\{-1\}\cup (-1,a_1) \cup \{a_1\}\cup ... \cup \{a_n\}\cup (a_n,1)\cup \{1\}$$
	and subdivide an open interval $(a_m,a_{m+1})$ in two halves:
	$$[-1,1]=\{-1\}\cup (-1,a_1) \cup \{a_1\}\cup ... \cup \{a_m\} \cup (a_m, a'_m) \cup \{a'_m\}\cup (a'_m,a_{m+1})\cup \{a_m\} \cup ... \cup \{a_n\}\cup (a_n,1)\cup \{1\}$$
Contracting each open subinterval in these decompositions gives maps to finite spaces 
	$\barRR \lra \vvLambda_n$ and $\barRR\lra \vvLambda_{n+1}$.
	Note that $\vvLambda_1=\vvLambda$ and $\vvLambda_2=\MM$. 
Contracting closed subintervals $[-1,a_m]$ and $[a_{m+1},1]$ gives maps to finite spaces 
$\barRR \lra \vvLambda$ and $\barRR\lra \MM$. 
These maps fit into a commutative diagram with a pull-back square on the right: 
	$$\xymatrix@C=2.39cm{
		[-1,1] \ar[r] \ar[rd] \ar@/^1pc/[rr]  \ar@/_0.3pc/[rrd] & \vvLambda_{n+1} \ar[r]  \ar[d]^{\therefore\MtooLlr} & \MM \ar[d]^{\MtooLlr}  \\
	                               & \vvLambda_n  \ar[r]  &           \vvLambda 
	}$$
As $r$-orthogonals are closed under pullbacks we have that 
	$\vvLambda_{n+1}\lra \vvLambda_n$ is in $\MtooLlr$.
Taking finer and finer decompositions we get a sequence of maps in $\MtooLlr$
	$$ ...\lra \vvLambda_{n+1} \lra \vvLambda_n \lra ... \lra \vvLambda $$ 
As $r$-orthogonals are closed under colimits, the following map is also in the class: 
	$$\vvLambda_\omega:=\mylim_{}  \left( ...\lra \vvLambda_{n+1} \lra \vvLambda_n \lra ... \lra \vvLambda \right)\lra \vvLambda  $$ 
Therefore the composition $\vvLambda_\omega\lra \vvLambda \lra \ptt$ lies in the orthogonal we consider: 
	$${\vvLambda_\omega\!\!\!\atop{\downarrow\atop \bullet}} \,\,\in\,\, \MtooLVtoolr$$
As orthogonals are closed under retracts, it only remains to prove that $[-1,1]$ is a retract of $\vvLambda_\omega$.
The maps $\tau_n: [-1,1]\lra \vvLambda_n$ associated with the decompositions give rise to 
a map $\tau_\omega: [-1,1]\lra \vvLambda_\omega$ by the universal property of colimits. 
Now, a standard argument using ``the idea of uniform convergence'' can be used to construct a continuous function $g:\vvLambda_\omega \lra [-1,1]$ 
	by setting $g(x):= \bigcap_n \cl (\tau_n\inv(\pr_n(x)))$ where $\pr_n:\vvLambda_\omega\lra\vvLambda_n$ is the obvious projection. 
Indeed, by construction we have $\tau_{n+1}\inv(g_{n+1}(x))\subseteq \tau_n\inv(g_n(x))$,
	and we may choose the maps $\tau_n:\barRR\lra \vvLambda_n$ such that the intervals $\tau_n\inv(g_n(x))$ become arbitrarily small. 
	To see that $g:\vvLambda_\omega \lra\barRR$ is continuous, note that for each open neighbourhood $U_t\ni t\in [-1,1]$ there is $n>1$ and a finite 
	open subset  $ U_n\subset \vvLambda_n$ such that $t\in \tau_n\inv(U_n)$ and  $ \cl ( \tau_n\inv(U_n)) \subset U_t$; 
	hence by construction of $g$ 
	we have that  $g(g_n\inv(U_n))\subset U_n \subset U_t$, hence there is an open neighbourhood of $x$ which maps inside $U_t\ni t:=g(x)$.    
%
\end{proof}
\begin{rema}\label{remaGone}\label{rema:perfectly-normal} 
	In general $g\circ \tau_n \neq g_n $; we only have that  $\tau_n(g(x))\in \cl_{\vvLambda_n}(g_n(x))$ for each $x\in X$. 
	In fact, $\emptyset\lra X \,\rtt\, [-1,1]\xra{\tau} \vvLambda$ iff $X$ is perfectly normal 
for the map $[-1,1]\xra{\tau} \vvLambda$
	corresponding to the decomposition $[-1,1]=\{-1\}\cup (-1,1)\cup \{1\}$ 
	(cf.~\cite[IX\S4,Exercise 7a]{Bourbaki}).
	As not each normal space is perfectly normal, we have that  
        $${\barRR \atop {\downarrow \atop {\displaystyle\vvLambda}}} \,\notin\, \MtooLVtoolr
	\text{ }\text{ yet }\text{ } 
	{{\,\,\,\vvLambda_\omega} \atop {\downarrow \atop {\displaystyle\vvLambda}}} \,\in\, \MtooLVtoolr
	$$
\end{rema}

\subsubsection{Corollary to Theorem 2} We now reformulate the corollary, which is the main theorem of this section.
\newline\noindent\includegraphicss[0.99]{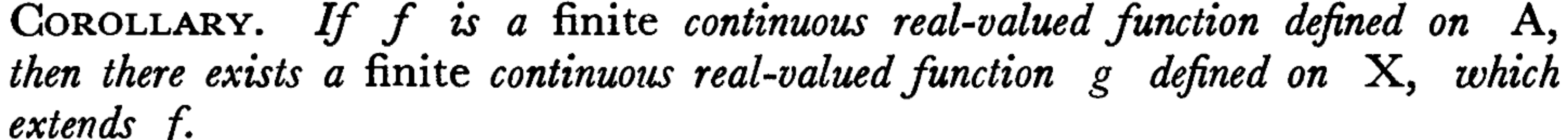}
\begin{enonce}{Corollary}\label{cor:Thm2}
        $${\bbbR \atop {\downarrow \atop {\ptt}}} \,\in\, \MtooLVtoolr$$
\end{enonce}
\begin{proof} We rewrite the proof of Bourbaki in a diagram chasing manner. For us it is more convenient to talk about $(-1,1)$ rather than $\bbbR$. 

	We need to construct an extension $g:X\lra (-1,1) \lra \ptt$ of an arbitrary function $f:A\lra (-1,1)$
	where $A\lra X$ satisfies some assumptions. 

	Up to some technicalities, the Bourbaki defines the extension $g:A\lra (-1,1)$ 
	as $g(x):=\frac{g_1(x)+g_2(x)}2$ 
	where $g_1:X\lra [-1,1]$ is an extension of $f:A\lra(-1,1)$, 
	and $g_2:X\lra [-1,1]$ is an extension of $\tilde g_1:A\cup g_1\inv(-1) \cup g_1\inv(1) \lra (-1,1)$  
	where $\tilde g_{1\,|A}=f$, and $\tilde g_{1\,|\, g_1\inv(-1) \cup g_1\inv(1) }=0$. 
%
%
%
%
%
%
%

We render the construction of $g_2$  as follows. 
	The two diagrams below are the same in different notation; the second diagram uses subscripts to indicate 
what the maps are in a somewhat informal manner. 
	$$\xymatrix@C=1.42cm{
		& (-1,1) \ar[d] \ar[rrr] & 			 &		    &		 \vvLambda_\omega\ar[d]\ar[r] & [-1,1]  \\
	A\ar[ru]|f \ar[r]\ar[d]	& [-1,1] \ar[r] \ar[rd]  & \MM \ar[r]\ar[d] & \VV \ar[r] & \MM \\
	X\ar@{-->}[ru]|{g_1}   \ar@/^2pc/@{-->}[rru] \ar@/^5pc/@{-->}[rrrru]   \ar@/_2pc/@{-->}[rrrruu]|{g_2}          &		& \vvLambda 	     & 
	}$$

	$$\phantom{.}\hskip-12pt\xymatrix@C=1.0039cm{
		(-1,1) \ar[r] \ar@/^2pc/[rrr] & [-1,1]	\ar[d]|(0.4){\pi_\MM}		 &		    &  \vvLambda_\omega \ar[r] \ar[d]|(0.4){\pi'_\MM}  &		 [-1,1] \\
	A\ar[u]|f \ar[r]\ar[d]	&    
			 {
				 \left\{\underset{\bblacksquare^a_{-1}}{}{\swarrow} \overset{\raisebox{1pt}{$ \bullet^u$}}{}{\searrow}
				 \underset{\bblacksquare^x}{}
				 {\swarrow}\overset{\raisebox{1pt}{$\bullet^v$}}{}{\searrow}\underset{\bblacksquare^b_1}{}\right\}
			 }
	\ar[r]\ar[d] & 
			 {
				 \left\{\overset{\raisebox{1pt}{$ \bullet^{u}$}}{}{\searrow}
				 \underset{\bblacksquare^{a=x=b}_{-1,1}}{}
				 {\swarrow}\overset{\raisebox{1pt}{$\bullet^{v}$}}{}\right\}
			 }
	\ar[r] & 
			 {
				 \left\{\underset{\bblacksquare}{}{\swarrow} \overset{\raisebox{1pt}{$ \bullet^{u}$}}{}{\searrow}
				 \underset{\bblacksquare^{a=x=b}_{-1,1}}{}
				 {\swarrow}\overset{\raisebox{1pt}{$\bullet^{v}$}}{}{\searrow}\underset{\bblacksquare}{}\right\}
			 }
	\\
	X\ar@/^2pc/@{-->}[ruu]|(0.6){g_1}  \ar@/^2pc/@{-->}[ru]|(0.35){g_1\circ\pi_\MM}    
	\ar@/^4pc/@{-->}[rrru]|(0.2){g_1\circ \pi_\MM\circ (\MM\to\VV\to \MM)=g'_2\circ \pi'_\MM}  \ar@/_2pc/@{-->}[rrruu]|{g'_2}  \ar@/_2pc/@{-->}[rrrruu]|{g_2}          &	
				   {                                                                                                                                                                                
					   \left\{\underset{\bblacksquare^a_{-1}}{}{\swarrow} \overset{\raisebox{1pt}{$ \bullet^{u=x=v}_{(-1,1)}$}}{}{\searrow}
					   \underset{\bblacksquare^b_{1}}{}\right\}
					   } &
	}$$
Construct $g_1:X\lra [-1,1]$ by $A\lra X \,\rtt\, [-1,1]\lra \ptt$ and get $\bar g_1: X\lra \MM$ 
as the composition $X \xra{g_1} [-1,1]\lra \MM$.
Then take the non-trivial endomorphism $\MM\lra \VV \lra \MM$; this gives another map $\bar g_2': X\lra \MM$ 
	such that $\bar g_1^{-1}(\bblacksquare^a)\cup \bar g_1^{-1}(\bblacksquare^b) \subset \bar {g'}_2\inv(\bblacksquare^{a=x=b})$. 
Construct $g'_2:X\lra \vvLambda_\omega$ by $ A\lra X \,\rtt\, \vvLambda_\omega \lra \MM$ and 
	let $g_2:X \lra [-1,1]$ be the composition $X \xra {g'_2} \vvLambda_\omega \lra [-1,1]$. Then by construction 
	$g_2:X\lra [-1,1]$ is an extension of $f:A\lra (-1,1)$, and so is $g(x):=\frac{g_1(x)+g_2(x)}2$.
	The inclusion $\bar g_1\inv(\bblacksquare^a)\cup \bar g_1\inv(\bblacksquare^b) \subset \bar {g'}_2\inv(\bblacksquare^{a=x=b})$
	implies $$ g_1\inv(-1)\cup \bar g_1\inv(1) =
	\bar g_1\inv(\bblacksquare^a)\cup \bar g_1\inv(\bblacksquare^b) \subset \bar {g'}_2\inv(\bblacksquare^{a=x=b})\subset g_2\inv((-1,1)).$$
	In other words, $-1<g_2(x)<1$ whenever  $g_1(x)=-1$ or $g_1(x)=1$ and therefore 
	$-1<g(x):=\frac{g_1(x)+g_2(x)}2<1$, as required. 
\end{proof}

\subsection{A definition of contractibility via finite spaces} 

\begin{enonce}{Theorem}\label{coro:contrCW} The following are equivalent for a separable completely
	metrisable space $X$ 
\bee
\item[(1)] $X\lra\{\bullet\}\,\in\,\Bigl\{ \MtoLLtoO \Bigr\}^{lr}$
\item[(2)] $X$ is a retract of $\bbbR^\omega$
\eee
Furthermore, if $X$ is also  
	an absolute neighbourhood retract
	for separable metric spaces,\footnote{In more detail this is described in a mathoverflow answer \url{https://mathoverflow.net/questions/436932/when-is-a-contractible-space-a-retract-of-the-hilbert-cube-or-bbb-r-omega} by Tyrone Cutler.} 
the above is also equivalent to each of the following:
\bee
\item[(3)] $X$ is contractible
\item[(4)] $X$ is homotopically trivial
\item[(5)] $X$ is an absolute retract for separable metric spaces
\item[(6)] $X$ is an absolute extensor for separable metric spaces
\eee
\end{enonce}
\begin{proof} 
	(3)$\Leftrightarrow$(4)$\Leftrightarrow$(5): This is the content of \cite[Theorem 4.2.20]{vMill}.\footnote{Note that \cite{vMill} uses the word ``space'' 
	to mean a separable metric space.}
	
	(5)$\Leftrightarrow$(6): \cite[Theorem $1.5.2$]{vMill}. 

	(2)$\implies$(1): ${}^{r}$-classes are closed under retracts and
	products, and by Corollary~\ref{cor:Thm2} $\bbbR\lra\{\bullet\}\,\in\,
	\Bigl\{ \MtoLLtoO \Bigr\}^{lr}$.

	(1)$\implies$(2): any separable completely metrisable space embeds into $\bbbR^\omega$ as a closed subset,\footnote{ 
	For a countable dense subset $(a_i)_{i\in\omega}$ of $X$, 
	a metric on $X$ defines an embedding $X\lra \bbbR^\omega, x\longmapsto \Bigl( y\mapsto \dist(x,a_i) \Bigr)$
	of $X$ as a 
	topological subspace.
	The metric being complete implies this is a closed embedding.}
	hence by Lemma~\ref{lem:rttMtoLLtoO}(4) $X\xra{\{\MtoLLtoO\}^{l}} \bbbR^\omega $,
	hence 
	$X\xra{\{\MtoLLtoO\}^{l}}  \bbbR^\omega \,\rtt\, X\xra{\{\MtoLLtoO\}^{lr}}  \{\bullet\}$, i.e.~$X$ is a retract of $\bbbR^\omega$.

	(5)$\implies$(2): any separable completely metrisable space embeds into separable metric space $\bbbR^\omega$ as a closed subset, 
	hence by the definition of an absolute retract $X$ is a retract of $\bbbR^\omega$. 

        (2)$\implies$(5): separable metric space $\bbbR^\omega$ is an absolute retract, and absolute retracts are closed under retracts.
%
%
%
%
%
%
\end{proof}

	\section{Compactness}\label{sec_Kompakt} 
We shall now use our concise notation for maps of finite topological spaces in \S\ref{not:mintsGE}  
to rewrite the 
Smirnov-Vulikh-Taimanov theorem \cite[$3.2.1$,p.$136$]{Engelking} giving sufficient conditions 
to extend a map to a compact Hausdorff space. We then observe in our reformulation 
all maps of finite spaces involved are closed and thereby proper, and that being proper is a lifting property
by \cite[I\S$10.2$,Theorem 1d)]{Bourbaki}.

This allows to reformulate the notions of compactness (for Hausdorff spaces) and of being proper (for maps of completely normal spaces). 

The focus of the exposition here is to convince the reader that there is finite combinatorics 
implicit in the definition of compactness.
\subsection{Reformulating the Smirnov-Vulikh-Taimanov theorem}\label{subsection_Taimanov}
We start by citing the theorem as it is presented in 
 \cite[$3.2.1$,p.$136$]{Engelking}. This is a topological counterpart of a standard fact in analysis that
 a uniformly continuous function to a complete metric space can always be extended from a dense subset
 to the whole space. 
\newline\noindent \includegraphics[width=\linewidth]{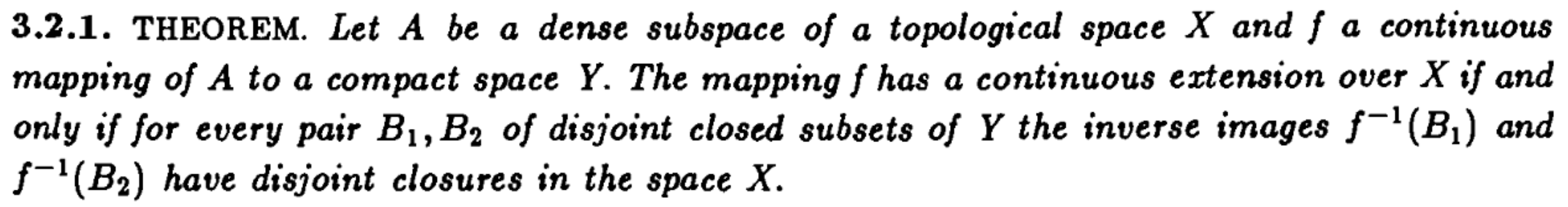}

\subsubsection{Reformulating assumption {\sf $A$ is a dense subspace of $X$}.}\label{sec:densesubset}
This means that the inclusion map $i:A\lra X$ is injective, the topology on $A$ is induced from $X$, and the image of $i:A\lra X$ is dense.
These conditions are represented by the following lifting diagrams involving simple counterexamples to these notions:
\begin{equation}
\begin{gathered}\label{def:Engel:inj}
	i:A\lra X  \text{ is injective iff } 
\vcenter{\hbox{\xymatrix
{  
	A \ar[d]|i\ar[r] & \{\star_a\leftrightarrow \star_b\} \ar[d]
	\\ X \ar[r] \ar@{-->}[ur] & \{\star_{a=b} \}  
	 } }} 
\end{gathered}\end{equation}
	\begin{equation}\begin{gathered}\label{def:Engel:ind}
		\text{ topology on }A\text{ is induced from }X\text{ iff }
		\vcenter{\hbox{\xymatrix
{  
		A \ar[d]|i\ar[r] & \left\{{\ensuremath{\bullet_o}}\rotatebox{-12}{\ensuremath{\to}}  \raisebox{-2pt}{\ensuremath{\bblacksquare_c}}\right\}\  \ar[d]
		\\ X \ar[r] \ar@{-->}[ur] & \{\bullet_{o=c}\}} }} 
\end{gathered}\end{equation}
	\begin{equation}\begin{gathered}\label{def:Engel:dense}
		\text{ the image of }i:A\lra X\text{ is dense }\text{ iff }
		\vcenter{\hbox{\xymatrix
		{       A \ar[d]|i \ar[r] & \{\bblacksquare_c\} \ar[d] \\
		X \ar@{-->}[ur]\ar[r] & \left\{{\ensuremath{\bullet_o}}\rotatebox{-12}{\ensuremath{\to}}  \raisebox{-2pt}{\ensuremath{\bblacksquare_c}}\right\}\  
		}}} 
\end{gathered}\end{equation}

\subsubsection{The conditions on $f:A\lra Y$.}\label{sec:EngelLambda}  To reformulate this we use the characteristic function $\chi_{B_1,B_2}:Y\lra
	\{\bblacksquare_a\raisebox{6pt}{}\rotatebox{12}{\ensuremath{\leftarrow}}\raisebox{2pt}{\ensuremath{\bullet_u}\raisebox{4pt}
	                {}\rotatebox{-13}{\ensuremath{\rightarrow}}} \bblacksquare_b\}$
of the two disjoint closed subsets $B_1,B_2\subset Y$:
	\begin{equation}
\begin{gathered}\label{def:Engel:BB}
	{
		\text{ for every pair }B_1,B_2\text{ of disjoint closed subsets of }Y\atop { \text{\normalsize the inverse images }\normalsize f^{-1}(B_1)\text{ 
	and }f^{-1}(B_2)\atop \text{\normalsize have disjoint closures in the space }X }
	}
	\text{ iff }
	\vcenter{\hbox{\xymatrix
{  
	A \ar[d]|i \ar[r]|f & Y \ar[d] \ar[r]  & \{\bblacksquare_a^{B_1}\raisebox{6pt}{}\rotatebox{12}{\ensuremath{\leftarrow}}\raisebox{2pt}{\ensuremath{\bullet_u}\raisebox{4pt}
			{}\rotatebox{-13}{\ensuremath{\rightarrow}}} \bblacksquare_b^{B_2}\} \ar[dl] 
	\\ X \ar[r] \ar@{-->}[urr] &   \{\bullet_{a=u=b}\}} }} 
\end{gathered}\end{equation}
		\begin{equation}
\begin{gathered}\label{def:Engel:fext}
	{
		\text{ The mapping }f\text{ has a continuous extension over }X
	}
	\text{ iff }
	\vcenter{\hbox{\xymatrix
{  
	A \ar[d]|i \ar[r]|f & Y \ar[d]  &
	\\ X \ar[r] \ar@{-->}[ur] &   \{\bullet_{a=u=b}\}} }} 
\end{gathered}\end{equation}

In particular, 
	\begin{equation}
\begin{gathered}\label{def:Engel:BBall}
	{
		\text{ every pair }B_1,B_2\text{ of disjoint closed subsets of }A
		\atop\text{ 
	have disjoint closures in the space }X 
	}
	\text{ iff }
	\vcenter{\hbox{\xymatrix
{  
	A \ar[d]|i  \ar[r]  & \{\bblacksquare_a^{B_1}\raisebox{6pt}{}\rotatebox{12}{\ensuremath{\leftarrow}}\raisebox{2pt}{\ensuremath{\bullet_u}\raisebox{4pt}
			{}\rotatebox{-13}{\ensuremath{\rightarrow}}} \bblacksquare_b^{B_2}\} \ar[d] 
	\\ X \ar[r] \ar@{-->}[ur] &   \{\bullet_{a=u=b}\}} }} 
\end{gathered}\end{equation}

We find it helpful to list the maps above in various notations
\footnote
{Our notation represents finite topological space as preorders or finite categories with each diagram commuting, \
			and is hopefully self-explanatory; see \S\ref{subsec:notationFP}  for details. 
	In short, 
	an arrow $o\to c$ indicates that $c\in \operatorname{cl}\,o$, and each point goes to ``itself''; 
	the list in $\{..\}$ after the arrow indicates new relations/morphisms added, thus in $\{o\to c\}\lra \{o=c\}$ 
	the equality indicates that the two points are glued together or that we added an identity morphism between $o$ and $c$.
	The notation in the 3rd line informal (red indicates new/added elements), and in the 4th line reminds of a computer syntax.}
as a table:
$$\hskip-42pt\begin{array}{ccccc}
	        \{a\raisebox{6pt}{}\rotatebox{12}{\ensuremath{\leftarrow}}\raisebox{2pt}{\ensuremath{u}\raisebox{4pt}
		        {}\rotatebox{-13}{\ensuremath{\rightarrow}}} b\}\lra\{a=u=b\}
			         &
					 \{a\llrra b\}\lra \{a=b\} &  \left\{{\ensuremath{o}}\rotatebox{-12}{\ensuremath{\to}}  \raisebox{-2pt}{c}\right\}\lra \{o=c\} &
					 \{c\}\lra \left\{{\ensuremath{o}}\rotatebox{-12}{\ensuremath{\to}}  \raisebox{-2pt}{c}\right\} &
					 \\ \text{(disjoint closures)}        & \text{(injective)}            & \text{(pullback topology)}  & \text{(dense image)}&
					 \\ \hskip3pt    \{a\raisebox{6pt}{\color{red}\bf\,=}\!\!\!\!\!\!\rotatebox{12}{\ensuremath{\leftarrow}}\raisebox{2pt}{\ensuremath{u}\raisebox{4pt}{\color{red}\bf\,=}\!\!\!\!\!\!\rotatebox{-13}{\ensuremath{\rightarrow}}} b\} &
					 \{a \leftrightarrow \!\!\!\!\!\!\raisebox{6pt}{{\color{red}\bf\!\!=\,\,}} b\}&
					 \{\raisebox{0pt}{\ensuremath{o}}\raisebox{6pt}{\color{red}\bf\,=}\!\!\!\!\!\!\rotatebox{-12}{\ensuremath{\to}}  \raisebox{-2pt}{c}\}&
						                 \{\raisebox{0pt}{\ensuremath{o}}{\color{red}\bf\rotatebox{-12}{\ensuremath{\to}}  \raisebox{-2pt}{c} }\}

								 \\      \verb|{a<-u->b}-->{a=u=v}|& \verb|{a<->b}-->{a=b}|& \verb|{o->c}-->{o=c}|& \verb|{c}-->{o->c}|&
\end{array}$$

\subsubsection{Reformulating the Smirnov-Vulikh-Taimanov theorem.}
The diagrams above lead us to the following reformulation.
\begin{prop}[($3.2.1$ Theorem)]\label{thm:Engel} Let $Y$ be Hausdorff compact. 
				Let	$i:A\lra X$ satisfy (\ref{def:Engel:inj}), (\ref{def:Engel:ind}), (\ref{def:Engel:dense}), i.e. 
	$$ {A\atop{{\downarrow}\atop X}} \,\,\in\,\, \left(
	{ {\{\star_a\leftrightarrow \star_b\}}\atop{\downarrow
	\atop{ \{\star_{a=b} \} } }}   \,\,\,
	{\left\{{\ensuremath{\bullet_o}}\rotatebox{-12}{\ensuremath{\to}}  \raisebox{-2pt}{\ensuremath{\bblacksquare_c}}\right\}\atop{\downarrow\atop{\{\bullet_{o=c}\}}}}
	\,\,\,{{{\{\bblacksquare_c\}}\atop{\downarrow
		\atop{\left\{{\ensuremath{\bullet_o}}\rotatebox{-12}{\ensuremath{\to}}  \raisebox{-2pt}{\ensuremath{\bblacksquare_c}}\right\}\  
		}}}}\right)^{l}
	$$ 
				Then for each $f:A\lra Y$ (\ref{def:Engel:BB}) and (\ref{def:Engel:fext}) are equivalent. 
In particular, 
	$$ {A\atop{{\downarrow}\atop X}} \,\,\in\,\, \left(
	{	\{\bblacksquare_a\raisebox{6pt}{}\rotatebox{12}{\ensuremath{\leftarrow}}\raisebox{2pt}{\ensuremath{\bullet_u}\raisebox{4pt}
				{}\rotatebox{-13}{\ensuremath{\rightarrow}}} \bblacksquare_b\} }\atop{\downarrow
	\atop{ \{\star_{a=u=b} \} } }   \,\,\,
		\right)^{l}
	\text{ implies } 
	 {A\atop{{\downarrow}\atop X}} \,\,\in\,\, \left(
	 Y\atop{\downarrow \atop \ptt} 
		\right)^{l}
	$$

\end{prop}
\begin{proof} The meaning/translation of the diagrams in the natural language in \S\ref{sec:densesubset}-\S\ref{sec:EngelLambda}  gives precisely the statement of 
\cite[$3.2.1$,p.$136$]{Engelking} cited above. \end{proof}

\subsubsection{Compactness in terms of the Quillen negation (orthogonals)}
Collecting the orthogonals together gives one implication of the following Corollary. 

\begin{coro}\label{coro:KtoOlr} A Hausdorff space $Y$ is  compact iff 
	\begin{equation}\label{KtoOlr} 
		{Y\atop{{\downarrow}\atop \ptt}} \,\,\in\,\, \left(
	{ {\{\star_a\leftrightarrow \star_b\}}\atop{\downarrow
	\atop{ \{\star_{a=b} \} } }}   \,\,\,
	{\left\{{\ensuremath{\bullet_o}}\rotatebox{-12}{\ensuremath{\to}}  \raisebox{-2pt}{\ensuremath{\bblacksquare_c}}\right\}\atop{\downarrow\atop{\{\bullet_{o=c}\}}}}
	\,\,\,
	{{{\{\bblacksquare_c\}}\atop{\downarrow
		\atop{\left\{{\ensuremath{\bullet_o}}\rotatebox{-12}{\ensuremath{\to}}  \raisebox{-2pt}{\ensuremath{\bblacksquare_c}}\right\}\  
	 }}}}%
	{{	\{\bblacksquare_a\raisebox{6pt}{}\rotatebox{12}{\ensuremath{\leftarrow}}\raisebox{2pt}{\ensuremath{\bullet_u}\raisebox{4pt}
				{}\rotatebox{-13}{\ensuremath{\rightarrow}}} \bblacksquare_b\} }\atop{\downarrow
	\atop{ \{\star_{a=u=b} \} } }}   \,\,\,
	\right)^{lr} \end{equation}
or, equivalently, iff 
	\begin{equation}\label{KtoOlrSimplified} 
		{K\atop{{\downarrow}\atop \ptt}} \,\,\in\,\, \left(
		{	\{\bblacksquare_a\raisebox{6pt}{}\rotatebox{12}{\ensuremath{\leftarrow}}\raisebox{2pt}{\ensuremath{\bullet_u}\raisebox{4pt}
				{}\rotatebox{-13}{\ensuremath{\rightarrow}}} \bblacksquare_b \leftrightarrow \filledstar_x \} }\atop{\Downarrow
	\atop{ 
		\left\{{\ensuremath{\bullet_o}}\rotatebox{-12}{\ensuremath{\to}}  \raisebox{-2pt}{\ensuremath{\bblacksquare_{a=u=b=x}}}\right\}
		} }
		\right)^{lr}
	\end{equation} 
\end{coro} 
\begin{proof} The combinatorial argument showing equivalence of the two reformulations is similar to that of Proposition~\ref{prop9}.
	We leave it to the reader to verify that
	each of the maps in the first equation is a retract of the map in the second, which in turn is a composition of pullbacks of the maps in the first equation. Now let us prove the main statement.

$\implies$: This follows formally from Proposition \ref{thm:Engel}.
Indeed, we need to show that $A \lra X \,\,\rtt\,\, Y \lra \ptt$ whenever 
	\begin{equation}\label{AtoXl}  
	{A\atop{{i\downarrow\,}\atop X}} \,\,\in\,\, \left(
	{ {\{\star_a\leftrightarrow \star_b\}}\atop{\downarrow
	\atop{ \{\star_{a=b} \} } }}   \,\,\,
	{\left\{{\ensuremath{\bullet_o}}\rotatebox{-12}{\ensuremath{\to}}  \raisebox{-2pt}{\ensuremath{\bblacksquare_c}}\right\}\atop{\downarrow\atop{\{\bullet_{o=c}\}}}}
	\,\,\,
	{{{\{\bblacksquare_c\}}\atop{\downarrow
		\atop{\left\{{\ensuremath{\bullet_o}}\rotatebox{-12}{\ensuremath{\to}}  \raisebox{-2pt}{\ensuremath{\bblacksquare_c}}\right\}\  
	 }}}}%
	{{	\{\bblacksquare_a\raisebox{6pt}{}\rotatebox{12}{\ensuremath{\leftarrow}}\raisebox{2pt}{\ensuremath{\bullet_u}\raisebox{4pt}
				{}\rotatebox{-13}{\ensuremath{\rightarrow}}} \bblacksquare_b\} }\atop{\downarrow
	\atop{ \{\star_{a=u=b} \} } }}   \,\,\,
		\right)^{l}
	\end{equation}
This means precisely that $i:A\lra X$ satisfies (\ref{def:Engel:inj}),
	(\ref{def:Engel:ind}), (\ref{def:Engel:dense}), and
	(\ref{def:Engel:BBall}).  Now,  (\ref{def:Engel:BBall}) implies that
	each $f:A\lra Y$ satisfies \eqref{def:Engel:BB} and therefore also
	(\ref{def:Engel:fext}), as required.

$\Longleftarrow$:
Following \cite[I\S$6.5$,Example,p.$62$]{Bourbaki} let us define {\em the
	topological space $X \sqcup_{\mathfrak F} \{\infty\}$ associated with a filter $\mathfrak F$} on the set of
	points of a topological space $X$. 
	Its set of points $X \sqcup  \{\infty\}$  is obtained by adjoining to $X$ a new point $\infty\not\in X$.
A subset of $X \sqcup \{\infty\}$  is defined to be {\em open} iff it is either an open subset of $X$ or of form $U\cup\{\infty\}$ 
				where $U\in \mathfrak F$ and $U$ is an open subset of $X$;
	equivalently, a subset of  $X \sqcup \{\infty\}$  is defined to be {\em closed} iff it is either of form $Z\cup\{\infty\}$ where $Z$ is a closed subset of $X$,
	or is a closed subset $Z$ of $X$ such that $X\setminus Z \,\in\,\mathfrak F$. 
%

	By \cite[I\S$7.1$, Def. 1]{Bourbaki} a filter $\mathfrak F$ converges to a point $x\in X$ iff 
	$U_x\in \mathfrak F$ for each neighbourhood $U_x \ni x$. The latter is equivalent to the continuity of the map 
$\id_x: X\sqcup_{\mathfrak F }\{\infty\}\lra X$ defined by $\id_x(\infty):=x$, and $\id_x(x):=x$ for all $x\in X$.
	This leads to the following reformulation of  \cite[I\S$7.1$, Def. 1]{Bourbaki}. A filter $\mathfrak F$ on $X$ is convergent iff 
	$$\xymatrix@C=1.39cm{X \ar[r]|{\id} \ar[d]   & X  \ar[d] \\  X\sqcup_{\mathfrak F}\{\infty\}\ar@{-->}[ru]  \ar[r] &\ptt }$$

	By \cite[I\S$9.1$, Def. 1(C$'$)]{Bourbaki} a topological space $X$ is quasi-compact iff each ultrafilter $\mathfrak U$ on $X$ is convergent.
	A verification shows that $X$ is a dense subset of $X \sqcup_{\mathfrak U} \{\infty\}$, and, moreover, 
		every pair $B_1,B_2$  of disjoint closed subsets of $X$ 
	has disjoint closures in $X$, i.e. 
	$$ {X\atop{{\downarrow}\atop X\sqcup_{\mathfrak U} \{\infty\}}} \,\,\in\,\, \left(
	{ {\{\star_a\leftrightarrow \star_b\}}\atop{\downarrow
	\atop{ \{\star_{a=b} \} } }}   \,\,\,
	{\left\{{\ensuremath{\bullet_o}}\rotatebox{-12}{\ensuremath{\to}}  \raisebox{-2pt}{\ensuremath{\bblacksquare_c}}\right\}\atop{\downarrow\atop{\{\bullet_{o=c}\}}}}
	\,\,\,
	{{{\{\bblacksquare_c\}}\atop{\downarrow
		\atop{\left\{{\ensuremath{\bullet_o}}\rotatebox{-12}{\ensuremath{\to}}  \raisebox{-2pt}{\ensuremath{\bblacksquare_c}}\right\}\  
	 }}}}%
	{{	\{\bblacksquare_a\raisebox{6pt}{}\rotatebox{12}{\ensuremath{\leftarrow}}\raisebox{2pt}{\ensuremath{\bullet_u}\raisebox{4pt}
				{}\rotatebox{-13}{\ensuremath{\rightarrow}}} \bblacksquare_b\} }\atop{\downarrow
	\atop{ \{\star_{a=u=b} \} } }}   \,\,\,
		\right)^{l}
	$$ 
	Hence, 
	$$ {K\atop{{\downarrow}\atop \ptt}} \,\,\in\,\, \left(
	{ {\{\star_a\leftrightarrow \star_b\}}\atop{\downarrow
	\atop{ \{\star_{a=b} \} } }}   \,\,\,
	{\left\{{\ensuremath{\bullet_o}}\rotatebox{-12}{\ensuremath{\to}}  \raisebox{-2pt}{\ensuremath{\bblacksquare_c}}\right\}\atop{\downarrow\atop{\{\bullet_{o=c}\}}}}
	\,\,\,
	{{{\{\bblacksquare_c\}}\atop{\downarrow
		\atop{\left\{{\ensuremath{\bullet_o}}\rotatebox{-12}{\ensuremath{\to}}  \raisebox{-2pt}{\ensuremath{\bblacksquare_c}}\right\}\  
	 }}}}%
	{{	\{\bblacksquare_a\raisebox{6pt}{}\rotatebox{12}{\ensuremath{\leftarrow}}\raisebox{2pt}{\ensuremath{\bullet_u}\raisebox{4pt}
				{}\rotatebox{-13}{\ensuremath{\rightarrow}}} \bblacksquare_b\} }\atop{\downarrow
	\atop{ \{\star_{a=u=b} \} } }}   \,\,\,
		\right)^{lr}
	$$ 
implies that for each ultrafilter $\mathfrak U$ on $K$
	\begin{equation}\label{DefCompUltras} 
	\xymatrix@C=1.39cm{K \ar[r]|{\id} \ar[d]   & K  \ar[d] \\  K\sqcup_{\mathfrak U}\{\infty\}\ar@{-->}[ru]  \ar[r] &\ptt }\end{equation} 
That is, each ultrafilter $\mathfrak U$ on $K$ converges, and  
	by \cite[I\S$9.1$, Def. 1(C$'$)]{Bourbaki} $K$ is quasi-compact.
\end{proof}

\subsection{The definition of compactness and properness in terms of ultrafilters as a lifting property}
	Bourbaki state almost explicitly the lifting property defining compactness in  \cite[I\S$9.1$, Def. 1($\mathrm C'\implies \mathrm C''$)]{Bourbaki} 
	while proving equivalence of various definitions of compactness:
\newline\includegraphicss[0.99]{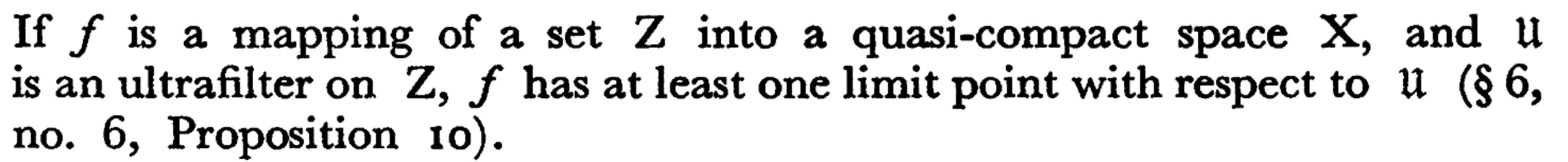}
	
	By  \cite[I\S$10.2$,Lemma $1$, also Theorem 1, Corollary 1]{Bourbaki} 
	a space $X$ is quasi-compact iff the map $X\lra\ptt$ is proper,
	and
	the same lifting property appears in the characterisation of properness \cite[I\S$10.2$, Theorem 1d)]{Bourbaki}. 
\newline\noindent\includegraphics[test]{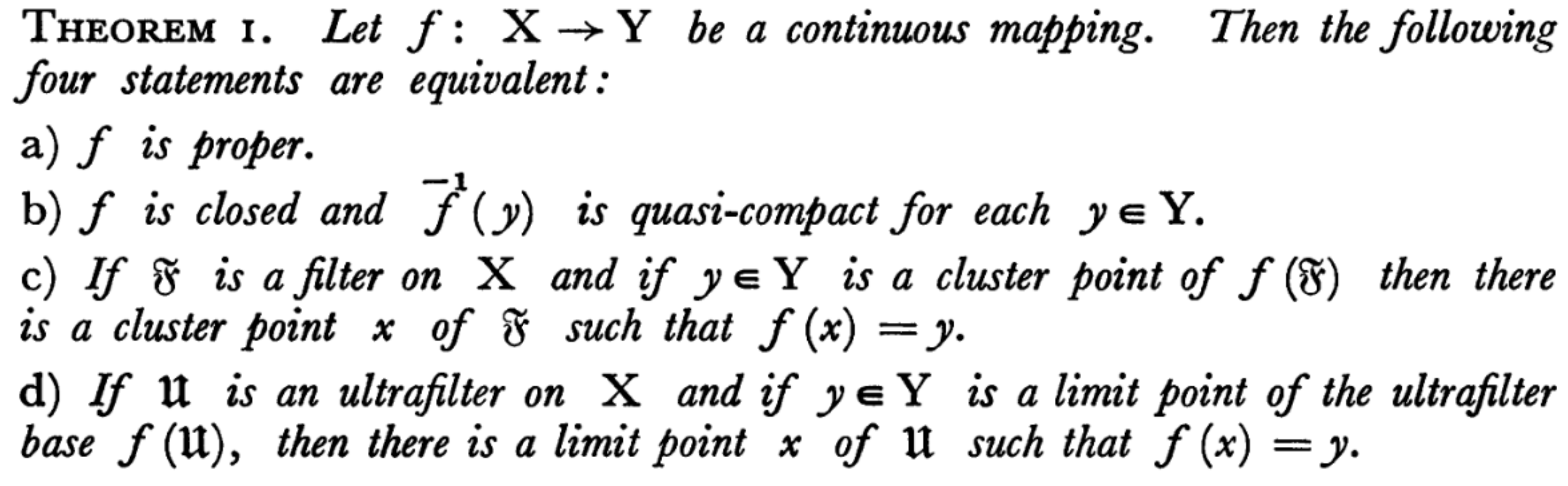}

The next proposition reformulates in our notation this remark and 
	\cite[I\S$10.2$, Theorem 1d)]{Bourbaki}, as lifting properties using the topological space associated with an ultrafilter. 
\begin{prop}[{\cite[I\S$10.2$, Theorem 1d)]{Bourbaki}}]\label{propo_Theorem1d}  
A topological space $X$ is quasi-compact iff 
	\begin{equation} {X\atop{{\downarrow}\atop \ptt}} \,\,\in\,\, \bigl\{\, Z \lra Z\sqcup_{\mathfrak U} \{\infty\} \,\,:\,\, \mathfrak U\text{ is an ultrafilter on a set }Z \,\bigr\}^r \label{ultrafiltersL}\end{equation} 
In fact,  
	\begin{equation} 
	\bigl\{\, Z \lra Z\sqcup_{\mathfrak U} \{\infty\} \,\,:\,\, \mathfrak U\text{ is an ultrafilter on a set }Z \,\bigr\}^r \label{ultrafiltersL}\end{equation} 
is the class of proper maps. 
\end{prop} 
\begin{proof} The proofs of both claims amounts to decyphering the notation. We spell out only the first proof.  
If $f$ is a mapping of a set $Z$ into a quasi-compact space $X$, and $\mathfrak U$
	is an ultrafilter on $Z$, $f$ has at least one limit point with respect to $\mathfrak U$ 
	\cite[I\S$9.1$, Def. 1($(C)\implies (C')$)]{Bourbaki}.  
The following diagram  expresses this property using
	the topological space 
	associated with a filter $\mathcal U$ 
	defined in \cite[I\S$6.5$,Example,p.$62$]{Bourbaki}. 
\begin{equation}\label{DefCompUltrasS} \xymatrix@C=1.39cm{Z \ar[r]|{f} \ar[d]   & X  \ar[d] \\  Z\sqcup_{\mathfrak U}\{\infty\}\ar@{-->}[ru]  \ar[r] &\ptt }\end{equation} 
	That is, 
	\begin{equation}  Z \lra Z \sqcup_{\mathfrak U} \{\infty\} \,\,\rtt\,\, X \lra \ptt \label{quasicompactUltrafilter} \end{equation}
		On the other hand, an ultrafilter $\mathfrak U$ on $Z:=X$ converges iff the mapping $f:=\id_X$ has at least one limit point with respect to $\mathfrak U$ by  \cite[I\S$7.3$, Proposition $7$]{Bourbaki},
	and thus \eqref{quasicompactUltrafilter}  implies that $X$ is quasi-compact  by \cite[I\S$9.1$, Def. 1($\mathrm C'$)]{Bourbaki}.
\end{proof}

\subsection{The definition of compactness in terms of the simplest counterexample}
\subsubsection{A class of finite closed maps implicit in the theorem}
			A verification shows that a map $g$ of finite topological spaces is {\em closed} iff 
			$\{\bullet_o\}\lra \left\{{\ensuremath{\bullet_o}}\rotatebox{-12}{\ensuremath{\to}}  \raisebox{-2pt}{\ensuremath{\bblacksquare_c}}\right\}
			\,\rtt\,
			g$, and that this lifting property holds for each of the maps in \eqref{KtoOlr}.  
Let $\Bigl\{ \{\bullet_o\}\lra \left\{{\ensuremath{\bullet_o}}\rotatebox{-12}{\ensuremath{\to}}  \raisebox{-2pt}{\ensuremath{\bblacksquare_c}}\right\} \Bigr\}^r_{<5}$
denote the subclass of $\Bigl\{ \{\bullet_o\}\lra \left\{{\ensuremath{\bullet_o}}\rotatebox{-12}{\ensuremath{\to}}  \raisebox{-2pt}{\ensuremath{\bblacksquare_c}}\right\} \Bigr\}^r$
consisting of maps of spaces of size less than $5$. 
\def\antiProper
	{{{\{\bblacksquare_c\}}\atop{\downarrow
		\atop{\left\{{\ensuremath{\bullet_o}}\rotatebox{-12}{\ensuremath{\to}}  \raisebox{-2pt}{\ensuremath{\bblacksquare_c}}\right\}\  
	 }}}}%
The above gives a cleaner proof of an analogue of Corollary~\ref{coro:KtoOlr}: 
\begin{coro}\label{coro:Engel:2}\label{coro:proper} A Hausdorff space $K$ is compact iff 
\begin{equation}\label{KtoOlr:2}  {K\atop{{\downarrow}\atop \ptt}} \,\,\in\,\, 
	\left( \left(
	{{{\{\bullet_o\}}\atop{\downarrow
		\atop{\left\{{\ensuremath{\bullet_o}}\rotatebox{-12}{\ensuremath{\to}}  \raisebox{-2pt}{\ensuremath{\bblacksquare_c}}\right\}\  
	 }}}}%
	\right)_{<5} 	\right)^{lr}
\end{equation}
Moreover,
		$$\left(  \left(
	{{{\{\bullet_o\}}\atop{\downarrow
		\atop{\left\{{\ensuremath{\bullet_o}}\rotatebox{-12}{\ensuremath{\to}}  \raisebox{-2pt}{\ensuremath{\bblacksquare_c}}\right\}\  
	 }}}}%
		\right)_{<5} \right)^{lr}$$
		is contained in the class of proper maps.
\end{coro}
\begin{proof} The second claim follows formally from the facts that the class of proper maps is defined by a lifting property, and 
	that each map of finite spaces lifting with respect to 
	$\{\bullet_o\}\lra \{ {\ensuremath{\bullet_o}}\rotatebox{-12}{\ensuremath{\to}}  \raisebox{-2pt}{\ensuremath{\bblacksquare_c}}\}$ 
	is closed and therefore proper. 
	The first claim follows from Corollary~\ref{coro:KtoOlr} and the fact that each map in   Equation~\eqref{KtoOlr} is proper.
\end{proof}

\subsubsection{Conclusions: a definition of compactness via the simplest counterexample ?} This leads to the following conjecture
defining compactness in terms of the simplest counterexample. 

\begin{conj}[Compactness via simplest counterexample]\label{conj:Engel}\label{conj:proper}
	$$\left( \left( 
	{{{\{\bullet_o\}\hphantom{\to\bblacksquare}}\atop{\downarrow
		\atop{\displaystyle \left\{{\ensuremath{\bullet_o}}\rotatebox{-12}{\ensuremath{\to}}  \raisebox{-2pt}{\ensuremath{\bblacksquare_c}}\right\}\  
	 }}}}%
		\right)_{<5} \right)^{lr}$$
		is the class of proper maps.
\end{conj}
\begin{proof}[Evidence] The proof of Theorem~\ref{thm:Engel} in \cite[Theorem  $3.2.1$]{Engelking} in fact shows that a proper map of completely normal spaces lies in the class. 
Proposition~\ref{propo_Theorem1d} implies that being proper is defined by a right lifting property, and
thus each map in the left-then-right negation of a class of proper map is necessarily proper.
\end{proof}

\subsubsection{Conclusions: a definition of compactness via a random example ?} 
The following corollary says that almost any example of non-surjective proper map ``complicated enough'' defines compactness for Hausdorff spaces. 
The condition ``complicated enough'' is purely combinatorial and means that certain two maps occur as retractions of $g$.
This makes it easy to give a lower bound on the probability that a random map of small spaces defines compactness. 

\begin{coro}\label{coro:KtoOlrSimplifiedFurther} 
Let $g:X\lra Y$ be a non-surjective closed map of finite spaces ``complicated enough'' in the sense that 
	(i) it has a fibre with two points bounded above but not below, and 
	(ii) its image is not a clopen subset.  

	Then a  Hausdorff space $K$ is  compact iff $K\lra\ptt \,\in \, g^{lr}$. 

\end{coro}
\begin{proof} A map of finite spaces is proper iff it is closed. 
If a map $g:X\lra Y$ is proper, then so in each map $g^{lr}$. 
On the other hand, 
Quillen negations are closed under retracts, 
and thus $g^{lr}$ contains all the retracts of $g$. 
If this includes each map in eq.~\eqref{KtoOlrSimplified} then by 
Corollary~\ref{coro:KtoOlr} for each compact Hausdorff space $K$
it holds $K\lra \ptt \,\in\, g^{lr}$. 

Finally, a verification shows a surjective map in 
 $   {{      \{\bblacksquare_a\raisebox{6pt}{}\rotatebox{12}{\ensuremath{\leftarrow}}\raisebox{2pt}{\ensuremath{\bullet_u}\raisebox{4pt}
	                                                                         {}\rotatebox{-13}{\ensuremath{\rightarrow}}} \bblacksquare_b\} }\atop{\downarrow
										                                                                                  \atop{ \{\star_{a=u=b} \} } }}   $
	is a retract of a map $g$ iff $g$ satisfies (i).
\footnote{The retraction sends each $x\in X$ 
	into the least upper bound (=colimit) 
	in $X'$ of $\{x': x\ra x'\}$. The map preserves the order because $x_1\lra x_2$ implies $\{x': x_1\lra x'\}\supset \{x':x_2\lra x'\}$
	and thus the same holds for the least upper bounds of these sets. We used that in our examples the least upper bound always exist,
	and that for no element $\{x': x\ra x'\}$ is empty.}

A verification shows that the map ${\{\bblacksquare_c\}}\atop{\downarrow
                \atop{\left\{{\ensuremath{\bullet_o}}\rotatebox{-12}{\ensuremath{\to}}  \raisebox{-2pt}{\ensuremath{\bblacksquare_c}}\right\}\
			 }}$ is a retract of any non-surjective closed map of finite spaces satisfying (ii).
The map $        {\left\{{\ensuremath{\bullet_o}}\rotatebox{-12}{\ensuremath{\to}}  \raisebox{-2pt}{\ensuremath{\bblacksquare_c}}\right\}\atop{\downarrow\atop{\{\bullet_{o=c}\}}}} $ is a retract of $   {{      \{\bblacksquare_a\raisebox{6pt}{}\rotatebox{12}{\ensuremath{\leftarrow}}\raisebox{2pt}{\ensuremath{\bullet_u}\raisebox{4pt}
					                                 {}\rotatebox{-13}{\ensuremath{\rightarrow}}} \bblacksquare_b\} }\atop{\downarrow
									         \atop{ \{\star_{a=u=b} \} } }}   $.
It is left to show that we need not consider $       { {\{\star_a\leftrightarrow \star_b\}}\atop{\downarrow
	\atop{ \{\star_{a=b} \} } }}$. The map $ A\lra B \,\rtt\,       {\left\{{\ensuremath{\bullet_o}}\rotatebox{-12}{\ensuremath{\to}}  \raisebox{-2pt}{\ensuremath{\bblacksquare_c}}\right\}\atop{\downarrow\atop{\{\bullet_{o=c}\}}}}$ iff the topology on $A$ is induced from $B$. In particular, each fibre of $A\lra B$ is indiscreet, and thus  
	any map $A\lra K$ to a Hausdorff space $K$ sends each fibre of $A\lra B$ to a single point. Hence, to construct a lifting 
	for the commutative square  for the lifting property $A\lra B \,\rtt\, K\lra\ptt$ it is enough to construct a lifting for $A'\lra B \,\rtt\, K\lra \ptt$
	where $A'$ is the quotient of $A$ making the map $A'\lra B$ is injective, which means precisely that 
	$A'\lra B \,\rtt\, { {\{\star_a\leftrightarrow \star_b\}}\atop{\downarrow \atop{ \{\star_{a=b} \} } }} $. 
	It is only left to verify that if $A\lra B$ has the left lifting property with respect to each map in Eq.~\eqref{KtoOlrSimplified} except possibly the first one,
	then so does $A'\lra B$, and thus $A'\lra B$ lifts with respect to each map in Eq.~\eqref{KtoOlrSimplified}. 
	Therefore $A\lra B \,\rtt\,K\lra \ptt$ as required.
\end{proof}


\section{The Brouwer fixed point theorem}

We reformulate the Brouwer fixed point theorem   
in terms of finite topological spaces
using our reformulations of 
contractibility and compactness.


	\begin{enonce*}{Theorem $5.1$}\label{thm:BrowerFP}
	Let $K$ be a separable metrisable space. Let $K \xra { \left\{\MtoLLtoO\right\}^{lr}} \ptt$ and 
	$K \xra { \left(\left\{ \{o\}\lra \left\{{\ensuremath{o}}\rotatebox{-12}{\ensuremath{\to}}  \raisebox{-2pt}{c}\right\} \right\}^r_{<5}\right)^{lr} } \ptt$.
	Then each endomorphism of $K$ has a fixed point.
\end{enonce*}
\begin{proof} Corollary~\ref{coro:proper} says that 
	$K \xra { \left(\left\{ \{o\}\lra \left\{{\ensuremath{o}}\rotatebox{-12}{\ensuremath{\to}}  \raisebox{-2pt}{c}\right\} \right\}^r_{<5}\right)^{lr} } \ptt$
	means that $K$ is compact. 
	A separable metrisable compact space embeds into $[0,1]^\omega$ as a closed subspace, and let $K\lra [0,1]^\omega$ denote such an embedding.
	By Lemma~\ref{lem:rttMtoL}(3) 
	$K\xra{ \left\{\MtoLLtoO\right\}^{l} }  [0,1]^\omega$, and thereby  
	$K\xra{ \left\{\MtoLLtoO\right\}^{l} }  [0,1]^\omega \,\rtt\, K  \xra { \left\{\MtoLLtoO\right\}^{lr}} \ptt$
	implies that $K$ is a retract of $[0,1]^\omega$. 
	By the Brouwer fixed point theorem each endomorphism of a retract\footnote{To see that this is implied by the Brouwer fixed point theorem for $[0,1]^\omega$, 
	take a fixed point of $[0,1]^\omega\lra K \lra K \lra [0,1]^\omega$.} 
	of $[0,1]^\omega$ has a fixed point,
	and the claim follows.
\end{proof}

	\begin{enonce*}{Conjecture 5.2}\label{conj:BrowerFP}
The theorem above holds for an arbitrary space $K$. \end{enonce*}

\def\MtoLLtoOL{\left\{\MtoLLtoO\right\}^l}

\begin{enonce*}{Remark 5.3}
	The following diagrams represent the statement of the conjecture above,
	We wonder if they are appropriate for any system of formalised mathematics. 
\begin{equation}\begin{gathered}\label{def:Brouwer:fp}
	\vcenter{\hbox{\xymatrix
	{       & \{\bullet\} \ar@{-->}[ld]|x \ar@{-->}[rd]|x  & \\ 
	K \ar[rr]  \ar[dr]_(0.42){ 
	\left(\left\{ \{o\} \longrightarrow \left\{ o
	                           \rotatebox{-12}{\ensuremath{\to}}\raisebox{-2pt}{\ensuremath{c}}
											 \right\}\right\}^r_{<5}\right)^{\rtt lr} 
	}  
	& & K \ar[dl]^(0.42){\{\MtoLLtoOL\}^{\rtt lr}} 
	\\  & \{\bullet\} &  } }}
	\ \ \ \ 
	\vcenter{\hbox{\xymatrix
	{       \{\bullet\} \ar@{-->}[d]   \\ 
	K \ar@(ur,dr)
	\ar[d]_(0.42){ 
	\left\{ \{o\} \longrightarrow \left\{ o
				   \rotatebox{-12}{\ensuremath{\to}}\raisebox{-2pt}{\ensuremath{$c$}  } 	\right\} 
				   \right\}^{\rtt lr} 
	}^  
	{\{\MtoLLtoOL\}^{\rtt lr}} 
	\\   \{\bullet\}   } }}
\end{gathered}\end{equation}
\end{enonce*}
%
%

\begin{flushright}
 \tiny 
 \begin{minipage}[t]{0.7\textwidth}
 \tiny {{\tiny Die Mathematiker sind eine Art Franzosen: Redet man zu ihnen, so 
 \"ubersetzen sie es in ihre Sprache, und dann ist es alsobald ganz etwas 
 anderes.\ } \\
 \tiny\ \ \,\,\,\,\, ------ Johann Wolfgang von Goethe. Maximen und Reflexionen. 
Aphorismen und Aufzeichnungen. Nach den Handschriften des Goethe- und Schiller-Archivs hg. von Max Hecker, Verlag der Goethe-Gesellschaft, 
Weimar 1907, 
Aus dem Nachlass, Nr.~1005, Uber Natur und Naturwissenschaft.
}\end{minipage}
\end{flushright}

\section{Notation for finite topological spaces\label{not:mintsGE}}
\label{subsec:notationFP} 
The exposition here is a slight modification of  \cite[\S$5.3$.1]{mintsGE}.
This notation lies at heart of the paper, and is perhaps the main contribution. 

\subsection{\label{app:top-notation} Finite topological spaces as preorders and as categories} 
A {\em topological space} comes with a {\em specialisation preorder} on its points: for
points $o,c \in X$,  $c \leqslant o$ iff $c \in \operatorname{cl} o$ ($c$ is in the {\em topological closure} of $o$).

The resulting {\em preordered set} may be regarded as a {\em category} whose
{\em objects} are the points of ${X}$ and where 
there is a unique {\em morphism} $o\lra c$ 
iff $c \in \operatorname{cl} o$.

For a {\em finite topological space} $X$, the specialisation preorder or
equivalently the corresponding category uniquely determines the space: a {\em
subset} of ${X}$ is {\em closed} iff it is
{\em downward closed}, or equivalently,
is a full subcategory such that there are no morphisms going outside the subcategory.

The monotone maps (i.e. {\em functors}) are the {\em continuous maps} for this topology.

\subsection{Notation for finite topological spaces and their maps} 
We denote a finite topological space by a list of the arrows (morphisms) in
the corresponding category; '$\leftrightarrow $' denotes an {\em isomorphism} and '$=$' denotes the {\em identity morphism}.  An arrow between two such lists
denotes a {\em continuous map} (a functor) which sends each point to the correspondingly labelled point, but possibly turning some morphisms into identity morphisms, thus gluing some points.

With this notation, we may display continuous functions for instance between the {\em discrete space} on two points, the {\em Sierpinski space}, the {\em antidiscrete space} and the {\em point space} as follows (where each point is understood to be mapped to the point of the same name in the next space):
$$
  \begin{array}{ccccccc}
	    \{a,b\}
	         &\longrightarrow&
		   \{a{\searrow}b\}
		        &\longrightarrow&
			  \{a\leftrightarrow b\}
			      &\longrightarrow&
			        \{a=b\}
				  \\
				    \text{(discrete space)}
				         &\longrightarrow&
					   \text{(Sierpinski space)}
					       &\longrightarrow&
					         \text{(antidiscrete space)}
						     &\longrightarrow&
						       \text{(single point)}
						         \end{array}
							 $$

   Each continuous map $A\lra B$ between finite spaces may be represented in this way; in the first list
   list relations between elements of $A$, and in the second list put relations between their images.
   However, note that this notation does not allow to represent {\em endomorphisms $A\lra A$}.
   We think of this limitation
   as a feature and not a bug: in a diagram chasing computation,
   endomorphisms under transitive closure lead to infinite cycles,
   and thus our notation has better chance to define a computable fragment of topology.

\subsection{Various conventions on naming points and depicting arrows}
While efficient, this notation is unconventional and requires some getting used to. For this reason, sometimes we employ 
more graphic notation where our notation is moved to subscripts, so to say:
points or objects are denoted by bullets $\bullet,\bblacksquare, \star,...$ with subscripts,
and the reader may think that the subscripts indicate where a point maps to. 
The shape of the bullet indicates whether the point is open, closed, or neither:
$\bullet$ stands for open points (which might also be closed),
$\bblacksquare$ stands for closed points, and $\star$ stands for points which are neither open or closed.
Importantly?, this notation makes visually apparent the shape of the preorder denoted.

Thus in this graphic notation we would write
$$
  \begin{array}{ccccccc}
	    \{\bullet_a,\bullet_b\}
	         &\longrightarrow&
		    \left\{{\ensuremath{\bullet_a}}\rotatebox{-12}{\ensuremath{\to}}  \raisebox{-2pt}{\ensuremath{\bblacksquare_b}}\right\}
		        &\longrightarrow&
			  \{\star_a\leftrightarrow \star_b\}
			      &\longrightarrow&
				\{\bullet_{a=b}\}
				 \\
				    \text{(discrete space)}
				         &\longrightarrow&
					   \text{(Sierpinski space)}
					       &\longrightarrow&
					         \text{(antidiscrete space)}
						     &\longrightarrow&
						       \text{(single point)}
						         \end{array}
							 $$
In $A \longrightarrow  B$, each object and each morphism in $A$ necessarily appears in ${B}$ as well; sometimes we avoid listing
the same object or morphism twice. Thus
both
$$
\{a\} \longrightarrow  \{a,b\}
  \phantom{AAA} \text{ and } \phantom{AAA}
  \{a\} \longrightarrow  \{b\}
  $$
  denote the same map from a single point to the discrete space with two points.

   In $\{a{\searrow}b\}$, the point $a$ is open and point ${b}$ is closed. We
   denote points by $a,b,c,..,U,V,...,0,1..$ to make notation reflect the intended meaning,
   e.g.~$X\lra \{U\searrow U'\}$ reminds us that the preimage of $U$ determines an open subset of $X$,
   $\{x,y\}\lra X$ reminds us that the map determines points $x,y\in X$, and $\{o\searrow c\}$ reminds that $o$ is open and $c$ is closed.

   Each continuous map $A\lra B$ between finite spaces may be represented in this way; in the first list
   list relations between elements of $A$, and in the second list put relations between their images.
   However, note that this notation does not allow to represent {\em endomorphisms $A\lra A$}.
   We think of this limitation
   as a feature and not a bug: in a diagram chasing computation,
   endomorphisms under transitive closure lead to infinite cycles,
   and thus our notation has better chance to define a computable fragment of topology.

\subsection{An example of our notation conventions} 
Both
   $$\{a{\swarrow}U{\searrow}x{\swarrow}V{\searrow}b\}\longrightarrow \{a{\swarrow}U=x=V{\searrow}b\}\text{ and }\{a{\swarrow}U{\searrow}x{\swarrow}V{\searrow}b\}\longrightarrow \{U=x=V\}$$
   denote the morphism gluing points $U,x,V$. More graphically the same morphism can be denoted both as: 
\begin{equation}\label{MtoL}\begin{gathered}
\xymatrix
{  
{
\left\{\underset{\bblacksquare_A}{}{\swarrow} \overset{\raisebox{1pt}{$ \bullet_{U}$}}{}{\searrow} 
\underset{\bblacksquare_{X}}{}
{\swarrow}\overset{\raisebox{1pt}{$\bullet_{V}$}}{}{\searrow}\underset{\bblacksquare_B}{}\right\}
}
 \ar[d]^{
}
         \\
{	 
\left\{\underset{\bblacksquare_A}{}{\swarrow} \overset{\raisebox{1pt}{$ \bullet_{U=X=V}$}}{}{\searrow} 
\underset{\bblacksquare_B}{}\right\}
}
}
	{\phantom{AAAAAAAAAAAAAA}}
\xymatrix
{  
{
\left\{\underset{\bblacksquare_a}{}{\swarrow} \overset{\raisebox{1pt}{$ \bullet_{u}$}}{}{\searrow} 
\underset{\bblacksquare_{x}}{}
{\swarrow}\overset{\raisebox{1pt}{$\bullet_{v}$}}{}{\searrow}\underset{\bblacksquare_b}{}\right\}
}
 \ar[d]^{
}
         \\
{	 
\left\{\underset{\bblacksquare_a}{}{\swarrow} \overset{\raisebox{1pt}{$ \bullet_{u=x=v}$}}{}{\searrow} 
\underset{\bblacksquare_b}{}\right\}
}
}\end{gathered}\end{equation}
Geometrically, this map is a ``finite model'' of the barycentric subdivision of the interval:
take the ``cell decomposition'' $[0,1]:=\{0\}\cup (0,1)\cup \{1\}$ and subdivide it as two cells:
$[0,1]:=\{0\}\cup (0\frac12)\cup \{\frac12\}\cup (\frac12,1)\cup \{1\}$. The ``finite model'' 
is obtained by contracting each open cell (i.e.~each open interval) to a point:
\newline\noindent\includegraphicss[0.5]{MtoL_with_boxes.png}
A verification shows that this map satisfies the lifting property $\Bbb S^n\lra \Bbb D^{n+1} \,\rtt\, \MtoL $, $n\geqslant 0$, defining trivial Serre fibrations.


\subsection{A summary of our notation} 
The following table is a summary of our notation. The line before last is how we would write our notation in computer code.
$$
  \begin{array}{ccccccc}
	    \{a,b\}
	         &\longrightarrow&
		   \{a{\searrow}b\}
		        &\longrightarrow&
			  \{a\leftrightarrow b\}
			      &\longrightarrow&
			        \{a=b\}
\\				 
	    \{a,b\}
	         &\longrightarrow&
		    \left\{{\ensuremath{a}}\rotatebox{-12}{\ensuremath{\to}}  \raisebox{-2pt}{\ensuremath{b}}\right\}
		        &\longrightarrow&
			  \{a\leftrightarrow b\}
			      &\longrightarrow&
			        \{a=b\}

\\				 
	    \{\bullet_a,\bullet_b\}
	         &\longrightarrow&
		    \left\{{\ensuremath{\bullet_a}}\rotatebox{-12}{\ensuremath{\to}}  \raisebox{-2pt}{\ensuremath{\bblacksquare_b}}\right\}
		        &\longrightarrow&
			  \{\star_a\leftrightarrow \star_b\}
			      &\longrightarrow&
				\{\bullet_{a=b}\}

\\ 
	    \verb|{a,b}|
	         &\longrightarrow&
		   \verb|{a->b}|
		        &\longrightarrow&
			  \verb|{a<->b}|
			      &\longrightarrow&
			        \verb|{a=b}|

				 \\
				    \text{(discrete space)}
				         &\longrightarrow&
					   \text{(Sierpinski space)}
					       &\longrightarrow&
					         \text{(antidiscrete space)}
						     &\longrightarrow&
						       \text{(single point)}
						         \end{array}
							 $$

\def\atob{\ensuremath{\left\{{\ensuremath{a}}\rotatebox{-12}{\ensuremath{\to}}  \raisebox{-2pt}{\ensuremath{b}}\right\}}}

\def\OtoC{\ensuremath{\left\{{\ensuremath{o}}\rotatebox{-12}{\ensuremath{\to}}  \raisebox{-2pt}{\ensuremath{c}}\right\}}}

\def\BtoSquare{\ensuremath{\left\{{\ensuremath{\bullet}}\rotatebox{-12}{\ensuremath{\to}}  \raisebox{-2pt}{\ensuremath{\bblacksquare}}\right\}}}

\def\LLambda{\ensuremath{
\left\{\underset{\bblacksquare}{}{\swarrow} \overset{\raisebox{1pt}{$ \bullet$}}{}{\searrow} 
\underset{\bblacksquare}{}\right\}}}

\def\LLambdaauxvb{\ensuremath{
\left\{\underset{\bblacksquare_a}{}{\swarrow} \overset{\raisebox{1pt}{$ \bullet_{u=x=v}$}}{}{\searrow} 
\underset{\bblacksquare_b}{}\right\}}}

\def\MMM{\ensuremath{
\left\{\underset{\bblacksquare}{}{\swarrow} \overset{\raisebox{1pt}{$ \bullet$}}{}{\searrow} 
\underset{\bblacksquare}{}
{\swarrow}\overset{\raisebox{1pt}{$\bullet$}}{}{\searrow}\underset{\bblacksquare}{}\right\}
}}

\def\MMMauxvb{\ensuremath{
\left\{\underset{\bblacksquare_a}{}{\swarrow} \overset{\raisebox{1pt}{$ \bullet_{u}$}}{}{\searrow} 
\underset{\bblacksquare_{x}}{}
{\swarrow}\overset{\raisebox{1pt}{$\bullet_{v}$}}{}{\searrow}\underset{\bblacksquare_b}{}\right\}
}}

\def\MMMABUVX{\ensuremath{
\left\{\underset{\bblacksquare_A}{}{\swarrow} \overset{\raisebox{1pt}{$ \bullet_{U}$}}{}{\searrow} 
\underset{\bblacksquare_{X}}{}
{\swarrow}\overset{\raisebox{1pt}{$\bullet_{V}$}}{}{\searrow}\underset{\bblacksquare_A}{}\right\}
}}

\section{Conjectures\label{sec_conjectures}} 
We formulate a few conjectures and open questions. 
At a cost of some repetion, we aim the exposition here to be self-contained as much as possible.

\subsection{Randomly generating a definition of compactness and contractibility}\label{sec_AI}
Our observations lead to a notation for topological properties of maps or spaces
so concise that the definitions of compactness, contractibility, and connectedness fit into two or four bytes. 
This notation makes explicit finite preorders implicit in these notions. Let us explain.

A concise and in some way intuitive notation for properties (=classes) of continuous maps 
is provided by iterated Quillen negations/orthogonals of maps of finite spaces:
{\em a word in two letters $l$,$r$, and a set of maps of finite topological spaces} 
represents a property (=class) of continuous maps. 
A property of maps applied to $X\lra\{\bullet\}$ or $\emptyset\lra X$ becomes a property of spaces,
and in way this notation is also able  define properties of spaces. 
Above we saw that to define connectedness, contractibility, and compactness, it is enough to consider one or two maps 
of finite spaces of size $\leqslant 5$ and $\leqslant 3$.  
\cite{LP1} gives a list of $20$ topological properties defined in this way
using a single map of spaces with $\leqslant 4$ and $\leqslant 3$ points, and 
\cite{MR} lists some $10$ natural properties  defined 
starting with the single map $\emptyset\lra\ptt$ using up to $7$ Quillen negations.
In particular, connectedness can be defined using a single map of spaces with two points 
(Proposition~\ref{propConnDisconn}), 
and compactness needs a single map of spaces with $4$ and $2$ points (Eq.~\eqref{KtoOlrSimplified}),
and contractibility needs a map of spaces with $5$ and $3$ points and a map from a space with $3$ points
to a singleton (Theorem~\ref{coro:contrCW}). 
A rough estimate of  the number of maps of preorders 
suggests that the definitions of these notions fit into $2$ bytes, or perhaps $4$.\footnote
{Let us bound the number $\#Maps_{p\to q}$  of maps from a preorder with $p$ elements to a preorder with $q$ elements
in terms of the number of labelled preorders with $p$ elements. 

Pick a preorder $Q$ with $q$ elements labelled by $1,..q$
and 
a preorder $P$ labelled by $1,..,p$. A subset (=increasing sequence)  $1\leq i_1\leq ... \leq i_{q'}\leq p$ with $q'\leqslant q-1$  
elements determines a (possibly not monotone) map $P\lra Q$ from $P$ into $Q$. 
Each map of unlabelled preorders with $p$ and $q$ elements can be constructed in this way.
Hence, the number of maps from a preorder with $p$ elements to a preorder with $q$ elements at is most 
the product of the number of labelled preorders with $p$ elements, the number of partitions of $p$ into $\leqslant q$ intervals,
and the number of (in fact, unlabelled) preorders with $q$ elements. 
Using the OEIS library (sequences 
\href{https://oeis.org/A001930}{A$001930$} and 
and \href{https://oeis.org/A000798}{A$000798$}) 
for $p=4$ and $q=2$ we get $\leqslant 355*4*3=4260\leqslant 2^{13}$, and 
for $p=5$ and $q=3$ we get $\leqslant 6942* 15*9=937170\approx 1000000\leqslant 2^{20}$
Thus, if we include the $lr$-suffix, the notions of connectedness and compactness fit into $2$ bytes, 
and contractibility may fit into $3$ or $4$ bytes.  
%
%
}

Fix a choice an encoding/notation as described above. This makes the following question precise.  
\begin{enonce}{Problem}[Ergo-logic of topology/AI]\label{prob_AI}
What is the probability that a random two or three bytes represent the notion of compactness, 
	contractibility, or connectedness (among ``nice'', e.g. Hausdorff, metrisable, etc spaces) ? 
	What is the probability that they represent a property of spaces or continuous maps 
	explicitly defined in say \cite[I,IX]{Bourbaki} ?
Is it non-negligible ? 
\end{enonce}
\begin{proof}[Evidence] By Corollary~\ref{coro:proper} and Corollary~\ref{coro:KtoOlr}  
	compactness is defined as left-then-right Quillen negation of any closed (=proper) map of finite topological spaces complicated enough;
	``complicated enough'' here means that the map has as retracts the maps in Eq.~\ref{KtoOlrSimplified}. Similarly, 
	contractibility is defined by any trivial Serre fibration complicated enough in a similar sense using the maps mentioned in Theorem~\ref{thm:Thm2}. 
Presumably the proportion of such maps is non-negligible. 
\end{proof}

Note that our notation makes it trivial to write a program which generates 
a complete definition of compactness or contractibility with non-negligible probabibily. 
We discuss this briefly in 
\cite[\S4.3, Question 4.12]{mintsGE}.

\begin{enonce}{Problem}[Ergo-learner/AI program for topology]\label{learner_AI}
Write a short program (kilobytes) following the 16 rules of ergo-learner
\cite[II\S25,p.168]{Ergobrain} 
which extracts maps of finite preorders when interacting with introductory texts on topology,
and  generates a list of maps of finite preorders corresponding to notions
defined in a typical introductory course of topology.
\end{enonce}

\subsection{Formalisation of topology in terms of finite preorders}

One wonders if this concise notation can be of use in formalisation of mathematics or teaching, 
in particular in an approach of \cite{Gowers}. In \cite[2.3]{mintsGE} we discuss 
how to rewrite the axioms of topology in rules for computation 
with labelled diagrams of finite preorders. We note that this point of view makes it clear 
that the axioms of topology can be seen as inference rules allowing to permute quantifiers 
$\forall\exists = \exists\forall$ in 
particular situations.

\begin{enonce}{Problem} Develop a proof system and a computer algebra system for topology in terms of Quillen negation/orthogonals and maps of finite topological spaces.
In particular, define a notion of ``diagram chasing'' with finite preorders which captures standard topological arguments involving contractibility and compactness. 
	In other words, define a (pseudo?-) category of ``formal'' topological spaces. 
\end{enonce}

Note that a finite preorder is a category, although a rather degenerate one (namely, such that there is at most one morphism between from one object to another).
Thus, we are looking for a definition of diagram chasing in a (pseudo?) category of categories. 

\begin{enonce}{Problem} Develop an introductory course using the notation in terms of finite preorders to introduce basics of topology and category theory at the same time.
\end{enonce}

The following is a more concrete problem along these lines. 

\begin{enonce}{Problem}[Simplicial approximation] Understand/interpret the simplicial approximation theorem as a diagram chasing rule. 
In particular, define a general diagram chasing rule which can be used for the infinite induction in our diagram chasing rendering of 
	Theorem $2$ of Urysohn \cite[IX\S$4.2$]{Bourbaki}, see Theorem~\ref{thm:Thm2}. 
\end{enonce}

\def\mathcalC{\{\emptyset\lra\ptt\}}
\def\mathcalCnobrackets{\emptyset\lra\ptt}
The identity in the following problem follows from the calculation of the orbit of $\mathcalC$ in \cite{MR}.
\begin{enonce}{Problem} Develop a computer algebra system able to express and prove that 
	$$\mathcalC^{rl}=\mathcalC^{rrrrll}=\mathcalC^{rllrrrll}$$
\end{enonce}

A more explicit problem along these lines is as follows, cf.~\cite[Question 4.13-14]{mintsGE}.

\begin{enonce}{Problem} Write down in terms of diagram chasing with preorders 
the standard argument showing that a compact Hausdorff space is necessarily normal. 
 
 Find a convenient  notation in which this argument becomes an intuitive calculation. 
\end{enonce}

\subsection{Interpreting our expressions for topological properties in other categories} 
In \S\ref{not:mintsGE} we define a notation for topological properties. It is tempting 
to intepret our expressions in other categories, for example so that the expression
for compactness and contractibility defines the same in categories of algebraic-geometric nature,
e.g.~the category of schemes. 

Note that our notation is based on maps of finite preorders
which can equivalently be seen as 
functors between categories, albeit rather degenerate ones. 
This allows to interpret our notation in the category of categories.

Our expressions can also be interepreted in the category of simplicial sets:
a preorder $P^\leqslant$ can also be viewed as a representable simplicial set 
$n^\leqslant \longmapsto \Hom_{\text{preorders}} (n^\leqslant, P^\leqslant)$.

\begin{enonce}{Problem} Calculate in the category of categories or in the category of simplicial sets
our expressions for $\pi_0$, compactness, and contractibility 
(i.e. the expressions in Proposition~\ref{prop9}, Corollary \ref{coro:KtoOlr} and \ref{coro:contrCW}).
Are the answers at all meaningful ? 
\end{enonce}

Perhaps  a more interesting problem is to interpret our notation in a category of algebraic nature,
such as that of modules, groups, or schemes. 
\cite[\S7]{DMG} suggests how to interpret the lifting property Eq.~\eqref{DefInjR} 
defining injectivity,
as the definition of monomorphism. 

\begin{enonce}{Problem} Find a meaningful interpretation of our expressions 
for $\pi_0$, compactness, and contractibility, in some other category,
e.g.~the category of schemes, modules, or finite groups of a fixed exponent.
\end{enonce}

\subsection{Iterated Quillen negation/orthogonals}

It appears that iterated Quillen negations/orthogonals have not been studied, e.g. we were unable to find 
a published calculation of a class of form $C^{ll}$ or $C^{rr}$ in any category. 
The following conjecture represents one of the first questions one may ask. 

\begin{conj} For each map $f$ of finite topological spaces there are only finitely many distinct 
	orthogonal classes of form $\{f\}^s$, $s\in \{l,r\}^n$, $n>0$.
\end{conj} 

In 
\cite{MR}  we 
verify this conjecture for $f$ being the simplest possible map $\emptyset\lra\ptt$,
and observe that most of the classes of this form 
are defined by properties introduced in an introductory course of topology 
--- surjective (l), injective (lrrr), subspace (rr), 
discrete (rl), having a section (lrr=lrrrllr), quotient (lrrrl), co-quotient (i.e. a surjective map with topology on the source induced from the target) (rrr), 
connected (rll), surjective on $\pi_0$ (rll=lrrrlll),  disjoint union (lrrrll), and a few others. 

\begin{conj}[Axiom M2]  For any finite set $P$ of maps of finite topological spaces and any string  $s\in \{l,r\}^n$, $n\geqslant 0$, 
	\label{conjM2} each morphism $f:X\lra Y$ decomposes both as 
	$$ X \xra{(P)^{sl}}\cdot \xra{(P)^{slr}} Y \text{ and as }  X \xra{(P)^{srl}}\cdot \xra{(P)^{sr}} Y ,$$
i.e.~$f$ factors both as
	$f=f_l\circ f_{lr}$ and $f=f_{rl}\circ f_{r}$  
	where $f_l\in P^{sl}$, $f_{lr} \in P^{slr}$, $f_{rl}\in P^{srl}$, $f_{r} \in P^{sr}$.
\end{conj}
\begin{proof}[Evidence] A verification shows this holds in a few easy examples, some of which correspond to basic facts in topology 
	explicitly stated in \cite{Bourbaki}. 
	 \cite[I\S$3.5$(Canonical decomposition of a continuous mapping)]{Bourbaki}
considers the canonical decomposition $f:X\lra Y$ as 
	$ f: X \xra{\phi} X/\approx_{f(x)=f(y)} \xra g f(X) \xra \psi Y$
	where $\phi:X\lra X/\approx_{f(x)=f(y)} $ is the canonical (surjective) mapping of X onto the quotient
	space $X/\approx_{f(x)=f(y)} $ by the equivalence relation $f(x)=f(y)$, 
	$\psi: f(X)\lra X$  is the canonical injection of the subspace $f(X)$ into $Y$,
	 and $g:X/\approx_{f(x)=f(y)} \lra f(X)$  is the bijection associated with $f$ (Set Theory, 
	  $\S5$, no. $3$). 
	By Proposition~\ref{PropQuotientl} $ \phi \in \left\{\, \{a\} \lra \{a\llrra b \} \, ,\,  \OcTooc \,\right\}^l $, 
	and the universal property of quotient spaces implies that $ \left\{\, \{a\} \lra \{a\llrra b \} \, ,\,  \OcTooc \,\right\}^{lr} $
	is the class of injective maps and thus contains $ X/\approx_{f(x)=f(y)} \xra g f(X) \xra \psi Y$.
	A verification shows that $\bigl\{\emptyset\lra\ptt\bigr\}^l=\bigl\{\{\star\}\lra\{\star\leftrightarrow\star\}\bigr\}^r\ni\phi\circ g$ is the class of surjections,
	$\bigl\{\emptyset\lra\ptt\bigr\}^{lr}=\bigl\{\{\star\}\lra\{\star\leftrightarrow\star\}\bigr\}^{rr}\ni\psi$ is the class of subsets.

	Proposition~\ref{propConnDisconn}\ref{propConnDisconnM2llr} shows this for $P:=\bigl\{ \{a,b\}\lra\ptt\bigr\}$.
	Proposition~\ref{prop9} shows that \cite[I\S$11.5$,Proposition~9]{Bourbaki} proves that each map $X\lra\ptt$ decomposes as 
	 $X \xra {(P_1)^l} X_{lr}
	          \xra {(P_1)^{lr}} \{\bullet\} $ 
		          for the class $P_1$ consisting of the following 4 morphisms of finite topological spaces:
			          $$ \{a,b\}\lra\{\bullet\}\,,\,  \{a\} \lra \{a\llrra b \} \, ,\,  \OcTooc \, 
				  ,       \left\{
					                   \raisebox{-2pt}{\ensuremath{\bblacksquare_a}\!\! \rotatebox{12}{\ensuremath{\leftarrow}}}\,
							                            {\ensuremath{\bullet_u}}\rotatebox{-12}{\ensuremath{\to}}  \raisebox{-2pt}{\ensuremath{\bblacksquare_b}}\right\}
									                     \lra
											                      \left\{{\ensuremath{\bullet_u}}\rotatebox{-12}{\ensuremath{\to}}  \raisebox{-2pt}{\ensuremath{\bblacksquare_{a=b}}}\right\}$$
The decompositions $A \lra A \sqcup Y \lra Y $
and $A \lra \operatorname{cl}_Y \Imm A \lra Y$ 
are of form 
$(\emptyset\lra\{\bullet\})^{lrrl}(\emptyset\lra\{\bullet\})^{lrr}$ and 
 $(\emptyset\lra\{\bullet\})^{rllrrll}(\emptyset\lra\{\bullet\})^{rllrrllr}$,
 resp., 
by \cite{MR}[Theorems~\ref{clrrl} and~\ref{clrr}, and Theorems~\ref{crllrrll} and~\ref{crllrrllr}], resp. 
\end{proof}

\subsection{Compactness}

Corollary~\ref{coro:proper} suggests the following conjecture. A calculation shows that $\Bigl\{ \{o\}\lra \left\{{\ensuremath{o}}\rotatebox{-12}{\ensuremath{\to}}  \raisebox{-2pt}{c}\right\} \Bigr\}^r_{<5}$ is the class of proper maps between spaces of size $<5$. Here by $P_{<5}$ we denote the subclass of $P$ consisting of maps of spaces with $<5$ points. 
\begin{conj} $\Bigl(\Bigl\{ \{o\}\lra \left\{{\ensuremath{o}}\rotatebox{-12}{\ensuremath{\to}}  \raisebox{-2pt}{c}\right\} \Bigr\}^r_{<5}\Bigr)^{lr}$ 
	is the class of proper maps.
\end{conj}

\begin{conj}\label{conj_proper}  Let $f$ denote a map of finite topological spaces. 
The class $\{f\}^{lr}$ is the class of all proper maps iff
$f$ is closed, not surjective, and satisfies the following:
\begin{enumerate}
\item $f$ has a fibre with two points bounded above but not below
\item $f$ has a fibre with two distinct topologically indistinguishable points
\item the image of $f$ is not both open and closed
\end{enumerate}
\end{conj}

\begin{conj} Being proper is maximal among $lr$-definable properties,
i.e.  for each class $Q$ it holds
$$(\text{proper maps}) \subsetneq Q^{lr} \implies Q^{lr}=(\text{all morphisms})$$
\end{conj}
\begin{proof}[Evidence] This can be shown for $Q$ consisting of maps of finite spaces.
In fact, for a class $Q$ of finite spaces the following implication holds:
$$(\text{proper maps of finite spaces}) \subsetneq Q^{lr} \implies Q^{lr}=(\text{all morphisms})$$
\end{proof}

The covering $\Bbb R=\bigcup\limits_{n\in \Bbb N} (-n,n)$ by open intervals $(-n,n)_{n\in \Bbb N}$
is a standard witness that the real line $\Bbb R$ is not compact. The corresponding lifting property gives
a reformulation of compactness often used in analysis: 

Each continuous real-valued function on a connected space $K$ is bounded 
iff 
$$ \emptyset \lra K \,\rtt\, \bigsqcup\limits_{n\in \Bbb N}\, (-n,n) \lra \Bbb R $$
This suggests the following problem.
\begin{enonce}{Problem} Find a map $f$ of finite topological spaces and a word $s\in \{l,r\}^{<\omega}$
such that 
a connected space $K$ is compact iff 
$${{\displaystyle \emptyset}\atop{\downarrow \atop \displaystyle  K}}
\in 
\{f\}^s
$$
\end{enonce} 

\subsection{Algebraic topology}
Recall $\MtoL$ denotes a particular map of finite topological spaces which is a trivial Serre fibration (and which is implicit in the definition of a normal space).
Perhaps one may want to replace it by some other trivial Serre fibration, such an the one implicit in the definition of 
a hereditary (=completely)  normal space \cite[IX\S4,Exercise 3b,p.$239$]{Bourbaki}
which happens to be the non-Hausdorff cone\footnote{In terms of preorders, this is the map obtained by adding a new least element both to the domain and the codomain. 
Topologically this is what you get by taking the usual cone of the domain and codomain, and then 
in each cone contracting each side to the base, i.e. $X\times [0,1)$ to $X$, in the cone 
$X\times [0,1] \,/\, X\times\{1\}$ of $X$.} 
over $\MtoL$. 

\subsubsection{Model structure}\label{sec_model_str}
We ask whether it is possible to define a model structure 
on the category of topological spaces
entirely 
in terms of maps of finite topological spaces. 

\begin{conj}[Model structure] A cellular map of finite CW complexes  is a trivial fibration iff it lies in $\bigl\{\MtoLLtoO\bigr\}^{lr}$.
	
	Any locally trivial fibre bundle with contractible  fibre over a paracompact space, and, more generally any numerable fibre bundle with contractible fibre, 
	and whose total space is a separable metric space, lies in  $\bigl\{\MtoLLtoO\bigr\}^{lr}$. 

	There is a model structure on $\Topp$ such that  $\bigl\{\MtoLLtoO\bigr\}^{lr}$ is its class of trivial fibrations. 
\end{conj}
\begin{proof}[Evidence] See \cite{V} for detailed speculations how define a model structure on $\Topp$. 
One precise question is whether there is  a model structure on $\Topp$ ``generated'' by maps of finite spaces
which are Serre (trivial) fibrations, i.e. the classes of (trivial) fibrations are 
$(C_{\operatorname{fin}})^{lr}$ and $(WC_{\operatorname{fin}})^{lr}$ 
where $(C_{\operatorname{fin}})$ and $(WC_{\operatorname{fin}})$ are the classes of maps of finite spaces which 
satisfy the lifting property definiting Serre fibrations, resp.~trivial Serre fibrations.  
\end{proof}  

\subsubsection{Brouwer fixed point theorem} 
We are not aware of a generalisation of the Brouwer fixed point theorem to non-Hausdorff spaces.  
	\begin{conj}[Brouwer fixed point theorem]\label{thm:BrowerFP}
	Let $K \xra { \left\{\MtoLLtoO\right\}^{lr}} \ptt$ and 
	$K \xra { \left(\left\{ \{o\}\lra \left\{{\ensuremath{o}}\rotatebox{-12}{\ensuremath{\to}}  \raisebox{-2pt}{c}\right\} \right\}^r_{<5}\right)^{lr} } \ptt$.
	Then each endomorphism of $K$ has a fixed point.
\end{conj}
\begin{proof}[Evidence] The conjecture holds for separable metric spaces. Indeed, by Theorem~\ref{coro:contrCW} 
$K$ is a retract of $\Bbb R^\omega$. A fixed point of its  endomorphism $f:K\lra K$ is a fixed point of
the endomorphism $\Bbb R^\omega\lra K \xra f K \lra \Bbb R^\omega$ of $\Bbb R^\omega$.  
\end{proof}
\begin{rema} The following diagrams represent the statement of the conjecture above.
	We wonder if they are appropriate for any system of formalised mathematics. 
\begin{equation}\begin{gathered}\label{def:Brouwer:fp}
	\vcenter{\hbox{\xymatrix
	{       & \{\bullet\} \ar@{-->}[ld]|x \ar@{-->}[rd]|x  & \\ 
	K \ar[rr]  \ar[dr]_(0.42){ 
	\left(\left\{ \{o\} \longrightarrow \left\{ o
	                           \rotatebox{-12}{\ensuremath{\to}}\raisebox{-2pt}{\ensuremath{c}}
											 \right\}\right\}^r_{<5}\right)^{ lr} 
	}  
	& & K \ar[dl]^(0.42){\{\MtoLLtoO\}^{ lr}} 
	\\  & \{\bullet\} &  } }}
	\ \ \ \ 
	\vcenter{\hbox{\xymatrix
	{       \{\bullet\} \ar@{-->}[d]   \\ 
	K \ar@(ur,dr)
	\ar[d]_(0.42){ 
	\left\{ \{o\} \longrightarrow \left\{ o
				   \rotatebox{-12}{\ensuremath{\to}}\raisebox{-2pt}{\ensuremath{$c$}  } 	\right\} 
				   \right\}^{ lr} 
	}^  
	{\{\MtoLLtoO\}^{ lr}} 
	\\   \{\bullet\}   } }}
\end{gathered}\end{equation}
\end{rema}

\subsubsection{Lebesgue covering dimension} 
	Let $\partial\Delta^n$ be a finite topological space weakly homotopy equivalent (``modelling'') the sphere $\Bbb S^n$
	such that 
	$\Bbb S^n \lra \partial\Delta^n$ is a trivial Serre fibration.  
\newcommand\HH{\mathrm{H}} \newcommand\cHH{\check{\mathrm{H}}}
	Let $\cHH_q(X),q\geqslant 0$ denote a \v Cech $q$-th homology group. We do not completely specify 
	the (co)homological theory because we are not sure which is 
	best appropriate for the conjecture below.
\begin{conj}[Lebesgue dimension] 
	For a finite CW complex $X$, \\ 
	$\cHH_q(X)=0$ for $q>n$ iff 
	 $$X\lra \ptt \,\in\, \bigl\{ \partial\Delta^n \lra\ptt, \MtoLLtoO \bigr\}^{lr}$$
\end{conj}
\begin{proof}[Evidence]
	For $n=0$ $\Bbb S^0=\{a,b\}$ and 
	$\partial\Delta^1:=\Bbb S^0=\{a,b\}$. An easy calculation shows that for a space $X$ with finitely many connected components, 
	$X\lra\ptt$ lies
	$\bigl\{ \{a,b\}\lra\ptt, \MtoLLtoO \bigr\}^{lr}$ iff for each connected component $X'$ of $X$ it holds
	$X'\lra\ptt$ lies
	$\bigl\{ \MtoLLtoO \bigr\}^{lr}$, i.e. is contractible by Theorem~\ref{coro:contrCW}. In other words, 
	such an $X$ is a disjoint union of contractible connected components, which for a ''nice'' enough space $X$ is equivalent to  $\cHH_q(X)=0$ for $q>0$, as required.

It is tempting to think that 
	$\Bbb S^n\lra\partial\Delta^n$ is a
	trivial Serre fibration and thereby is in 
	$ \bigl\{ \partial\Delta^n \lra\ptt, \MtoLLtoO \bigr\}^{lr} $,
	and thus so is $\Bbb S^n\lra \ptt$ and hence at least there are some examples of $X$ in the class 
	with $H_n(X)\neq 0$.

	Let $\dim X $ denote the {\em Lebesgue (covering) dimension}, see \cite[\S$3.1$,footnotes (19),(20)]{V}
	for suggestions how to define it in terms of the left lifting property. 
For a normal space $X$,
	$\dim X \leqslant n$ iff  $A\to X \rtt \mathbb{S}^n\to \ptt$
	for every closed subset $A\subset X$. 
	Theorem~\ref{thm:Thm2} replaces in Theorem $2$ of Urysohn \cite[IX\S$4.2$]{Bourbaki} 
	the assumption ``$A$ is a closed subset of a normal space $X$'' by 
	$A\lra X \,\rtt\, \MtoLLtoO$, thus it is tempting to do the same here.

	Also, $\dim X \leqslant n$ implies that $\cHH_q(X,A)=0$, $q>n$, for any closed subset $A$ of $X$, 
	and an argument ``gluing in holes in $X$ of low dimension'' should give that
	$X \xra { \bigl\{ \partial\Delta^n \lra\ptt, \MtoLLtoO \bigr\}^{lr} } \ptt$ implies $H_q(X)=0$ for $q>n$.
\end{proof}

\begin{rema} To get a reasonable definition of $\pi_0(X)$ it was not enough to consider only the morphism 
	$\{a,b\}\lra \ptt$, i.e.~$\Bbb S^0\lra\ptt$, but rather we needed to consider 
	more complicated morphisms (with $\leqslant 3$ points)
	implicit in \cite[I\S$11.5$,Proposition $9$]{Bourbaki}, see Proposition~\ref{prop9}.
        So perhaps also in the Problem above we are missing some non-obvious morphisms such as those  
	implicit in the definition of the Lebesgue covering dimension
	and order of coverings, see \cite[\S$3.1$,footnotes (19),(20)]{V}.
\end{rema}

\section{Appendix}

In \S\ref{AppTietzeLemmaNoOrthogonals} 
we rewrite the proof of Theorem~\ref{thm:Thm2}
above using explicit diagram chasing arguments
instead of properties of orthogonals (negations).
In \S\ref{App:Fibrations} we speculate on how to define fibrations
using the notion of a microfibration.
In \S\ref{App:mintsGEexamples} we use very simple examples to explain how to read off finite combinators
from the text of Bourbaki.

\subsection{A diagram-chasing rendering of the proof of the Tietze Lemma}
\label{TietzeLemmaNoOrthogonals}\label{AppTietzeLemmaNoOrthogonals}

We rewrite the proof of Theorem~\ref{thm:Thm2} above using explicit diagram chasing arguments
instead of properties of orthogonals (negations).
In this form it is easier to see its relation to the exposition in \cite[IX\S$4.2$,Theorem 2]{Bourbaki}. 

The exposition below is self-contained. 

\begin{enonce}{Theorem}[Urysohn]\label{thm:Thm2:diagrams} 
	$${{\displaystyle \barRR}\atop{\displaystyle \downarrow\atop\displaystyle \{\bullet\}}} \,\in\, \MtooLVtoolr $$
\end{enonce}
\begin{proof} 
	Let $A\lra X$ be an arbitrary map such that $A\lra X \,\rtt\, \MtoL$. We need to show that $A\lra X \,\rtt\, \barRR \lra \{\bullet\}$. 
Let $f:X\lra \barRR$ be a map. 

	Decompose the interval into smaller intervals as 
	$$[-1,1]=\{-1\}\cup (-1,a_1) \cup \{a_1\}\cup ... \cup \{a_n\}\cup (a_n,1)\cup \{1\}$$
	and subdivide an open interval $(a_m,a_{m+1})$ in two halves:
	$$[-1,1]=\{-1\}\cup (-1,a_1) \cup \{a_1\}\cup ... \cup \{a_m\} \cup (a_m, a'_m) \cup \{a'_m\}\cup (a'_m,a_{m+1})\cup \{a_m\} \cup ... \cup \{a_n\}\cup (a_n,1)\cup \{1\}$$
Contracting each open subinterval in these decompositions gives maps to finite spaces 
	$\barRR \lra \vvLambda_n$ and $\barRR\lra \vvLambda_{n+1}$.
	Note that $\vvLambda_1=\vvLambda$ and $\vvLambda_2=\MM$. 
Contracting closed subintervals $[-1,a_m]$ and $[a_{m+1},1]$ gives maps to finite spaces 
$\barRR \lra \vvLambda$ and $\barRR\lra \MM$. 
	These maps fit into a commutative diagram with a pull-back square on the right:\footnote{Here we construct the lifting by hand instead of using that orthogonals are closed under colimits and compositions as in the proof of Theorem~\ref{thm:Thm2}.}
	$$\xymatrix@C=2.39cm{
		[-1,1] \ar@/_1pc/[rd] \ar@/_1pc/[rdd] \ar@/^1pc/[rrd]  \ar@/^1.3pc/[rrdd] &&  \\ 
		& \vvLambda_{n+1} \ar[r]  \ar[d]^{\therefore\MtooLlr} & \MM \ar[d]^{\MtooLlr}  \\
	                               & \vvLambda_n  \ar[r]  &           \vvLambda 
	}$$
	Hence, given 
	a map $g_n:X\lra \vvLambda_n$ making commutative the square 
	in the diagram below
        $$\xymatrix@C=2.39cm{
		                A \ar[r]|f  \ar[d] & [-1,1] \ar[r] 
				& \vvLambda_{n+1} \ar[d] \\ 
				        X  \ar[rr]|{g_n} \ar@{-->}[urr]|{g_{n+1}} 
					&                              & \vvLambda_n \\
					        }$$
we construct a map $g_{n+1}:X\lra \vvLambda_{n+1}$ by first using the lifting property $A\lra X \,\rtt\, \MtoL$ and then the pullback property:
	$$\xymatrix@C=2.39cm{
		A \ar[r]|f  \ar[d] & [-1,1] \ar[r] \ar[rd] \ar@/^1pc/[rr]  \ar@/_0.3pc/[rrd] & \vvLambda_{n+1} \ar[r]  \ar[d] & \MM \ar[d] \\
	X  \ar@/_0pc/[rr]|{g_n} \ar@/_1pc/@{-->}[urr]|{g_{n+1}}  \ar@/_1pc/@{-->}[urrr] &                              & \vvLambda_n  \ar[r]  &           \vvLambda 
	}$$
Thus by induction we can construct a sequence of maps $g_n:X\lra \vvLambda_n$, $n\geqslant 1$, 
	$$\xymatrix@C=1.39cm{
		A \ar[r]|f  \ar[d] & [-1,1] \ar[r] & ... \ar[r] & \vvLambda_{n+1} \ar[r] & \vvLambda_n \ar[r] & ... \ar[r] & \vvLambda_2=\MM \ar[d] \\
	X  \ar@/_1pc/@{-->}[rrrru]|{g_n} \ar@/_1pc/@{-->}[urrr]|{g_{n+1}}  \ar@/_1pc/@{-->}[urrrrrr]|{g_2}   \ar[rrrrrr]|{g_1}     &&& &                              &       &           \vvLambda_1=\vvLambda 
	}$$
given a map $g_1:X\lra \vvLambda_1$ fitting into the diagram  
	\begin{equation}\begin{gathered}\label{gOne} \xymatrix@C=3.39cm{
		A \ar[r]|f  \ar[d] & [-1,1] \ar[r] & 
		\vvLambda_1=\vvLambda \ar[d] \\ 
	X  
		\ar@{-->}[urr]|{g_1} \ar[rr] & & \{\bullet\} 
	}\end{gathered}\end{equation}
The diagram is a lifting square, hence such a map exists whenever $A \xra f B \,\rtt\, \vvLambda\lra \ptt$. 
	Hence, we can construct such a sequence of maps $g_n:X\lra\vvLambda_n$ whenever $A \xra f B \,\in\,  \Bigl\{\MtoLLtoO\Bigr\}^{l}$.

Now, a standard argument using ``the idea of uniform convergence'' can be used to construct a continuous function $g:X\lra [-1,1]$ 
	by setting $g(x):= \bigcap_n \cl (\tau_n\inv(g_n(x)))$. Indeed, by construction we have $\tau_{n+1}\inv(g_{n+1}(x))\subseteq \tau_n\inv(g_n(x))$,
	and we may choose the maps $\tau_n:\barRR\lra \vvLambda_n$ such that the intervals $\tau_n\inv(g_n(x))$ become arbitrarily small. 
	To see that $g:X\lra\barRR$ is continuous, note that for each open neighbourhood $U_t\ni t\in [-1,1]$ there is $n>1$ and a finite 
	open subset  $ U_n\subset \vvLambda_n$ such that $t\in \tau_n\inv(U_n)$ and  $ \cl ( \tau_n\inv(U_n)) \subset U_t$; 
	hence by construction of $g$ 
	we have that  $g(g_n\inv(U_n))\subset U_n \subset U_t$, hence there is an open neighbourhood of $x$ which maps inside $U_t\ni t:=g(x)$.    
%
\end{proof}


\subsection{A definition of a fibration ?}\label{App:Fibrations}  
We try to rephrase the definition of a {\em microfibration} in terms of finite topological spaces and Quillen negation. 
More details can be found in 
\cite{V} which speculates how to define a model structure on the category of topological spaces in terms of finite topological spaces.

Recall that a map $f:X\lra Y$ is called a {\em microfibration} iff for any closed cofibration $A\lra B$ and any commutative square with sides $A\lra B$ and $f:X\lra Y$
 there is a lifting $h:U\lra X$ defined on an open neighbourhood $A\subset U\subset B$; often spaces involved are assumed to be ``nice'' in some sense. .
$$\xymatrix@C=3.39cm{ & A \ar[r]\ar[d]|{\text{(closed cofibration)}} & X\ar[d]^{\therefore\text{(microfibration)}}  
\\ U \ar@{-->}[r]|{\text{(open subset)}} \ar@/_1pc/@{-->}[urr]  &B \ar[r] & Y } $$
\begin{conj} A cellular map $f:X\lra B$ of finite CW complexes is a fibration iff it decomposes as
	$$X \xrightarrow { \left( \displaystyle {\left\{\bullet_o \rotatebox{-12}{\ensuremath{\to}}\raisebox{-2pt}{\ensuremath{\filledstar_a\leftrightarrow \filledstar_b}}\right\}}\atop
	{\downarrow\atop{ \displaystyle \left\{{\bullet_o\!=\!\filledstar_a\leftrightarrow \filledstar_b}\right) }} \right\}^l }
  X_{{}_p\!\!\searrow}\underset{B}{}
  \xrightarrow 
	{ \left( \displaystyle {\left\{\bullet_o \rotatebox{-12}{\ensuremath{\to}}\raisebox{-2pt}{\ensuremath{\filledstar_a\leftrightarrow \filledstar_b}}\right)}\atop
		{\downarrow\atop{ \displaystyle \left\{{\bullet_o\!=\!\filledstar_a\leftrightarrow \filledstar_b}\right\} }} \right)^{lr}, \MtooLVtoolrScriptsize  }
	B
	,$$
	i.e.~$f$ decomposes as $f=f_lf_{lr}$ where 
	$f_l\in \left\{ \left\{o \rotatebox{-12}{\ensuremath{\to}}\raisebox{-2pt}{\ensuremath{a\leftrightarrow b}}\right\}\longrightarrow
	  \left\{{o\!=\!a\leftrightarrow b}\right\} \right\}^l $ and 
	  $f_{lr}\,\in\, \left\{ \left\{o \rotatebox{-12}{\ensuremath{\to}}\raisebox{-2pt}{\ensuremath{a\leftrightarrow b}}\right\}\longrightarrow
	          \left\{{o\!=\!a\leftrightarrow b}\right\} \right\}^{lr} \cap \left\{\MtoLLtoO\right\}^{lr}  $.

	A locally constant map over a paracompact space, and, more generally, a
numerable fibre bundle\footnote{
	Recall that a map $f:X\lra B$  
	is called {\em a numerable fibre bundle} iff there exists a numerable
	covering $\{V_i \}_{x\in \omega}$ of the {\em  base} $B$ 
	such that its restriction $f_{ V_i}$ (the part of $f$
	over $V_x$) is trivial for every $i\in \omega$ \cite[Def.7.1]{Dold}. 
	\cite[p.$460$]{DoldICM} says that ``The applications to bundle theory are based on the notion of numerable
	bundle. ... The class  of numerable bundles seems to
	have the right size for a satisfactory bundle theory...''. 
	} of separable metric spaces  always decomposes in this way. 
\end{conj}
\begin{proof}[Evidence] A verification shows that 
$\left\{ \left\{o \rotatebox{-12}{\ensuremath{\to}}\raisebox{-2pt}{\ensuremath{a\leftrightarrow b}}\right\}\longrightarrow
	          \left\{{o\!=\!a\leftrightarrow b}\right\} \right\}^l $
		  is the class of open maps $A\lra B$ where
		                    topology on $A$ is induced from $B$. In particular, an injective map in this class 
				    is ``represented'' by an open subset. 

	Define the {\em non-Hausdorff mapping cylinder} of a map $f:X\lra B$, denoted by $ X_{{}_f\!\!\searrow}\underset{B}{}$, to be 
  the disjoint union $X \sqcup B$ equipped with the following topology:
  an open subset is either an open subset of $X$,
  or the union of an open subset of $B$ and its preimage in $X$. A motivation for the terminology is that 
	this is the quotient of 
	the (usual) mapping cylinder $X\times [0,1] \sqcup B / \{(x,1)\approx f(x)\}_{x\in X} $ by $X\times [0,1)$ for compact Hausdorff $X$ and $B$. 
	The non-Hausdorff mapping cylinder fits into a sequence
  $$X \xrightarrow { \left\{ \left\{o \rotatebox{-12}{\ensuremath{\to}}\raisebox{-2pt}{\ensuremath{a\leftrightarrow b}}\right\}\longrightarrow
  \left\{{o\!=\!a\leftrightarrow b}\right\} \right\}^l }
  X_{{}_f\!\!\searrow}\underset{B}{}
  \xrightarrow {\left\{ \left\{o \rotatebox{-12}{\ensuremath{\to}}\raisebox{-2pt}{\ensuremath{a\leftrightarrow b}}\right\}\longrightarrow
	\left\{{o\!=\!a\leftrightarrow b}\right\} \right\}^{lr}   } B
	.$$
	The lifting property $A\lra B \,\rtt\, X_{{}_f\!\!\searrow}\underset{B}{} \lra B$ says precisely that 
	for a square commuting on an open subset  $U\subset A$ of $A$, one can find a lifting $h:V\lra X$ on an open subset $U\subset V\subset B$ of $B$.
	The latter property is similar to the defining property of a microfibration. Finally, one hopes that for nice enough maps 
	being a microfibration is equivalent to being a fibration, in the same way as being a Serre or Hurwitcz fibration are equivalent. 
\end{proof}

\subsection{Transcribing Bourbaki (dense and $T_0$) into diagram chasing}\label{App:mintsGEexamples}
The goal of this appendix taken from \cite{mintsGE} 
is explain on {\em very} simple examples how extract 
{\em diagram chasing and finite combinatorics}
 implicit in the {\em text} of basic topological definitions in \cite{Bourbaki}.

Here we transcribe line by line very carefully a couple of vary basic and simple definitions 
stated in \cite{Bourbaki}. We then explain how the resulting reformulations
remind us of category theory.

\%subsection{Dense subspaces and Kolmogoroff $T_0$ spaces.}
We shall now transcribe the definitions of {\em dense} and {\em Kolmogoroff $T_0$} spaces.

\subsubsection{``$A$ is a dense subset of $X$.''}
By definition \href{http://mishap.sdf.org/mints-lifting-property-as-negation/tmp/Bourbaki_General_Topology.djvu}{[Bourbaki, I\S$1.6$, Def.12]},
\newline\noindent\includegraphics[width=\linewidth]{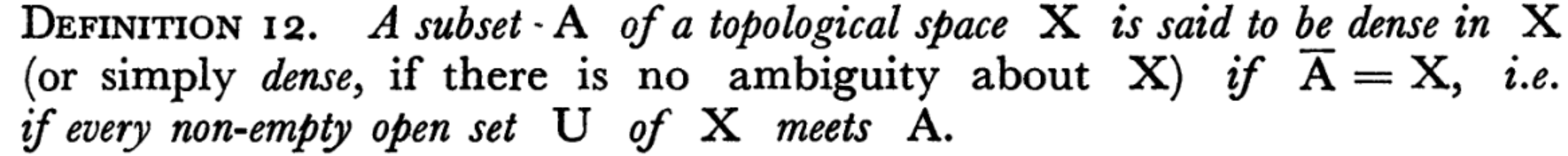}
 Let us transcribe this by means of the language of arrows.

 {\sf A subset $A$ of a topological space $X$} is an arrow $A\lra X$.
 (Note we are making a choice here: there is an alternative translation analogous to the one used in the next sentence).
 An {\sf open subset $U$ of $X$} 
 is an arrow $X\lra \{ U \ra U' \}$ ; here  $\{ U \ra U' \}$ denotes the topological space
 consisting of  one open point $U$ and one closed point $U'$; by the arrow $\ra$ we mean that
 that $U'\in cl(U)$.
  {\sf Non-empty}: a subset $U$ of $X$ is {\em empty} iff
  the arrow  $X\lra \{ U \ra U' \}$ factors as  $X\lra  \{ U' \}\lra \{ U \ra U' \}$ ;
  here the map $ \{ U' \}\lra \{ U \ra U' \}$ is the obvious map sending $U'$ to $U'$.
  {\sf set $U$ of $X$ meets $A$}: $U\cap A =\emptyset$ iff the arrow $A\lra X\lra \{ U \ra U' \}$ factors as
  $A\lra \{ U' \} \lra \{ U \ra U' \}$.

  Collecting above (Figure 1c), we see that
  a map $A\xra f X $ has dense image iff
  $$ A \xra f X \rtt  \{ U' \}\lra \{ U \ra U' \}$$

  Note a little miracle:
  $\{ U' \}\lra \{ U \ra U' \}$ is the simplest map whose image isn't dense.
  We'll see it happen again.

%
%
%

	 \subsubsection{Kolmogoroff spaces, axiom $T_0$.} By definition \href{http://mishap.sdf.org/mints-lifting-property-as-negation/tmp/Bourbaki_General_Topology.djvu}{[Bourbaki,I\S1, Ex.2b; p.$117$/122]},
  \newline\noindent\includegraphics[width=\linewidth]{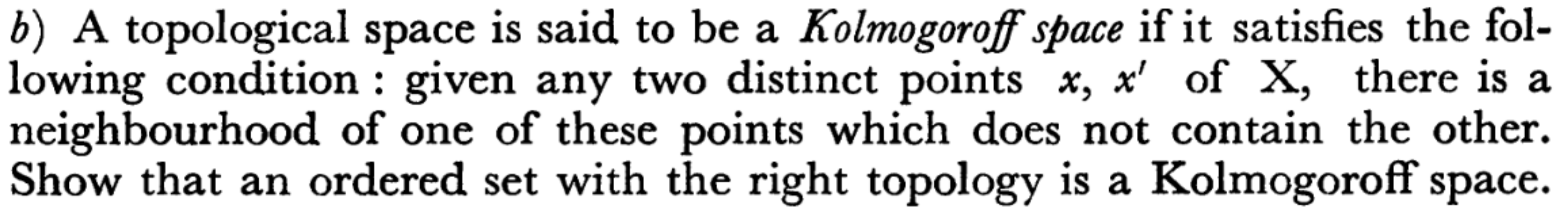}
  Let us transcribe this. {\sf given any two ...~points  $x$, $x'$ of $X$}: given a map $\{x,x'\}\xra f X$.
  {\sf two {\em distinct} points}:  the map $\{x,x'\}\xra f X$ does not factor through a single point,
  i.e.  $\{x,x'\}\lra  X$ does not factor as $\{x,x'\}\lra \{ x=x' \} \lra  X$.
  The negation of the sentence {\sf there is a neighbourhood which does not contain the other} defines a topology on the set $\{x,x'\}$: indeed,
   the antidiscrete topology on the set $\{x,x'\}$ is the only topology with 
   the property that
   {\sf there is [no]  neighbourhood of one of these points which does not contain the other}.
   Let us denote by  $\{x \llrra x'\}$ the antidiscrete space consisting of $x$ and $x'$.
   Now we note that the text implicitly defines
   the space $\{x \llrra x'\}$, and the only way to use it is to consider
   a map $\{x \llrra x'\}\xra f X$  instead of the map $\{x,x'\}\xra f X$.


   Collecting above (see Figure~1d), we see that {\em a topological space $X$ is said to be a {\em Kolmogoroff} space
   iff any map $\{x \llrra x'\}\xra f X$ factors as $\{x \llrra x'\}\lra \{ x=x' \} \lra  X$.}

   Note another little miracle: it also reduces to  orthogonality of morphisms
   $$ \{x \llrra x'\}\lra \{ x=x' \} \rtt X\lra \{ x=x' \} $$
   and $\{ x\llrra x' \}$ is the simplest non-Kolmogoroff space.

\subsubsection{Separation axiom $T_1$.} 

\includegraphics[test]{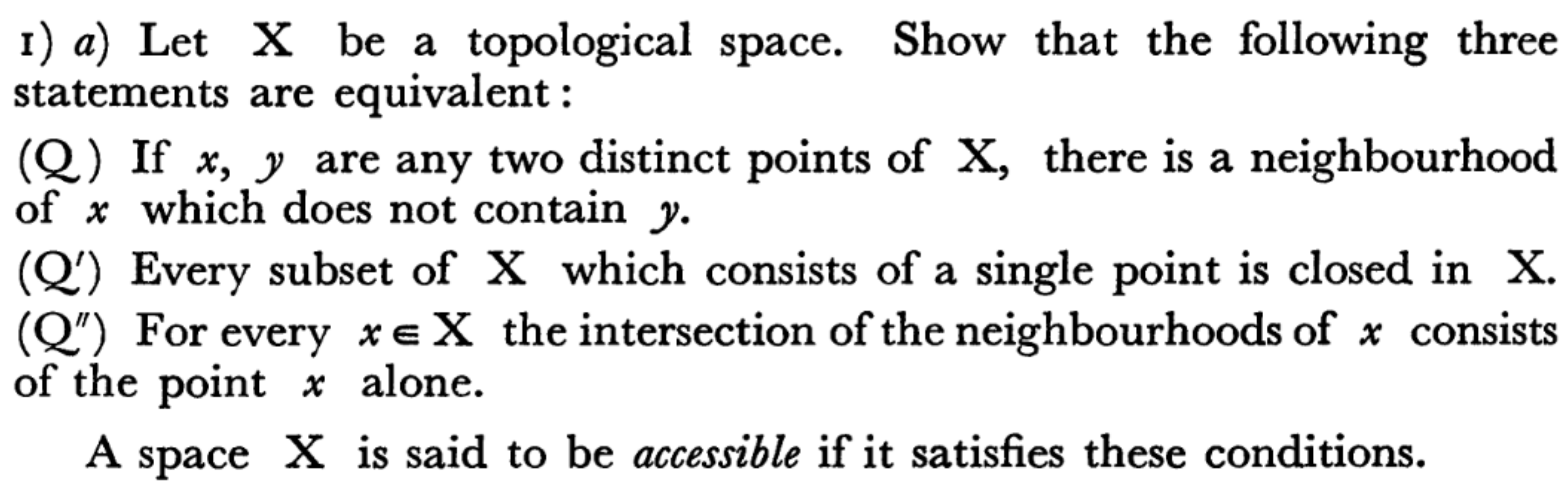}

{\tt If $x,y$ are any two \Remove{distinct} points of $X$}: take a map $\{x,y\}\lra X$. 
The point $x$ and $y$ are {\em distinct} iff the map $\{x,y\}\lra X$ does not factor as
$\{x,y\}\lra \{x=y\} \lra X$. {\tt there is a neighbourhood of $x$ which does not contain $y$}
is almost a definition of topology on $\{x,y\}$: there is a unique topology on $\{x,y\}$
such that, in this topology, there is a neighbourhood of $x$ which does not contain $y$
but not for $y$, i.e.~there is no neighbourhood of $y$ which does not contain $x$.
This topology is the topology with $x$ the open point, and $y$ the closed point.
Thus, we reformulate (Q) as:
\begin{quote} (Q$'''$)\,\,\, Each continuous mapping $\{x\ra y\}\lra X$ factors as 
	$\{x\ra y\}\lra \{x=y\}\lra X$. 
\end{quote}

Note yet another little miracle: it also reduces to  orthogonality of morphisms
   $$ \{x \ra y\}\lra \{ x=y \} \rtt X\lra \{ x=y \} $$
   and $\{ x\ra y \}$ is the simplest non-accessible space.

   \def\rrt#1#2#3#4#5#6{\xymatrix{ {#1} \ar[r]|{} \ar@{->}[d]|{#2} & {#4} \ar[d]|{#5} \\ {#3}  \ar[r] \ar@{-->}[ur]^{}& {#6} }}
   \begin{figure}
	   \begin{center}
		   \large
		   $ (a)\ \xymatrix{ A \ar[r]^{i} \ar@{->}[d]|f & X \ar[d]|g \\ B \ar[r]|-{j} \ar@{-->}[ur]|{{\tilde j}}& Y }$
		   $(b)\  \rrt  {\{\}}  {} {\{\bullet\}}  X {\therefore(surj)} Y $
		   $(c)\  \rrt  {A}  {\therefore(dense)} {B}  {\{U'\}} {} {\{U\ra U'\} } $
		   $(d)\  \rrt  {\{ x\llrra x' \}}  {} { \{ x=x' \}}  X {\therefore(T_0)} {\{ x=x' \}} $
	   \end{center}
	   \caption{\label{fig1}
	   Lifting properties. Dots $\therefore$ indicate free variables and what property of these variables is being defined;
	   in a diagram chasing calculation, "$\therefore(dense)$" reads as:
	   given a (valid) diagram, add label $(dense)$ to the corresponding arrow.\newline
	    (a) The definition of a lifting property $f\rtt g$: for each $i:A\lra X$ and $j:B\lra Y$
	    making the square commutative, i.e.~$f\circ j=i\circ g$, there is a diagonal arrow $\tilde j:B\lra X$ making the total diagram
	    $A\xra f B\xra {\tilde j} X\xra g Y, A\xra i X, B\xra j Y$ commutative, i.e.~$f\circ \tilde j=i$ and $\tilde j\circ g=j$.
	     (b) $X\lra Y$ is surjective 
	      (c) the image of $A\lra B$ is dense in $B$
	       (d) $X$ is Kolmogoroff/$T_0$
	       }
   \end{figure}

\subsubsection{Finite topological spaces as categories.}
Our notation  $\{ U' \}\lra \{ U \ra U' \}$ and $\{x \llrra x'\}\lra \{ x=x'
\}$ suggests that {\em we reformulated the two topological properties of being dense and Kolmogoroff
in terms of diagram chasing in (finite)  categories}. And indeed, we may think
of finite topological spaces as categories and of continuous maps between them as {\em functors},
as follows; see Appendix~\ref{app:top-notation} for details and a definition of
our notation for finite topological spaces and maps between them.

A {\em topological space} comes with a {\em specialisation preorder} on its points: for
points $x,y \in X$,  $x \leqslant y$ iff $y \in \operatorname{cl} x$ ($y$ is in the {\em topological closure} of $x$).
 The resulting {\em preordered set} may be regarded as a {\em category} whose
 {\em objects} are the points of ${X}$ and where there is a unique {\em morphism} $x{\ra}y$ iff $y \in cl x$.

 For a {\em finite topological space} $X$, the specialisation preorder or
 equivalently the corresponding category uniquely determines the space: a {\em
 subset} of ${X}$ is {\em closed} iff it is
 {\em downward closed}, or equivalently,
 it is a subcategory such that there are no morphisms going outside the subcategory.

 The monotone maps (i.e. {\em functors}) are the {\em continuous maps} for this topology.

 We denote a finite topological space by a list of the arrows (morphisms) in the
 corresponding category; '$\leftrightarrow $' denotes an {\em isomorphism} and
 '$=$' denotes the {\em identity morphism}.  An arrow between two such lists
 denotes a {\em continuous map} (a functor) which sends each point to the
 correspondingly labelled point, but possibly turning some morphisms into
 identity morphisms, thus gluing some points.



\end{document}